\documentclass[11pt]{amsart}

\usepackage{epsfig}
\usepackage{amsmath}
\usepackage{amssymb}
\usepackage{amscd}
\usepackage{latexsym}
\usepackage{mathtools}
\usepackage{tabularx}
\usepackage{blkarray}
\usepackage{a4wide}
\usepackage[usenames]{color}
\usepackage{mathdots}
\usepackage{mleftright}
\usepackage{subfigure}
\usepackage{url}
\usepackage{tikz,tikz-cd}
\usetikzlibrary{decorations.markings}
\usetikzlibrary{calc}
\usepackage[matrix,arrow,tips,curve]{xy}
\usepackage{pb-diagram,pb-xy}
\usepackage{verbatim}
\usepackage{mathrsfs}
\usepackage{color}
\usepackage{scalerel}
\usepackage{mathabx}
\usepackage{enumitem}

\setlist[enumerate]{label=\textnormal{(\arabic*)}}

\definecolor{cadmiumgreen}{rgb}{0.0, 0.42, 0.24}
\usepackage[
colorlinks, citecolor=cadmiumgreen,
pdfauthor={Omid Amini, Lucas Gierczak}, 
pdftitle={Combinatorial linear series},
pdfstartview ={FitV},
]{hyperref}

\usepackage{todonotes} 

\newcommand{\out}{\ssub{\mathrm{out}}!}
\newcommand{\outval}{\ssub{\mathrm{outval}}!}

\newcommand{\Trop}{\mathrm{trop}} 
\newcommand{\trop}{{^{\hspace{-.02cm}\scaleto{\Trop}{4.4pt}}}} 
\newtheorem{thm}{Theorem}[section]

\newtheorem{lemma}[thm]{Lemma}

\newtheorem{propdefi}[thm]{Proposition - Definition}
\newtheorem{prop}[thm]{Proposition}

\newtheorem{claim}[thm]{Claim}
\newtheorem{cor}[thm]{Corollary}

\newtheorem{question}[thm]{Question}
\newtheorem*{question*}{Question}

\theoremstyle{definition}
\newtheorem{defii}[thm]{Definition}
\newenvironment{defi}
	{\pushQED{\qed}\defii}
	{\popQED\enddefii}
	\newtheorem{convention}[thm]{Convention}
\newenvironment{conv}
	{\pushQED{\qed}\convention}
	{\popQED\endconvention}
\newtheorem{remm}[thm]{Remark}
\newenvironment{remark}
	{\pushQED{\qed}\remm}
	{\popQED\endremm}
\newtheorem{exx}[thm]{Example}
\newenvironment{example}
	{\pushQED{\qed}\exx}
	{\popQED\endexx}
\newtheorem{factt}[thm]{Fact}
\newenvironment{fact}
	{\pushQED{\qed}\factt}
	{\popQED\endfactt}

\newcommand{\Red}{\mathrm{Red}}

\newcommand{\Q}{\mathbb{Q}}

\newcommand{\Z}{\ssub{\mathbb{Z}}!}

\newcommand{\R}{\mathbb{R}}

\newcommand{\N}{\mathbb{N}}

\newcommand{\ord}{\ssub{\mathrm{ord}}!}
\newcommand{\slope}{\ssub{\mathrm{slope}}!}
\newcommand{\Rat}{\ssub{\mathrm{Rat}}!}

\newcommand{\Div}{\operatorname{Div}}

\newcommand{\fS}{\ssub{\mathfrak S}!}

\renewcommand{\div}{\mathrm{div}}

\renewcommand{\k}{\kappa}
\newcommand{\K}{\mathbf{K}}

\newcommand{\Spec}{\operatorname{Spec}}
\newcommand{\g}{\mathfrak g}
\newcommand{\W}{{\scaleto{W}{4pt}}}

\newcommand{\an}{\operatorname{an}}

\newcommand{\T}{\ssub{\mathrm{T}}!}
\newcommand{\TF}{\mathbb{T}}
\DeclareMathOperator\val{\mathbf{val}}

\newcommand{\grdeg}{{\delta}}

\newcommand{\ee}{{\mathrm{e}}}
\newcommand{\vv}{{\mathrm{v}}}

\newcommand{\alg}{{^{\hspace{-.02cm}\mathrm{clls}}}}

\newcommand{\rest}[1]{\raisebox{-1pt}{$\vert$}_{#1}}

\newcommand{\st}{\bigm|} 

\newcommand{\Ical}{\mathcal I}

\newcommand{\rhostandard}{\rho^{\scaleto{\mathrm{st}}{4pt}}}

\newcommand{\rquot}[2]{#1\big/#2} 

\usepackage{manfnt}

\newcommand{\cube}{\mbox{\,\mancube}}

\newcommand{\comp}[1]{\overbar[.5]{#1}} 

\newcommand{\mg}{\mathscr M} 

\makeatletter
\newsavebox\myboxA
\newsavebox\myboxB
\newlength\mylenA

\newcommand*\overbar[2][0.65]{%
	\sbox{\myboxA}{$\m@th#2$}%
	\setbox\myboxB\null
	\ht\myboxB=\ht\myboxA%
	\dp\myboxB=\dp\myboxA%
	\wd\myboxB=#1\wd\myboxA
	\sbox\myboxB{$\m@th\overline{\copy\myboxB}$}
	\setlength\mylenA{\the\wd\myboxA}
	\addtolength\mylenA{-\the\wd\myboxB}%
	\ifdim\wd\myboxB<\wd\myboxA%
		\rlap{\hskip .9\mylenA\usebox\myboxB}{\usebox\myboxA}%
	\else
		\hskip -0.5\mylenA\rlap{\usebox\myboxA}{\hskip 0.5\mylenA\usebox\myboxB}%
	\fi}

\newcommand{\threestars}{(\hspace{-.15cm}\begin{tabular}{c}$\ast$\\[-1.6ex] $\ast\ast $\end{tabular}\hspace{-.15cm})}

\NewDocumentCommand{\ssub}{O{0pt} O{0pt} O{.9} m t! e{_^}}{
	#4%
	\IfValueT{#6}{
		\IfBooleanTF{#5}{\sb{\hspace{#1}\scaleobj{#3}{#6}}}{\sb{#6}}
	}
	\IfValueT{#7}{
		\IfBooleanTF{#5}{\sp{\hspace{#2}\scaleobj{#3}{#7}}}{\sp{#7}}
}
}

\ExplSyntaxOn
\NewDocumentCommand{\tossub}{o o m}{
	\expandafter\let\csname old\cs_to_str:N #3\endcsname#3
	\renewcommand#3%
	{\ssub[#1][#2]{\csname old\cs_to_str:N #3\endcsname}}
}
\ExplSyntaxOff

\tossub\omega
\tossub\Rat
\tossub\phi
\newcommand{\wtGamma}{\ssub[-1pt]{\widetilde \Gamma}!}
\newcommand{\ssM}{\ssub{M}!}
\newcommand{\mgg}[1]{\ssub{\mathscr M}!_{#1}}
\newcommand{\mgbarg}[1]{\comp{\ssub{\mg}!_{#1}}} 

\newcommand{\valuation}{\mathfrak v}
\newcommand{\vR}{\mathrm{R}}
\newcommand{\fm}{\mathfrak{m}}
\newcommand{\rmC}{\mathrm{C}}
\newcommand{\sstildef}{\ssub{\tilde f}!}

\newcommand{\sstildeH}{\ssub{\widetilde H}!}

\newcommand{\ssa}{\ssub{a}!}
\newcommand{\ssb}{\ssub{b}!}
\newcommand{\ssx}{\ssub{x}!}
\newcommand{\ssy}{\ssub{y}!}
\newcommand{\ssn}{\ssub{n}!}
\newcommand{\ssrho}{\ssub{\rho}!}
\newcommand{\ssvarepsilon}{\ssub{\varepsilon}!}
\newcommand{\ssi}{\ssub{i}!}
\newcommand{\ssI}{\ssub{I}!}
\newcommand{\ssJ}{\ssub{J}!}
\newcommand{\ssd}{\ssub{d}!}
\newcommand{\ssdelta}{\ssub{\partial}!}
\newcommand{\sse}{\ssub{e}!}
\newcommand{\ssE}{\ssub{E}!}

\newcommand{\ssf}{\ssub{f}!}
\newcommand{\ssg}{\ssub{g}!}
\newcommand{\ssbarf}{\ssub{\bar f}!}

\newcommand{\ssh}{\ssub{h}!}
\newcommand{\ssj}{\ssub{j}!}
\newcommand{\ssell}{\ssub{\ell}!}
\newcommand{\sslambda}{\ssub{\lambda}!}
\newcommand{\ssu}{\ssub{u}!}
\newcommand{\ssv}{\ssub{v}!}
\newcommand{\ssc}{\ssub{c}!}
\newcommand{\ssC}{\ssub{C}!}
\newcommand{\ssD}{\ssub{D}!}
\newcommand{\sssim}{\ssub{\sim}!}
\newcommand{\ssnsim}{\ssub{\nsim}!}

\newcommand{\ssleq}{\ssub{\leq}!}
\newcommand{\ssle}{\ssub{<}!}

\newcommand{\sst}{\ssub{t}!}

\newcommand{\ssfilt}{\ssub{\filt}!}
\newcommand{\ssm}{\ssub{m}!}
\newcommand{\ssS}{\ssub{S}!}
\newcommand{\sss}{\ssub{s}!}
\newcommand{\ssp}{\ssub{p}!}
\newcommand{\ssmathcalP}{\ssub{\mathcal P}!}
\newcommand{\ssq}{\ssub{q}!}
\newcommand{\ssX}{\ssub{X}!}
\newcommand{\ssY}{\ssub{Y}!}
\newcommand{\ssH}{\ssub{H}!}
\newcommand{\ssR}{\ssub{R}!}

\newcommand{\ssL}{\ssub{L}!}

\newcommand{\ssA}{\ssub{\mathbb A}!}

\newcommand{\varT}{T}

\newcommand{\e}{\ssub{\mathfrak e}} 

\newcommand{\unde}{\ssub[1pt]{\underline{\e}}!}
\newcommand{\undc}{\ssub{\underline c}!}
\newcommand{\undb}{\ssub{\underline b}!}
\newcommand{\unda}{\ssub{\underline a}!}
\newcommand{\undx}{\ssub{\underline x}!}
\newcommand{\undy}{\ssub{\underline y}!}
\newcommand{\undz}{\ssub{\underline z}!}
\newcommand{\undw}{\ssub{\underline w}!}
\newcommand{\ssP}[1]{\ssub{P}!_{\hspace{-2pt}#1}}
\newcommand{\unds}{\ssub{\underline s}!}

\newcommand{\undq}{\ssub{\underline q}!}

\newcommand{\rk}{\ssub{\mathbf{r}}!}

\newcommand{\VS}{\mathrm{H}}
\newcommand{\filt}{\mathrm{F}}

\newcommand{\subfaceeq}{\preceq}
\newcommand{\supfaceeq}{\succeq}
\newcommand{\subface}{\prec}
\newcommand{\supface}{\succ}

\renewcommand{\curve}{\mathbf{C}} 
\newcommand{\varD}{\mathbf{D}}
\newcommand{\varE}{\mathbf{E}}

\newcommand{\varH}{\mathbf{H}}
\newcommand{\varf}{\ssub{\mathbf{f}}!}
\newcommand{\varg}{\mathbf{g}}
\newcommand{\varh}{\mathbf{h}}

\newcommand{\hcube}[2]{\ssub[0pt,2pt]{\cube}_{#1}^{^{#2}}}

\renewcommand{\setminus}{\smallsetminus}
\newcommand{\matricube}{\ssub{\mathscr M}!}
\newcommand{\unimatricube}{\ssub{\mathscr U}!}
\newcommand{\undr}{\underline r}

\newcommand{\Fl}{\mathscr F}

\title{Limit linear series: combinatorial theory}

\author{Omid Amini}
\address{CNRS - CMLS, \'Ecole polytechnique, Institut polytechnique de Paris}
\email{\href{omid.amini@polytechnique.edu}{omid.amini@polytechnique.edu}}

\author{Lucas Gierczak}
\address{CMLS, \'Ecole polytechnique, Institut polytechnique de Paris}
\email{\href{lucas.gierczak-galle@polytechnique.edu}{lucas.gierczak-galle@polytechnique.edu}}

\date{\today}
\begin{document}

\begin{abstract}
	We develop a purely combinatorial theory of limit linear series on metric graphs. This will be based on the formalisms of hypercube rank functions and slope structures. We provide a full classification of combinatorial limit linear series of rank one, and discuss connections to other concepts in tropical algebra and combinatorial algebraic geometry.
\end{abstract}

\maketitle

\setcounter{tocdepth}{1}

\tableofcontents

\section{Overview}
	
	One of the longstanding open questions regarding the asymptotic geometry of curves is the problem of degeneration of linear series on smooth curves of given genus when they approach the boundary of their corresponding moduli space. That is, to fully describe all the possible limits of linear series of given rank and degree when smooth curves degenerate to singular ones. This question was studied in a series of works by Eisenbud and Harris \cite{EHBull, EH86, EH87-KD, EH87-WP, EH87mon} for which they managed to provide a satisfactory answer in the case the limit curve is of compact type, and used this to make major progress in the study of curves. For curves of pseudo-compact type, these results were generalized by Esteves--Medeiros~\cite{EM02} (for curves with two components) and by Osserman~\cite{Os-pseudocompact, Oss19, He19}. The case of rank zero linear series in the pluricanonical systems is studied in recent work~\cite{BCGGM18, BCGGM19, MUW20, TT22}.
 	
	\smallskip
	
	Tropical geometry provides a modern perspective on degeneration methods in algebraic geometry and a new approach to classical questions in algebraic geometry. Developing the mathematics behind the tropical approach usually requires the introduction of new combinatorial structures, and it has become apparent now that from the viewpoint of applications, it is enough in many cases to understand the geometry behind these combinatorics. Two such examples are given in the development of a combinatorial theory of divisors on graphs and metric graphs~\cite{baker2007riemann, BJ16}, and in the development of tropical and combinatorial Hodge theories~\cite{AHK18, amini2020hodge}.
	
	\smallskip
	
	In a previous work~\cite{AB15}, the first named author and Baker introduced linear series on hybrid objects called \emph{metrized complexes} and used them to recover and partially generalize the Eisenbud--Harris theory of limit linear series. In a subsequent work~\cite{Ami14}, the formalism of slope structures on metric graphs was introduced as a way to describe the limiting behavior of Weierstrass points on degenerating families of curves. Slope structures were used in a recent work of Farkas--Jensen--Payne~\cite{FJP20} in the study of the geometry of moduli space of curves.
	
	The aim of this paper is to take the tropical approach one step further by introducing a purely combinatorial theory of linear series of arbitrary rank and degree on metric graphs. This can be regarded as a combinatorial theory of limit linear series.
	
	More precisely, we aim to draw relevant combinatorial properties of tropicalizations of linear series, regarding:
	
	\begin{itemize}
		\item the slopes taken by the tropicalizations of functions;
		
		\item the vectors of slopes taken by these functions around points;
		
		\item the tropical dependence between these functions, and
		
		\item the topological properties of tropicalizations,
	\end{itemize}
	in order to develop a formalism of linear series on metric graphs. This will be based on two ingredients: hypercube rank functions and slope structures. In the rest of this introduction, we give an overview of the setup and the results.
	
	\subsection{Degeneration problem for linear series}
		
		Recall that a linear series $\g^r_d$ on a projective curve $Y$ is by definition a vector subspace of dimension $r + 1$ of the space of global sections $H^0(Y, L)$ of a line bundle $L$ of degree $d$ on $Y$.
		
		Let $\mgg{g}$ be the moduli space of smooth projective curves of genus $g$, and $\mgbarg{g}$ its Deligne--Mumford compactification. Let $X$ be a stable curve of genus $g$ over an algebraically closed field $\kappa$ and $x$ the corresponding point in $\mgbarg{g}$. The degeneration problem for linear series can be informally stated as follows.
		
		\begin{question} \label{ques:main}
			Describe all the possible limits of linear series over any sequence of smooth projective curves of genus $g$ when their corresponding points in $\mgg{g}$ converge to $x$.
		\end{question}
	
	\subsection{Metric graphs and their divisor theory}
		
		Metric graphs arise as tropical limits of one-parameter families of smooth projective curves.
		
		We denote by $\R_+$ the set of positive real numbers. Recall that a metric graph $\Gamma$ is a compact metric space isomorphic to the metric realization of a pair $(G, \ell)$ consisting of a finite graph $G = (V, E)$ and a length function $\ell \colon E \to \R_+$: this is obtained by associating to each edge $e$ a copy of the interval $\ssub{\Ical}!_e = [0, \ssell_e]$, with the two extremities identified with those of $e$, and then further identifying the ends of different intervals corresponding to a same vertex $v$. The quotient topology on $\Gamma$ is metrizable by the path metric. The pair $(G, \ell)$ is called a \emph{model} of $\Gamma$.
		
		In the context related to the degeneration of algebraic curves, a metric graph $\Gamma$ with model $(G = (V, E), \ell)$ is endowed with a function $\g \colon V \to \Z_{\geq 0}$ associating to each vertex $v$ the genus of some algebraic curve $C_v$ represented by $v$. Such a triple $(G, \ell, \g)$ will be called an \emph{augmented metric graph}. However, we will mostly handle non-augmented metric graphs in this paper, and make comments about the relevance of the genus function in the theory developed here.
		
		The set of \emph{rational functions} on $\Gamma$ is denoted by $\Rat(\Gamma)$, and, by definition, consists of all continuous piecewise affine functions $f: \Gamma \rightarrow \mathbb R$ with integral slopes. The tropicalization of rational functions on curves gives rise to rational functions on metric graphs.
		
		As in the algebraic setting, rational functions on metric graphs are linked to divisors. A divisor $D$ on a metric graph $\Gamma$ is a finite formal sum with integer coefficients of points of $\Gamma$. 
		For any rational function $f \in \Rat(\Gamma)$, the corresponding divisor of zeros and poles is defined by
		\[ \div(f) \coloneqq \sum_{x} \ord_x(f) \, (x), \qquad \textrm{with } \ord_x(f) \coloneqq - \sum_{\nu \in \T_x(\Gamma)} \mathrm{slope}_{\nu}(f), \]
		where $\T_x(\Gamma)$ is the set of outgoing unit tangent vectors to $\Gamma$ at $x$, and $\mathrm{slope}_{\nu}(f)$ is the slope of $f$ along $\nu$ at $x$. A divisor obtained in this way is called \emph{principal}, and two divisors $D_1$ and $D_2$ whose difference $D_1 - D_2$ is principal are called \emph{linearly equivalent}.
		
		Each divisor $D$ gives rise to a line bundle on $\Gamma$ whose space of global sections is denoted by $\Rat(D)$. Concretely, $\Rat(D)$ is the set of $f \in \Rat(\Gamma)$ such that $\div(f) + D$ is effective (that is, has only non-negative coefficients). Contrary to the algebraic setting, $\Rat(D)$ is not a vector space. Nevertheless, Baker and Norine discovered a way to associate a combinatorial notion of rank to $\Rat(D)$~\cite{baker2007riemann, MZ08, gathmann2008riemann}. This is defined as the maximum integer among $-1$ and integers $r \geq 0$ such that for all points $\ssx_1, \dots, \ssx_r$ in $\Gamma$, the divisor $D - (\ssx_1) - \cdots - (\ssx_r)$ is linearly equivalent to an effective divisor. We refer to the survey paper by Baker and Jensen~\cite{BJ16} for details, extensions, and several applications.
	 	
	 	\smallskip
	 	
	 	The results of the present paper are motivated by the question of describing the tropical limits of linear series, when the metric graph arises as the tropical limit of a one-parameter family of smooth proper curves. In the next two sections, we describe two combinatorial structures that allow to approach this question.
	
	\subsection{Rank functions on hypercubes}
		
		We first describe a combinatorial way of encoding the intersection patterns of a flag arrangement.
		
		Let $r$ be a non-negative integer. We set $[r] \coloneqq \{0, \dots, r\}$. For a positive integer $\grdeg$, the hypercube $\hcube{r}{\grdeg}$ of dimension $\grdeg$ and width $r$ is the product $[r]^\grdeg$.
		
		We define a partial order $\subfaceeq$ on $\hcube{r}{\grdeg}$ where for a pair of elements $\unda, \undb$ in $\hcube{r}{\grdeg}$, we write $\unda \subfaceeq \undb$ whenever for every $j \in \{1, \dots, \grdeg\}$, we have $\ssa_j \leq \ssb_j$. 
		We define two operations $\vee$ and $\wedge$ by taking the maximum and the minimum coordinate-wise, respectively: for $\unda = (a_1, \dots, a_\grdeg)$ and $\undb = (b_1, \dots, b_\grdeg)$ in $\hcube{r}{\grdeg}$,
		\[ \unda \vee \undb \coloneqq (\max(a_1, b_1), \dots, \max(a_\grdeg, b_\grdeg)), \qquad \unda \wedge \undb \coloneqq (\min(a_1, b_1), \dots, \min(a_\grdeg, b_\grdeg)). \]
		A function $f: \hcube{r}{\grdeg} \rightarrow \mathbb Z$ is called \emph{supermodular} if for every two elements $\unda$ and $\undb$, we have
		\[ f(\unda) + f(\undb) \leq f(\unda \vee \undb) + f(\unda \wedge \undb). \]
		
		A function $\rho: \hcube{r}{\grdeg} \rightarrow \mathbb Z$ is called a \emph{(hypercube) rank function} if it is supermodular and, in addition, satisfies the following conditions:
		
		\begin{enumerate}[label=(\arabic*)]
			\item For every $i \in \{1, \dots, \grdeg\}$, and each $t \in [r]$, we have $\rho\mleft(t \, \unde_i\mright) = r - t$.
			
			\item $\rho$ is non-increasing with respect to the partial order of $\hcube{r}{\grdeg}$, that is, if $\unda \subfaceeq \undb$, then $\rho(\undb) \leq \rho(\unda)$.
			
			\item The values of $\rho$ are in the set $\{-1, 0, 1, \dots, r\}$.
		\end{enumerate}
		
		The geometric situation to have in mind in order to interpret hypercube rank functions is a vector space $\VS$ of dimension $r + 1$ over some field $\k$, and a collection of $\grdeg$ complete flags $\ssfilt^\bullet_1, \dots, \ssfilt^{\bullet}_\grdeg$. That is, for $j = 1, \dots, \grdeg$, $F^\bullet_j$ consists of a chain of vector subspaces
		\[ \VS = \ssfilt_j^0 \supsetneq \ssfilt_j^1 \supsetneq \dots \supsetneq \ssfilt_j^{r - 1} \supsetneq \ssfilt_j^r \supsetneq (0), \]
		with $\mathrm{codim}(\ssfilt_j^i) = i$. In this case, the function $\rho \colon \hcube{r}{\grdeg} \rightarrow \mathbb Z$ defined by
		\[ \rho(a_1, \dots, a_\grdeg) \coloneqq \dim_\k \mleft(\ssfilt_1^{\ssa_1} \cap \dots \cap \ssfilt_\grdeg^{\ssa_\grdeg}\mright) -1 \]
		is a rank function. Hypercube rank functions appearing in this way are called \emph{representable}. 
		On a smooth projective curve, given a finite dimensional vector space of rational functions, taking the orders of vanishing at a point leads to a complete flag (see Section~\ref{subsubsec:geometric_rank_functions}). Such flags appear naturally in the degeneration of linear series, see Section~\ref{subsec:linear_series_specialization}.
	
		We say that a point $\unda$ of the hypercube is a \emph{jump} of the rank function $\rho$ if we have $\ssrho_v(\unda + \unde_j) < \ssrho_v(\unda)$ whenever $\unda + \unde_j$ is in the hypercube. Here, $\unde_j$ is the point whose $j$-th coordinate is equal to one and whose other coordinates are all zero. In the representable case, the jumps correspond to extremal points beyond which the dimension drops. It turns out that the jumps uniquely determine the rank function.
		
	\subsection{Slope structures on graphs and metric graphs}
		
		Slope structures encode the information regarding possible slopes of functions arising from tropicalizations of linear series.
		
		For a graph $G = (V, E)$, denote by $\mathbb E$ the set of all possible orientations we can put on edges of $G$, that is, each edge $\{u, v\} \in E$ gives rise to two oriented edges $uv$ and $vu$ in $\mathbb E$. We suppose that $G$ is simple, that is, it has no loop (an edge with identical extremities) and no parallel edges (edges with the same extremities). 
		Let $\mathbb E_v \subseteq \mathbb E$ be the set of oriented edges $uv$ in $G$, for all $\{u, v\} \in E$, and denote by $d_v = |\mathbb E_v|$ the valence of $v$.
		
		A \emph{slope structure $\fS$ of width $r$} on $G$, or simply an \emph{$r$-slope structure}, is the data of
		
		\begin{itemize}
			\item For each oriented edge $e = uv \in \mathbb E$ of $G$, a collection $S^e$ of $r + 1$ integers $\sss_0^e < \sss_1^e < \dots < \sss_r^e$ subject to the requirement that $\sss^{uv}_i + \sss^{vu}_{r - i} = 0$ for each edge $\{u, v\} \in E$.
			
			\item For each vertex $v$ of $G$, a rank function $\ssrho_v$ on the hypercube $\hcube{r}{\ssd_v}$.
		\end{itemize}
		
		We denote by $S^v \subseteq \prod_{e \in \mathbb E_v} S^e$ the set of all points $\sss_{\unda} = (\sss_{\ssa_1}, \dots, \sss_{\ssa_\grdeg})$ for $\unda = (\ssa_1, \dots, \ssa_{\grdeg})$ a jump of the rank function $\ssrho_v$. Since jumps determine the rank function, using the above notation, we then write $\fS = \mleft\{S^v; S^e\mright\}_{v \in V, e \in \mathbb E}$.
		
		Let $\Gamma$ be a metric graph. By an $r$-slope structure on $\Gamma$, we mean an $r$-slope structure $\fS$ on a simple graph model $G = (V, E)$ of $\Gamma$ that we naturally extend to any point of $\Gamma$ by associating to every point $x \in \Gamma \setminus V$ the so-called \emph{standard rank function}, see Section~\ref{subsec:rank_function_induced_flags}. More precisely, for every point $x \in \Gamma \setminus V$ and each outgoing unit tangent vector $\nu \in \T_x(\Gamma)$, we define $\ssS^\nu = \ssS^{uv}$, where $uv$ is the unique oriented edge of $G$ which is parallel to $\nu$, and define $\ssS^x \subseteq \ssS^{uv} \times \ssS^{vu}$ as the set of all pairs $(\sss_i^{uv}, \sss_j^{vu})$ with $i + j \leq r$.
		
		Taking into account the slope structure $\fS$, we can define a relevant notion of rational functions. A function $f$ in $\Rat(\Gamma)$ is said to be compatible with $\fS$ if the two conditions $(i)$ and $(ii)$ below are verified. First,	
		\begin{itemize}
			\item[$(i)$] for any point $x \in \Gamma$ and any tangent direction $\nu \in \T_x(\Gamma)$, the outgoing slope of $f$ along $\nu$ lies in $S^\nu$.
		\end{itemize}
		Denote by $\ssdelta_x(f)$ the vector in $\prod_{\nu \in \T_x(\Gamma)} S^\nu$ which consists of outgoing slopes of $f$ along $\nu \in \T_x(\Gamma)$. Then the second condition is:
		
		\begin{itemize}
			\item[$(ii)$] for any point $x \in \Gamma$, the vector $\ssdelta_x(f)$ belongs to $S^x$.
		\end{itemize}
		
		\medskip
		
		We denote by $\Rat(\Gamma, \fS) \subset \Rat(\Gamma)$, or simply $\Rat(\fS)$ if $\Gamma$ is understood, the set of rational functions on $\Gamma$ compatible with $\fS$. Endowed with the operations $c \odot f \coloneqq f + c$ and $f \oplus g \coloneqq \min(f, g)$ for every $f, g \in \Rat(\fS)$ and $c \in \R$, $\Rat(\fS)$ has the structure of a tropical semimodule. Moreover, it is naturally endowed with the norm $\| \cdot \|_\infty$.
	
		If $D$ is a divisor on $\Gamma$, we denote by $\Rat(D, \fS) \coloneqq \Rat(D) \cap \Rat(\fS)$ the set of all $f \in \Rat(\fS)$ such that $D + \div(f)$ is effective. This is a sub-semimodule of $\Rat(\fS)$ (see Proposition~\ref{prop:structure_linear_series}).
	
	\subsection{Admissible semimodules and combinatorial limit linear series} 
		
		A semimodule $M \subseteq \Rat(D)$ is called \emph{admissible of rank $r$} if it is closed for the topology induced by $\| \cdot \|_\infty$ (equivalently, for that of point-wise convergence), and if there exists an $r$-slope structure $\fS$ such that $M \subseteq \Rat(D, \fS)$ and the following holds:
		
		\smallskip
		
		$(**)$ For every effective divisor $E$ of degree $r$, there exists $f \in M$ such that
		
		\begin{enumerate}[label=(\arabic*)]
			\item For every point $x \in \Gamma$, $\ssrho_x(\ssdelta_x(f)) \geq E(x)$; and in addition,
			
			\item $D + \div(f) - E \geq 0$.
		\end{enumerate}
		
		In the degeneration picture for linear series in a one-parameter family of smooth projective curves, viewing $x$ inside the Berkovich analytification of the generic fiber in the family, the first condition reflects the dimension counts underlying the reduction at $x$ of rational functions in the linear series, imposing vanishing conditions along the incident branches. The second one is the analogue of the Baker--Norine rank condition in this setting.
		
		In order to define a notion of linear series, we use moreover the concept of tropical rank introduced by Jensen and Payne in their work on applications of tropical divisor theory to the study of the geometry of generic curves~\cite{JP14, JP16} (for the definition, see Section~\ref{subsec:structure_rat}).
		
		\smallskip

		A \emph{(combinatorial limit) linear series} of rank $r$ and degree $d$, more simply called a $\g^r_d$, is a pair $(D, M)$ consisting of a divisor $D$ of degree $d$ and an admissible semimodule $M \subseteq \Rat(D)$ of rank $r$ which is moreover \emph{finitely generated} and has \emph{tropical rank} $r$.
		
		Here, we mean by ``finitely generated'' that there exist finitely many elements of $M$ which generate $M$ using the tropical operations of scalar addition and minimum; the tropical rank is defined as the maximal number of tropically independent elements of $M$. See Sections~\ref{subsec:finite_generation_closedness} and~\ref{subsec:tropical_rank} for more details.
		
		If $M$ is a $\g^r_d$, the \emph{linear system} $|M|$ is the space of all effective divisors $E$ on $\Gamma$ of the form $D + \div(f)$ for $f \in M$.
		
		\smallskip
		
		We also define \emph{refined linear series} to be those linear series which in addition verify the following stronger version of $(**)$:
		
		\medskip
		
		$\threestars$ For any effective divisor $E$ on $\Gamma$ of degree $s\leq r$, there exists a linear series $M_E$ of rank $r - s$ associated to $(D, \fS_E)$ with $\fS_E$ a slope substructure of $\fS$ of rank $r - s$, such that for every function $f \in M_E$, we have
		
		\begin{enumerate}[label=(\arabic*)]
			\item For every point $x \in \Gamma$, $\ssrho_x(\ssdelta_x(f)) \geq E(x)$; and in addition,
			
			\item $D + \div(f) - E \geq 0$.
		\end{enumerate}
	
	\subsection{Realization property for slope vectors}
		
		An important feature of admissible semimodules is the following realization property for jumps, proved in Section~\ref{subsec:realization_jumps}.
		
		\begin{thm}[Realization of slope vectors] \label{thm:jumps_realized_intro}
			Let $(D, \fS)$ be a pair consisting of a divisor $D$ and a slope structure $\fS$ of width $r$ on $\Gamma$, and let $M \subseteq \Rat(D, \fS)$ be an admissible semimodule of rank $r$. Let $v$ be a point of $\Gamma$ and let $\unda$ be a jump of $\ssrho_v$. Then, there exists $f \in M$ such that $\ssdelta_v(f) = \sss_{\unda}$.
		\end{thm}
		
		One immediate consequence is that $\fS$ can be entirely retrieved from $M$.
		
		\begin{cor} \label{cor:slope_structure_retrieved_linear_series_intro}
			The data of an admissible semimodule $M \subset \Rat(D)$ determines the slope structure $\fS$ uniquely.
		\end{cor}
	
	\subsection{Finiteness of slope structures}
		
		Let $\Gamma$ be a metric graph, let $D$ be a divisor on $\Gamma$. Let $G = (V, E)$ be a combinatorial graph underlying $\Gamma$ and supporting $D$. Let $M \subset \Rat(D)$ be an admissible semimodule with underlying slope structure $\fS$ defined on some model of $\Gamma$ (possibly different from $G$). We prove in Proposition~\ref{prop:non_increasing_slope_vectors} that the vector of slopes in $\fS$ is non-increasing along each edge. That is, as we move from one extremity of an edge to the other, the coordinates of the vector $(\sss_0^\nu, \dots, \sss^\nu_r)$ do not increase. This property turns out to be crucial in proving the following finiteness theorem (see Section~\ref{subsec:slope_vectors_non_increasing} for a more general result).

		\begin{thm}[Finiteness of slopes structures] \label{thm:finiteness_slope_structures_intro}
			For each integer $r$, there are finitely many subdivisions $\ssH_1, \dots, \ssH_k$ of $G$, and finitely many slope structures $\fS_1, \dots, \fS_k$ of rank $r$ defined on them, respectively, such that every admissible semimodule $M \subset \Rat(D)$ with underlying slope structure $\fS$ has a combinatorial model $\ssH_j$ among $\ssH_1, \dots, \ssH_k$ such that $\fS = \fS_j$.
		\end{thm}
		
		This can be regarded as a first result in the direction of defining the moduli space of $\g^r_d$'s over the moduli space of tropical curves of given genus.
	
	\subsection{Classification of $\g^1_d$'s}
		
		In the case $r = 1$, we prove in Section~\ref{sec:g1d} that the data of a $\g^1_d$ on $\Gamma$ is equivalent to the data of a finite harmonic map to a tree.
		
		\begin{thm}[Classification of $\g^1_d$'s on metric graphs] \label{thm:g1d_intro}
			Let $(D, M)$ be a $\g^1_d$ on $\Gamma$ with $D$ a divisor of degree $d$. Suppose that the constant functions are in $M$. Then, there exist a tropical modification $\alpha \colon \widetilde \Gamma \longrightarrow \Gamma$ of $\Gamma$ and a finite harmonic morphism $\varphi \colon \widetilde \Gamma \rightarrow T$ of degree $d$ to a metric tree $T$ such that $M$ is the preimage of the unique $\g^1_1$ on $T$, restricted to $\Gamma$.
		\end{thm}
		
		Using the smoothing theorems proved in~\cite{amini2015lifting1, amini2015lifting2} for finite harmonic morphisms to trees, we deduce the following smoothing theorem for combinatorial $\g^1_d$'s (see below for the tropicalization, and Section~\ref{subsec:linear_series_specialization}).
		
		\begin{thm}[Smoothing theorem for $\g^1_d$'s] \label{thm:g1d_smoothing_intro}
			A $\g^1_d$ $(D, M)$ on $\Gamma$ is smoothable, that is, it is the tropicalization of a $\g^1_d$ from a smooth curve.
		\end{thm}
		
		The question of the existence of harmonic morphisms to metric trees of smallest degree, called \emph{geometric gonality}, is thoroughly studied by Draisma and Vargas~\cite{DV21}, and by Cool and Draisma~\cite{CD18}. The above theorem provides an algebraic characterization of geometric gonality, as the least integer $d$ such that the metric graph admits a $\g^1_d$, see~\cite{DV21bis} and the references there for the gonality of metric graphs.
		
		In order to prove the above theorem, we generalize in Section~\ref{sec:reduced_divisors} the theory of reduced divisors to the setting of linear series, consider the map defined by reduced divisors in Section~\ref{sec:g1d}, and connect it to the tropical rank to conclude.
	
	\subsection{Tropicalization}
		
		We now discuss the connection between linear series on algebraic curves and their combinatorial counterparts. We assume familiarity with the Berkovich theory of algebraic curves, see Section~\ref{sec:geometric_tropicalization} and~\cite{BJ16} for more details.
		
		Let $\K$ be an algebraically closed field with a non-trivial non-Archimedean valuation $\val$ and $\curve$ be a smooth proper curve over $\K$. We assume that $\K$ is complete with respect to $\val$ and we denote by $\k$ the residue field of $\K$, which is also algebraically closed. Denote by $\curve^{\an}$ the Berkovich analytification of $\curve$. Let $\Gamma \subset \curve^{\an}$ be a metric graph skeleton of $\curve^{\an}$. For each point $x \in \Gamma$, let $\valuation_x$ be the valuation on the function field $\K(\curve)$ defined by $x$. For each nonzero rational function $\varf \in \K(\curve)$, the tropicalization of $\varf$, denoted $\Trop(\varf) \colon \Gamma \to \R$, is defined by $\Trop(\varf)(x) = \valuation_x(f)$ for all $x \in \Gamma$. This is a piecewise affine function on $\Gamma$ with integral slopes.
		
		Let $\varD$ be a divisor of degree $d$ and rank at least $r$ on $\curve$, and $(\mathcal O(\varD), \varH)$ be a $\g^r_d$ on $\curve$, so $\varH$ is an $(r + 1)$-dimensional vector subspace of $H^0(\curve, \mathcal O(\varD))$, the space of global sections of the line bundle associated to $\varD$. We can identify $\varH$ with a subspace of $\K(\curve)$ of dimension $r + 1$. We define the \emph{tropicalization}
		\[ M \coloneqq \Trop(\varH) = \mleft\{\Trop(\varf) \st \varf \in \varH \setminus \{0\}\mright\}. \]
		 The following theorem is proved in Section~\ref{subsec:linear_series_specialization}.
		
		\begin{thm}[Tropicalization of linear series] \label{thm:specialization_linear_series_grd_intro}
			Notation as above, let $(\mathcal O(\varD), \varH)$, $\varH \subseteq H^0\mleft(\curve, \mathcal O(\varD)\mright) \subset \K(\curve)$, be a $g^r_d$ on $\curve$. Let $\Gamma$ be a skeleton of $\curve^{\an}$. The slopes of rational functions $F$ in $M$ along edges in $\Gamma$ yield a well-defined slope structure $\fS$ on $\Gamma$. Let $D$ be the tropicalization of $\varD$ to $\Gamma$. Then, $M \subseteq \Rat(D, \fS)$ is a refined $\g^r_d$ on $\Gamma$.
		\end{thm}

		In particular, the tropicalization $\Trop(\varH)$ is finitely generated (see Proposition~\ref{prop:finite_generation}).

	\subsection{Discussion of applications}
		
		The formalism of this paper has applications to the geometry of curves. In particular, the equidistribution theorem proved in~\cite{Ami14} is a consequence of the formalism of slope structures and the behavior of reduced divisors in a given combinatorial linear series. The results of~\cite{FJP20} also use slope structures and the notion of tropical independence, the underlying concepts of the materials presented in this paper. In a joint work of the authors with Harry Richman~\cite{AGR23}, we apply the formalism of this paper to associate a Weierstrass weight to each connected component of the naive Weierstrass locus of a given divisor on a metric graph. This solves a problem posed by Matt Baker from his original work on the specialization of linear series from curves to graphs~\cite{baker2008specialization}. Using these ideas, we explain the discrepancy between the naive counting of Weierstrass points on metric graphs in the work of Richman~\cite{richman2024distribution} and the correct count of multiplicities.
		
		The hypercube rank functions considered here are linked to matricubes and a theory of combinatorial flag arrangements developed in a companion work~\cite{AG24}, which provide a generalization of the theory of matroids. For more information about the connection, we refer to Section~\ref{subsec:definition}, and more specifically to Proposition~\ref{prop:duality}.
		
		\begin{remark}
			When this paper was in the final stages of preparation, we learned of a notion of tropical linear series on metric graphs developed simultaneously and independently by Dave Jensen and Sam Payne~\cite{JP22}, with some similarities to the framework introduced and studied here. There are also meaningful differences between the two approaches. One of the interesting features of the work~\cite{JP22} which does not appear in this paper is the valuated matroid property that allows to put global restrictions on the space of relations between elements of the tropical semimodule, and to define maps to tropical linear spaces, as in the classical setting. The recursive condition on the existence of sublinear series is imposed on tangent vectors, which add flexibility to the proofs. Clearly, understanding the precise connection between the two approaches, and making a combination of the properties discovered in both works, should lead to further progress in the development of tropical methods and their applications to the study of the geometry of curves and their moduli spaces.
		\end{remark}
	
	\subsection*{Acknowledgments}
		
		We thank Eduardo Esteves, St\'ephane Gaubert, and Noema Nicolussi for discussions related to the content of this work. We are grateful to Harry Richman for valuable comments on an earlier draft, for the idea of the example given in Section~\ref{subsec:realizability_question}, and for collaboration on related questions. Special thanks to Dave Jensen and Sam Payne for sharing a draft of their work, and for interesting exchanges on the content of the two papers. In particular, we thank them for bringing to our attention the fact that the tropicalization of a finite dimensional vector space is a finitely generated semimodule, which allows to conclude the proof of Proposition~\ref{prop:finite_generation}, as well as St\'ephane Gaubert for further references. Finally, it is a pleasure to thank Matt Baker for his inspiring work and questions, as well as for previous collaboration on related subjects.
	
	\subsection*{Further notation} \label{subsec:basic_notions}
		
		For a subset $S \subseteq \R^d$, $d \in \N$, and an element $a \in \R^d$, we define $S + a$ as the set of all elements $b + a$ for $b \in S$. For the rest of the article, $\Lambda$ will denote a divisible subgroup of $\R$ -- a group is called \emph{divisible} when multiplication by every positive integer $n$ is surjective. Examples of such $\Lambda$ are $\Q$ or $\R$ itself.
		
		A \emph{$\Lambda$-metric graph}, or \emph{$\Lambda$-rational metric graph}, is a metric graph $\Gamma$ whose edge lengths are $\Lambda$-rational, that is, lie in $\Lambda$. A point $x \in \Gamma$ is said to be \emph{$\Lambda$-rational} if its distances to the endpoints of its incident edges are in $\Lambda$. In the case $\Gamma$ is $\Lambda$-rational, a divisor whose support is made up of $\Lambda$-rational points is said to be a \emph{$\Lambda$-rational} divisor. We also denote by $\Rat!_\Lambda(\Gamma)$ the set of functions of $\Rat(\Gamma)$ which only change slope at $\Lambda$-rational points of $\Gamma$ and which take a value in $\Lambda$ on some (equivalently, on every) vertex $v$ of $\Gamma$.
		
 		We denote by $\TF$ the semifield of tropical numbers, which is the set $\R$ endowed with the internal operations of tropical addition $\oplus \coloneqq \min$ and tropical multiplication $\odot \coloneqq +$. A (divisible) subgroup $\Lambda$ endowed with the operations $\oplus$ and $\odot$ defines a sub-semifield of $\TF$.

\section{Rank functions on hypercubes} \label{sec:rank}
	
	\subsection{Definition} \label{subsec:definition}
		
		Let $r$ be a non-negative integer and $[r] = \{0, 1, \dots, r\}$. For a positive integer $\grdeg$, the \emph{hypercube} $\hcube{r}{\grdeg}$ of \emph{dimension} $\grdeg$ and \emph{width} $r$ is the product $[r]^\grdeg$. We denote the elements of $\hcube{r}{\grdeg}$ by vectors $\unda = (\ssa_1, \dots, \ssa_\grdeg)$, for $0 \leq \ssa_1, \dots, \ssa_\grdeg \leq r$.
		
		The hypercube $\hcube{r}{\grdeg}$ is endowed with a natural partial order $\subfaceeq$ defined by declaring $\unda \subfaceeq \undb$ for elements $\unda = (\ssa_1, \dots, \ssa_\grdeg)$ and $\undb = (\ssb_1, \dots, \ssb_\grdeg)$ in $\hcube{r}{\grdeg}$, if $\ssa_j \leq \ssb_j$ for all $j \in \{1, \dots, \grdeg\}$. The smallest and largest elements of $\hcube{r}{\grdeg}$ with respect to this partial order are $\underline 0 \coloneqq (0, \dots, 0)$ and $\undr = (r, \dots, r)$, respectively. Moreover, there is a lattice structure on $\hcube{r}{\grdeg}$, where the two operations of join $\vee$ and meet $\wedge$ are defined by
		\[ 
			\unda \vee \undb = (\max(\ssa_1, \ssb_1), \dots, \max(\ssa_\grdeg, \ssb_\grdeg)), \quad \unda \wedge \undb = (\min(\ssa_1, \ssb_1), \dots, \min(\ssa_\grdeg, \ssb_\grdeg)) \quad \forall \unda, \undb \in \hcube{r}{\grdeg}.
		\]
		
		A function $f \colon \hcube{r}{\grdeg} \rightarrow \mathbb Z$ is called \emph{supermodular} if for every pair of elements $\unda$ and $\undb$, we have
		\[ f(\unda) + f(\undb) \leq f(\unda \vee \undb) + f(\unda \wedge \undb). \]
		If the inequalities above are all reversed, then we say that $f$ is \emph{submodular}. We define the \emph{conjugate} of $f$, denoted by $\ssbarf$, as the integer-valued function on the hypercube $\hcube{r}{\grdeg}$ given by
		\[ \ssbarf(\unda) \coloneqq r - f(\unda) \qquad \forall \, \unda \in \hcube{r}{\grdeg}. \]
		Note that $\bar{\ssbarf} = f$. Moreover, $\ssf$ is supermodular, resp. submodular, if, and only if, $\ssbarf$ is submodular, resp. supermodular.
		
		In this paper, we will be working with a special kind of supermodular function on $\hcube{r}{\grdeg}$. For each integer $i \in \{1, \dots, \grdeg\}$, denote by $\unde_i$ the vector whose coordinates are all zero except the $i$-th coordinate, which is equal to one. For $t \in [r]$, the vector $t \, \unde_i$ lies in $\hcube{r}{\grdeg}$.
		
		\begin{defi}[Hypercube rank function] \label{def:rank_function}
			A function $\rho \colon \hcube{r}{\grdeg} \rightarrow \mathbb Z$ is called a \emph{rank function} if it satisfies the following conditions:
			
			\begin{enumerate}[label=(HR\arabic*)]
				\item \label{def:axiom_rank_function_1} for every $i \in \{1, \dots, \grdeg\}$, and each $t \in [r]$, we have $\rho\mleft(t \, \unde_i\mright) = r - t$.
				
				\item \label{def:axiom_rank_function_2} $\rho$ is non-increasing with respect to the partial order of $\hcube{r}{\grdeg}$, that is, if $\unda \subfaceeq \undb$, then $\rho(\undb) \leq \rho(\unda)$.
				
				\item \label{def:axiom_rank_function_3} $\rho$ is supermodular.
				
				\item \label{def:axiom_rank_function_4} The values of $\rho$ are in the set $\{-1, 0, 1, \dots, r\}$.
			\end{enumerate}
			The integer $\rho(\unda)$, for $\unda \in \hcube{r}{\grdeg}$, is called the \emph{rank} of $\unda$. A hypercube $\hcube{r}{\grdeg}$ endowed with a rank function will be called a \emph{ranked hypercube}.
		\end{defi}
		
		\begin{remark} \label{rem:link_coordinate_rank}
			The above properties imply that if $\unda \in \hcube{r}{\grdeg}$ has rank $j$, then $\ssa_i \leq r - j$ for all $i \in \{1, \dots, \grdeg\}$. In particular, $\underline 0$ is the only element of rank $r$ in $\hcube{r}{\grdeg}$.
		\end{remark}
		
		The following proposition relates hypercube rank functions to matricubes, defined in~\cite{AG24}.
		
		\begin{prop} \label{prop:duality}	
			Let $\rho \colon \hcube{r}{\grdeg} \to \Z$ be an integer-valued function, and $\rk = \bar \rho$ be its conjugate. The following are equivalent:
			\begin{itemize}
				\item $\rho$ is a hypercube rank function in the sense of Definition~\ref{def:rank_function}.
				
				\item $\rk$ is the rank function of a simple matricube of rank $r$ or $r + 1$, in the sense of~\cite{AG24}.
			\end{itemize}
		\end{prop}
		
		\begin{proof}
			 Properties~\ref{def:axiom_rank_function_1}-\ref{def:axiom_rank_function_2}-\ref{def:axiom_rank_function_3} are equivalent to the axioms (R1$^*$)-(R2)-(R3), respectively, in the definition of matricubes in~\cite{AG24}. Property~\ref{def:axiom_rank_function_4} for $\rho$ is equivalent to requiring $\matricube$ being of rank $r$ or $r + 1$, via~\cite[Proposition~2.4]{AG24}.
		\end{proof}
		
		\begin{remark}
			It turns out that working with hypercube rank functions, instead of their conjugate matricube rank functions, considerably simplifies the mathematical expressions appearing in the treatment of combinatorial linear series. For this reason, we prefer them. We do not use any non-trivial result from~\cite{AG24}.
		\end{remark}
		
		Although this is not used in the following, we mention that, combining the above result with~\cite[Theorem 7.1]{AG24}, we get the following corollary.
		
		\begin{cor}
			The data of a hypercube rank function $\rho$ on $\hcube{r}{\grdeg}$ is equivalent to the data of a permutation array in the terminology of~\cite{eriksson2000combinatorial, eriksson2000decomposition}.
		\end{cor}
		
		In order to give examples in low dimension, we choose the following notational convention.
		
		\begin{conv}[Cases $\grdeg = 1, 2, 3$] \label{conv:notation_rank_function}
			In this article, for $\grdeg = 1$, a function on $\hcube{r}{1} = [r]$ is described by a tuple with $r + 1$ entries $(\sst_0, \dots, \sst_r)$, which means that the value of the function on the $i$-th entry of $\hcube{r}{1}$ is $\sst_i$.
			
			In the same way, for $\grdeg = 2$, a function on $\hcube{r}{2}$ will often be described by an array of size $(r + 1) \times (r + 1)$, $\ssub{(\sst_{ij})}!_{0 \leq i,j \leq r}$, which means that the function takes value $\sst_{ij}$ on $(i, j) \in \hcube{r}{2}$. We choose the convention that the first direction is horizontal, the second direction is vertical, and the origin is the bottom left-hand corner.
			
			When $\grdeg = 3$, a function defined on $\hcube{r}{3}$ will be specified by $r + 1$ arrays $\ssR_0, \dots, \ssR_r$ of size $(r + 1) \times (r + 1)$, where $\ssR_k$ describes the values of the function on $\hcube{r}{2} \times \{k\} \subseteq \hcube{r}{3}$.
		\end{conv}
		
		Here are two examples of two-dimensional rank functions, with $r = 3$ and $r = 4$ respectively.
		\[
			\begin{pmatrix}
				0 & 0 & 0 & 0 \\
				1 & 1 & 0 & 0 \\
				2 & 2 & 1 & 0 \\
				3 & 2 & 1 & 0
			\end{pmatrix}
			\qquad \qquad
			\begin{pmatrix}
				0 & 0 & 0 & 0 & -1 \\
				1 & 1 & 1 & 0 & -1 \\
				2 & 2 & 2 & 1 & 0 \\
				3 & 2 & 2 & 1 & 0 \\
				4 & 3 & 2 & 1 & 0
			\end{pmatrix}
		\]
	
	\subsection{Rank functions induced by complete flags} \label{subsec:rank_function_induced_flags}
		
		Let $r$ be a non-negative integer, and let $\VS$ be a vector space of dimension $r + 1$ over some field $\k$. A complete flag of $\VS$ consists of a chain of vector subspaces
		\[ \VS = \ssfilt^0 \supsetneq \ssfilt^1 \supsetneq \dots \supsetneq \ssfilt^{r - 1} \supsetneq \ssfilt^r \supsetneq \ssfilt^{r + 1} = (0), \]
		where, for each $i \in [r + 1]$, $\ssfilt^i$ is a vector subspace of codimension $i$ in $\VS$.
		
		Let $\grdeg$ be a positive integer, and let $\ssfilt_1^\bullet, \dots, \ssfilt_\grdeg^{\bullet}$ be a collection of $\grdeg$ complete flags of $\VS$. Define the function $\rho \colon \hcube{r}{\grdeg} \rightarrow \mathbb Z$ by
		\begin{equation} \label{eq:def_rank_function}
			\rho(\ssa_1, \dots, \ssa_\grdeg) \coloneqq \dim_\k \mleft(\ssfilt_{1}^{\ssa_1} \cap \dots \cap \ssfilt_{\grdeg}^{\ssa_\grdeg} \mright) - 1.
		\end{equation}
		
		\begin{prop}
			The function $\rho$ defined in~\eqref{eq:def_rank_function} is a rank function on $\hcube{r}{\grdeg}$.
		\end{prop}
		
		\begin{proof}
			This is~\cite[Proposition~2.6]{AG24}. Let $\unda$ and $\undb$ be two points of $\hcube{r}{\grdeg}$, and let $\undx \coloneqq \unda \wedge \undb$ and $\undy \coloneqq \unda \vee \undb$. We have an injection
			\[ \rquot{\bigl(\ssfilt_1^{\ssa_1} \cap \dots \cap \ssfilt_{\grdeg}^{\ssa_\grdeg}\bigr)}{\bigl(\ssfilt_1^{\ssy_1} \cap \dots \cap \ssfilt_\grdeg^{\ssy_\grdeg}\bigr)} \hookrightarrow \rquot{\bigl(\ssfilt_1^{\ssx_1} \cap \dots \cap \ssfilt_\grdeg^{\ssx_\grdeg}\bigr)}{\bigl(\ssfilt_1^{\ssb_1} \cap \dots \cap \ssfilt_\grdeg^{\ssb_\grdeg}\bigr)}, \]
			from which, comparing the dimensions, we get $\rho(\unda) - \rho\mleft(\undy\mright) \leq \rho(\undx) - \rho(\undb)$. This proves the supermodularity of $\rho$. Properties \ref{def:axiom_rank_function_1}-\ref{def:axiom_rank_function_2}-\ref{def:axiom_rank_function_4} in Definition \ref{def:rank_function} are trivially verified.
		\end{proof}
		
		A hypercube rank function $\rho$ on $\hcube{r}{\grdeg}$ is called \emph{representable over a field $\k$} if it comes from a collection of $\grdeg$ complete flags $\ssfilt_1^\bullet, \dots, \ssfilt_\grdeg^\bullet$ as above, in a $\k$-vector space $\VS$ of dimension $r + 1$.
		
		\subsubsection{Standard rank functions}
			
			The simplest kind of rank function is the following.
			
			\begin{defi}[Standard rank function] \label{def:standard_rank_function}
				The \emph{standard rank function} of dimension $\grdeg$ and width $r$ is the rank function $\rhostandard$ on $\hcube{r}{\grdeg}$ given by
				\[ \rhostandard(\unda) \coloneqq \max\{-1, r - \ssa_1 - \cdots -\ssa_\grdeg \} \qquad \forall \unda \in \hcube{r}{\grdeg}. \qedhere \]
			\end{defi}
			
			For instance, the standard rank function of dimension $2$ and width $4$ is given by
			\[
				\begin{pmatrix}
					0 & -1 & -1 & -1 & -1 \\
					1 & 0 & -1 & -1 & -1 \\
					2 & 1 & 0 & -1 & -1 \\
					3 & 2 & 1 & 0 & -1 \\
					4 & 3 & 2 & 1 & 0
				\end{pmatrix}
			\]
			
			Via Proposition~\ref{prop:duality}, a rank function $\rho$ on $\hcube{r}{\grdeg}$ is standard if, and only if, its conjugate $\rk = \bar \rho$ is the rank function of the uniform matricube $\unimatricube_{\undr, r + 1}$, see~\cite{AG24} for the definition.
			
			\begin{prop}
				Every rank function $\rho$ on the hypercube $\hcube{r}{\grdeg}$ dominates the standard rank function $\rhostandard$. That is, for every $\unda \in \hcube{r}{\grdeg}$, we have $\rho(\unda) \geq \rhostandard(\unda)$.
			\end{prop}
			
			\begin{proof}
				It will be enough to show that for every $\unda = (\ssa_1, \dots, \ssa_\grdeg) \in \hcube{r}{\grdeg}$, we have $\rho(\unda) \geq r - \ssa_1 - \cdots - \ssa_\grdeg$. This can be proved by induction on $\grdeg$, using the supermodularity of $\rho$.
			\end{proof}
			
			\begin{remark}
				The standard rank function is induced by complete flags, over an infinite field, which are in \emph{general relative position}, that is, whose intersection patterns have the smallest possible dimensions.
			\end{remark}
		
		\subsubsection{Geometric rank functions} \label{subsubsec:geometric_rank_functions}
			
			Let $C$ be a smooth proper curve over an algebraically closed field $\k$. Let $\k(C)$ be the function field of $C$, and let $\VS \subset \k(C)$ be a vector subspace of rational functions of dimension $r + 1$ over $\k$. Each $\k$-point $p$ of $C$ gives a complete flag $\ssfilt_p^\bullet$ of $\VS$ by considering the orders of vanishing at $p$ of functions in $\VS$. Define $\ssS_{p} \coloneqq \mleft\{\ord_{p}(f) \, \st \, f \in \VS \setminus \{0\}\mright\}$. This is a finite set of cardinality $r + 1$. Denote by $\sss_0^p < \dots < \sss_r^p$ its elements, enumerated in increasing order. The flag $\ssfilt_p^\bullet$ is defined by setting, for $j \in [r]$,
			\[ \ssfilt_p^j \coloneqq \mleft\{f \in \VS \setminus \{0\} \, \st \, \ord_{p}(f) \geq \sss_j^p \mright\} \cup \{0\}. \]
			Each $\ssfilt_p^j$ has codimension $j$ in $\VS$.
			
			Let now $\grdeg$ be a positive integer, and let $A = \{\ssp_1, \dots, \ssp_\grdeg\}$ be a collection of $\grdeg$ distinct $\k$-points on $C$. By the construction above, each point $\ssp_i$ leads to a complete flag $\ssfilt_{i}^\bullet$. 
			This leads to a rank function $\rho$ on the hypercube $\hcube{r}{\grdeg}$ using~\eqref{eq:def_rank_function} in Section~\ref{subsec:rank_function_induced_flags}.
			
			A rank function $\rho$ on $\hcube{r}{\grdeg}$ that arises from the above construction for a curve $C$ over an algebraically closed field $\k$ is called \emph{geometric}.
	
	\subsection{Jumps of a rank function}
		
		For $\unda \in \hcube{r}{\grdeg}$, let $\ssI_{\unda}$ be the set of all $i \in \{1, \dots, \grdeg\}$ such that $\unda + \unde_i \in \hcube{r}{\grdeg}$.
		
		\begin{prop} \label{prop:inequality_one}
			Let $\rho$ be a rank function on $\hcube{r}{\grdeg}$. For an element $\unda \in \hcube{r}{\grdeg}$ and $i\in \ssI_{\unda}$, we have $\rho(\unda) - 1 \leq \rho(\unda + \unde_i) \leq \rho(\unda)$.
		\end{prop}
		
		\begin{proof}
			The first inequality results from the supermodularity property applied to the vectors $\unda$ and $(\ssa_i + 1) \, \unde_i$, using Property~\ref{def:axiom_rank_function_1} in Definition~\ref{def:rank_function}; the second inequality comes from the non-increasing property of $\rho$.
		\end{proof}
		
		The proposition leads to the following definition.
		
		\begin{defi}[Jumps of a rank function] \label{def:jumps}
			Let $\rho$ be a rank function on $\hcube{r}{\grdeg}$. A point $\unda$ of $\hcube{r}{\grdeg}$ is called a \emph{jump} for $\rho$ if
			\begin{enumerate}
				\item \label{defi:jump_1} $\rho(\unda) \geq 0$, and
				
				\item \label{defi:jump_2} for every $i \in \{1, \dots, \grdeg\}$ such that $\unda + \unde_i$ belongs to $\hcube{r}{\grdeg}$, we have $\rho(\unda + \unde_i) = \rho(\unda) - 1$.
			\end{enumerate}
			We denote by $\ssJ_\rho$ the set of jumps of $\rho$.
		\end{defi}
		
		Here are three rank functions of dimension two, with $r = 3$ for the first one, and $r = 4$ for the second and third ones. The jumps of each rank function are depicted in blue.
		\[
			\begin{pmatrix}
				0 & 0 & 0 & \color{blue} 0 \\
				1 & \color{blue} 1 & 0 & 0 \\
				2 & \color{blue} 2 & \color{blue} 1 & 0 \\
				\color{blue} 3 & 2 & 1 & 0
			\end{pmatrix}
			\qquad \qquad
			\begin{pmatrix}
				0 & 0 & 0 & \color{blue} 0 & -1 \\
				1 & 1 & \color{blue} 1 & 0 & -1 \\
				2 & 2 & \color{blue} 2 & \color{blue} 1 & \color{blue} 0 \\
				\color{blue} 3 & 2 & 2 & 1 & 0 \\
				\color{blue} 4 & \color{blue} 3 & 2 & 1 & 0
			\end{pmatrix}
			\qquad \qquad
			\begin{pmatrix}
				\color{blue} 0 & -1 & -1 & -1 & -1 \\
				\color{blue} 1 & \color{blue} 0 & -1 & -1 & -1 \\
				\color{blue} 2 & \color{blue} 1 & \color{blue} 0 & -1 & -1 \\
				\color{blue} 3 & \color{blue} 2 & \color{blue} 1 & \color{blue} 0 & -1 \\
				\color{blue} 4 & \color{blue} 3 & \color{blue} 2 & \color{blue} 1 & \color{blue} 0
			\end{pmatrix}
		\]
		
		The set of jumps of a rank function has the following properties.
		
		\begin{prop} \label{prop:jump_stability} 
			The set of jumps $\ssJ_\rho$ of a rank function $\rho$ on $\hcube{r}{\grdeg}$ is stable under $\wedge$.
		\end{prop}
		
		\begin{proof}
			Let $\unda, \undb \in \ssJ_\rho$ and let $\undc = \unda \wedge \undb$. The non-increasing property of $\rho$ ensures that $\rho(\undc) \geq 0$. Let $i \in \{1, \dots, \grdeg\}$ be such that $\undc + \unde_i$ belongs to $\hcube{\undr}{}$. We have to show that $\rk(\undc + \unde_i) = \rk(\undc) + 1$. By symmetry, we can suppose that $\ssa_i \leq \ssb_i$, that is, $\ssc_i = \ssa_i$. Since $\unda$ is a jump, we have $\rk(\unda + \unde_i) = \rk(\unda) + 1$. We conclude by applying the supermodularity property to the points $\unda + \unde_i$ and $\undc$.
		\end{proof}
		
		\begin{prop} \label{prop:jumps_graded_poset}
			The set of jumps $\ssJ_\rho$ is a graded poset. The grading is given by the conjugate $\bar \rho$. In particular, if $\unda \subface \undb$ are two distinct jumps of $\rho$, then we have $\rho(\unda) > \rho(\undb)$.
		\end{prop}
		
		\begin{proof}
			It is sufficient to show that if $\unda \subface \undb$ are two distinct jumps of $\rho$ with $\rho(\undb) \leq \rho(\unda) - 2$, then there exists a jump $\unda \subface \undc \subface \undb$ such that $\rho(\undc) = \rho(\unda) - 1$. Let $i$ be an index such that $\ssa_i < \ssb_i$, and let $\undx \coloneqq \unda + \unde_i$. Since $\unda$ is a jump, we have $\rho(\undx) = \rho(\unda) - 1$. Let $\undc$ be the meet of all the points $\undc'$ of $\hcube{r}{\grdeg}$ verifying $\undc' \supfaceeq \undx$ and such that either $\undc'$ is a jump or $\undc' = \undr$. Since $\undx \subfaceeq \undb$ and $\undb$ is a jump, $\undc$ is also a jump by virtue of Proposition~\ref{prop:jump_stability}, and moreover $\unda \subface \undc \subfaceeq \undb$. It is easy to see, by induction on the $\ell_1$-norm of $\undx$, that $\rho(\undc) = \rho(\undx) = \rho(\unda) - 1$, which implies that $\undc \subface \undb$. This shows that $\undc$ is as desired and that the poset $\ssJ_\rho$ is graded. By construction, the grading is given by $\bar\rho$.
		\end{proof}
		
		The following fact will be useful in the sequel.
		
		\begin{fact} \label{fact:max_only_jump}
			Let $\rho$ be a rank function on $\hcube{r}{\grdeg}$. If $\undr$ is a jump for $\rho$, then it is the only element of $\ssJ_\rho$ having some coordinate equal to $r$.
		\end{fact}
		
		\begin{proof}
			For the sake of a contradiction, let $\undx \neq \undr$ be an element of $\ssJ_\rho$ with some coordinate equal to $r$, say the first one. By~\ref{defi:jump_1} in Definition~\ref{def:jumps}, we have $\rho(\undr) \geq 0$. The inequality $\undx \supfaceeq r \, \unde_1$ implies that $\rho(\undx) \leq 0$, and Proposition~\ref{prop:jumps_graded_poset} implies $\rho(\undr) < \rho(\undx) \leq 0$, a contradiction.
		\end{proof}
		
		\begin{remark}
			Let $\matricube$ be the matricube on the ground set $\hcube{r}{\grdeg}$ defined by the conjugate $\rk = \bar \rho$. The set of jumps $\ssJ_\rho$ coincides with the set of flats $\Fl(\matricube)$ of $\matricube$ if $\rk(\matricube) = r$ and with $\Fl(\matricube) \setminus \{\undr\}$ if $\rk(\matricube) = r + 1$. Propositions~\ref{prop:jump_stability} and~\ref{prop:jumps_graded_poset} can be deduced from the analogue properties of flats of matricubes. For the sake of completeness, we provided the short proofs of these results. Note that contrary to the set of flats of matricubes, the set of jumps is not necessary a lattice (this happens only in the case $\rk(\matricube) = r + 1$).
	\end{remark}
	
	\subsection{Partition Lemma} \label{subsec:partition_lemma}
		
		In this section, we prove a result about hypercube rank functions which turns out to be useful in the sequel.
		
		Let $\rho$ be a rank function on $\hcube{r}{\grdeg}$. The point $\underline 0$ is the only point of $\hcube{r}{\grdeg}$ whose image by $\rho$ is $r$ (Remark~\ref{rem:link_coordinate_rank}). Besides, the set of jumps $\ssJ_\rho$ of $\rho$ contains the point $\underline 0$ (because $\rho(\unde_i) = r - 1$ for all $i$). Every jump of $\rho$ of rank $r - 1$ has only coordinates equal to zero or one (Remark~\ref{rem:link_coordinate_rank}), among which at least one is equal to one. For each $\unda \in \ssJ_\rho$ such that $\rho(\unda) = r - 1$, denote by $\ssP{\unda}$ the subset of $\{1, \dots, \grdeg\}$ consisting of all the indices $i$ with $\ssa_i = 1$, that is, $\ssP{\unda}$ is the \emph{support} of $\unda$. Denote by $\ssmathcalP_\rho$ the collection of all sets $\ssP{\unda}$ for $\unda \in \ssJ_\rho \setminus \{\underline 0\}$ verifying $\rho(\unda) = r - 1$. We have the following proposition.
		
		\begin{lemma}[Partition Lemma] \label{lem:partition}
			Notation as above, $\ssmathcalP_{\rho}$ provides a partition of $\{1, \dots, \grdeg\}$.
		\end{lemma}
		
		\begin{proof}
			We first prove that the elements of $\ssmathcalP_\rho$ are pairwise disjoint. Let $\unda$ and $\undb$ be two distinct elements of $\ssJ_\rho \setminus \{\underline 0\}$ with $\rho(\unda) = \rho(\undb) = r - 1$. Since $\unda \not \subfaceeq \undb$ and $\undb \not \subfaceeq \unda$, we have $\rho(\unda \vee \undb) \leq r - 2$. Using the supermodularity property for $\unda$ and $\undb$, we get $\rho(\unda \wedge \undb) \geq r$ and therefore $\rho(\unda \wedge \undb) = r$. This forces $\unda \wedge \undb = \underline 0$, from which we can conclude that $\ssP{\unda} \cap \ssP{\undb} = \varnothing$.
			
			It remains to prove that the sets $\ssP{\unda}$ cover $\{1, \dots, \grdeg\}$. For an $i \in \{1, \dots, \grdeg\}$, we need to show the existence of $\unda \in \ssJ_\rho$ with $\rho(\unda) = r - 1$ and $\ssa_i = 1$. We define $\unda$ as the meet of all the jumps $\unda' \in \ssJ_\rho$ such that $\unda' \supfaceeq \unde_i$. As in the proof of Proposition~\ref{prop:jumps_graded_poset}, $\unda$ is well-defined, belongs to $\ssJ_\rho$ and has rank $\rk(\unda) = \rk(\unde_i) = r - 1$. The non-increasing property of $\rho$ implies that $\ssa_i \leq 1$, as desired.
		\end{proof}
		
		\begin{remark} \label{rem:partition_1}
			In the case $r = 1$, in the construction above, the condition $\rho(\unda) = r - 1$ is automatic. This will be crucially used in Section~\ref{sec:g1d}.
		\end{remark}
		
		Other elementary results about rank functions will be given in Section~\ref{subsec:realization_jumps}.
		
\section{Slope structures} \label{sec:slope_structure}
	
	In the following sections, we define combinatorial linear series on metric graphs with the help of an auxiliary data called a \emph{slope structure}. A slope structure is the data of a family of hypercube rank functions of given width $r$, parametrized by the points of the metric graph, of varying dimension given by the valences of points, and verifying a finiteness condition.
	
	\subsection{Slope structures on graphs} \label{subsec:slope_structures}
		
		Let first $G = (V, E)$ be a simple graph. We denote by $\mathbb E$ the set of all the orientations of edges of $G$, so that for an edge $\{u, v\}$ in $E$, we have two orientations $uv, vu \in \mathbb E$. For an oriented edge $e = uv \in \mathbb E$, we call $u$ the tail and $v$ the head of $e$. We denote by $\overline e = vu$ the oriented edge in $\mathbb E$ with reverse orientation. For a vertex $v \in V$, we denote by $\mathbb E_v \subseteq \mathbb E$ the set of oriented edges in $\mathbb E$ which have tail $v$, that is, all $vu \in \mathbb E$ for edges $\{v, u\} \in E$.
		
		A \emph{slope structure} $\fS = \mleft\{\ssS^v; \ssS^e\mright\}_{v \in V, e \in \mathbb E}$ \emph{of width $r$} on $G$, or simply an \emph{$r$-slope structure}, is the data of
		
		\begin{enumerate}[label=(SLS\arabic*)]
			\item For every oriented edge $e = uv \in \mathbb E$ of $G$, a collection $S^e$ of $r + 1$ integers
			\[
				\sss_0^e < \sss_1^e < \dots < \sss_r^e
			\]
			subject to the requirement that $\sss_i^{uv} + \sss_{r - i}^{vu} = 0$ for every edge $\{u, v\} \in E$.
			
			\item For every vertex $v$ of $G$, a rank function $\ssrho_v$ on the hypercube $\hcube{r}{\ssd_v}$.
		\end{enumerate}
		
		If $\ssJ_{\ssrho_v}$ denotes the set of jumps of $\ssrho_v$ (see Definition~\ref{def:jumps}), we denote by $\ssS^v \subseteq \prod_{e \in \mathbb E_v} S^e$ the set of all points $\sss_{\unda}$ for $\unda \in \ssJ_{\ssrho_v}$.
		
		Here, for a point $\unda = (\ssa_e)_{e \in \mathbb E_v}$ of the hypercube, the element $\sss_{\unda} \in \prod_{e \in \mathbb E_v} \ssS^e$ denotes the point in the product which has coordinate at $e \in \mathbb E_v$ equal to $\sss_{\ssa_e}^e$. In other words, $\ssS^v$ fits into the following natural commutative diagram:
		\[
			\begin{tikzcd}
				\ssJ_{\ssrho_v} \arrow[d, "\unda \mapsto \sss_{\unda}"] \arrow[r, hook] & \hcube{r}{\ssd_v} \arrow[d, "\unda \mapsto \sss_{\unda}"] \\
				\ssS^v \arrow[r, hook] & \Pi_{e \in \mathbb E_v} \ssS^e
			\end{tikzcd}
		\]
		
		We will sometimes need to separate the data relative to edges and the data relative to vertices. In this case, we will denote by $\fS^\ee$ the data of a set of prescribed slopes on each edge, and by $\fS^\vv$ the data of a rank function for each vertex.
	
	\subsection{Slope structures on metric graphs} \label{subsec:slope_structure_metric_graphs}
		
		Let now $\Gamma$ be a metric graph. By an \emph{$r$-slope structure} on $\Gamma$ we mean an $r$-slope structure $\fS$ on a simple graph model $G = (V, E)$ of $\Gamma$, extended to each point of $\Gamma$ as follows.
		
		For every point $x$ and each outgoing unit tangent vector $\nu \in \T_x(\Gamma)$, there exists a unique oriented edge $uv$ of $G$ which is parallel to $\nu$. Define $\ssS^\nu = \ssS^{uv}$. Also, for every point $x \in \Gamma \setminus V$ in the interior of an edge $\{u, v\}$, define $\ssrho_x$ to be the standard rank function on $\hcube{r}{2}$. In particular, $\ssS^x \subseteq \ssS^{uv} \times \ssS^{vu}$ can be identified with the set of all pairs $(\sss_i^{uv}, \sss_j^{vu})$ with $i + j \leq r$. We call the collection $\mleft\{\ssS^x; \ssS^{\nu} \, \st x \in \Gamma, \, \nu \in \T_x(\Gamma)\mright\}$ a \emph{slope structure of width $r$}, or simply an \emph{$r$-slope structure} on $\Gamma$. We denote it by $\fS_\Gamma$, or simply $\fS$, if there is no risk of confusion. We extend the notation $\fS^\ee$ and $\fS^\vv$ in the natural way. Note that a slope structure on a metric graph can arise from choices of slope structures on different graph models of $\Gamma$.
		
		\begin{example} \label{ex:slope_structure}
			We give an example of a 2-slope structure on a circle. Consider the metric graph $\Gamma$ depicted below, with edges of arbitrary positive lengths.
			
			\begin{figure}[h!]
				\centering
				\begin{tikzpicture}[scale=0.7]
					\coordinate (A) at (0,0);
					\coordinate (B) at (3,0);
					
					\draw[postaction=decorate,decoration={markings,mark=at position 1/2 with {\arrow[red]{>}; }}]
					(A) .. controls (1,0.8) and (2,0.8) .. (B);
					
					\draw[postaction=decorate,decoration={markings,mark=at position 1/2 with {\arrow[red]{>};}}]
					(A) -- (B);
					
					\foreach \c in {A, B} {
						\fill (\c) circle (2.5pt);
						}
					
					\node[left] at (A) {$u$};
					\node[right] at (B) {$v$};
				\end{tikzpicture}
			\end{figure}
			
			Let $\sss_0 < \sss_1 < \sss_2$ and $\sss_0' < \sss_1' < \sss_2'$ be two sets of distinct integers. We define $\fS^\ee$ by allowing slopes $\sss_0 < \sss_1 < \sss_2$ on the top edge and $\sss_0' < \sss_1' < \sss_2'$ on the bottom edge. We define $\fS^\vv$ by choosing the rank functions at $u$ and $v$ to be given by the array
			\[ \begin{pmatrix}
				0 & 0 & \color{blue} 0 \\
				1 & \color{blue} 1 & 0 \\
				\color{blue} 2 & 1 & 0
			\end{pmatrix} \]
			(jumps depicted in blue). The rank functions at all the other points of $\Gamma$ are standard. This fully describes a $2$-slope structure $\fS$ on $\Gamma$.
		\end{example}
	
	\subsection{Rational functions compatible with a slope structure} \label{subsec:rational_functions_compatible}
		
		We define a notion of rational function on metric graphs compatible with a slope structure.
		
		Let $\Gamma$ be a metric graph and let $\fS = \mleft\{\ssS^x; \ssS^\nu \st {x \in \Gamma, \, \nu \in \T_x(\Gamma)}\mright\} = \mleft(\fS^\ee, \fS^\vv\mright)$ be a slope structure of width $r$ on $\Gamma$. Recall that we denote by $\Rat(\Gamma)$ the set of continuous piecewise affine functions $f \colon \Gamma \rightarrow \mathbb R$ with integral slopes. For each point $x \in \Gamma$, and each $\nu \in \T_x(\Gamma)$, we denote by $\slope_{\nu}(f)$ the slope of $f$ at $x$ along $\nu$.
		
		A function $f$ in $\Rat(\Gamma)$ is said to be \emph{compatible} with $\fS$ if the two conditions $(i)$ and $(ii)$ below are verified:
		
		\begin{itemize}
			\item[$(i)$] for each point $x \in \Gamma$ and each $\nu \in \T_x(\Gamma)$, the outgoing slope of $f$ along $\nu$ lies in $\ssS^\nu$.
		\end{itemize}

		Denote by 
		\[ \ssdelta_x(f) \in \prod_{\nu \in \T_x(\Gamma)} \ssS^\nu, \qquad \ssdelta_x(f)(\nu) = \slope_\nu(f) \quad \forall\, \nu \in \T_x(\Gamma) \]
		the vector in $\prod_{\nu \in \T_x(\Gamma)} \ssS^\nu$ with $\nu$-coordinate consisting of the outgoing slope of $f$ along $\nu \in \T_x(\Gamma)$. Then, the second condition is:
		
		\begin{itemize}
			\item[$(ii)$] for every point $x \in \Gamma$, the vector $\ssdelta_x(f)$ belongs to $\ssS^x$.
		\end{itemize}
		
		We denote by $\Rat(\Gamma, \fS)$, or simply $\Rat(\fS)$ if there is no risk of confusion, the space of rational functions on $\Gamma$ compatible with $\fS$. We also denote by $\Rat\mleft(\Gamma, \fS^\ee\mright)$ or $\Rat\mleft(\fS^\ee\mright)$ the space of rational functions satisfying $(i)$.
		
		If $\Gamma$ is $\Lambda$-rational, we define the spaces $\Rat_\Lambda(\Gamma, \fS)$, $\Rat_\Lambda(\Gamma, \fS^\ee)$ and $\Rat_\Lambda(\Gamma, \fS^\vv)$ accordingly, adding the constraint that $f(x)$ is in $\Lambda$ for all $\Lambda$-rational points of $\Gamma$.
	
	\subsection{Slope substructures} \label{subsec:slope_substructures}
		
		Let $\Gamma$ be a metric graph and let $\fS = \mleft\{\ssS^x; \ssS^\nu \st {x \in \Gamma, \, \nu \in \T_x(\Gamma)}\mright\}$ be an $r$-slope structure on $\Gamma$.
		
		A \emph{slope substructure} of $\fS$ of width $s \leq r$ is a slope structure
		\[ \fS' = \mleft\{\ssS'^{x}; \ssS'^{\nu} \st {x \in \Gamma, \, \nu \in \T_x(\Gamma)}\mright\} \]
		of width $s$ on $\Gamma$ such that for every $x \in \Gamma$, we have:
		
		\begin{enumerate}
			\item for every $\nu \in \T_x(\Gamma)$, the set of prescribed slopes $\ssS^{'\nu}$ is a subset of $\ssS^\nu$;
			
			\item the set of prescribed vectors of slopes $\ssS^{'x}$ is a subset of $\ssS^x$.
		\end{enumerate}
		
		Note that if $\fS'$ is a slope substructure of $\fS$, then the inclusion $\Rat(\fS') \subseteq \Rat(\fS)$ holds. We will see in Corollary~\ref{cor:equivalence_slope_substructure} that in the case of interest to us, the converse will be also true.
	
	\subsection{Divisors on a metric graph and their rank} \label{subsec:divisors}
		
		A divisor $D$ on a metric graph $\Gamma$ is a finite formal sum over $\Z$ of points of $\Gamma$, that is, $D = \sum_{i \in I} \ssn_i \, (\ssx_i)$ with $\ssn_i \in \Z$ and distinct points $\ssx_i \in \Gamma$, for a finite set $I$. The coefficient of a point $x$ of $\Gamma$ in $D$ is denoted by $D(x)$.
		A divisor $D$ is called \emph{effective}, written $D \geq 0$, if $D(x) \geq 0$ for all $x \in \Gamma$. For any rational function $f \in \Rat(\Gamma)$, the corresponding divisor is denoted by
		\[ \div(f) \coloneqq \sum_{x \in \Gamma} \ord_x(f) \, (x), \qquad \textrm{where } \ord_x(f) \coloneqq - \sum_{\nu \in \T_x(\Gamma)} \slope_{\nu}(f). \]
		A divisor obtained in this way is called \emph{principal}. In the case $\Gamma$ is $\Lambda$-rational, a divisor whose support is made up of $\Lambda$-rational points is said to be a \emph{$\Lambda$-rational} divisor. Notice that the space $\Rat_\Lambda(\Gamma)$ defined in Section~\ref{subsec:basic_notions} can be redefined as the set of functions of $\Rat(\Gamma)$ such that $\div(f)$ is $\Lambda$-rational.
		
		We have the following elementary fact.
		
		\begin{prop} \label{prop:functions_equal_divisor}
			For $f, g \in \Rat(\Gamma)$, $\div(f) = \div(g)$ if, and only if, $f - g$ is constant on $\Gamma$.
		\end{prop}
		
		Two divisors $D_1$ and $D_2$ are called \emph{linearly equivalent} if their difference $\ssD_1 - \ssD_2$ is principal. The Baker--Norine rank $r(D)$ of a divisor $D$ is defined as the maximum integer among $-1$ and the integers $r \geq 0$ such that for all points $\ssx_1, \dots, \ssx_r$ in $\Gamma$, the divisor $D - (\ssx_1) - \cdots - (\ssx_r)$ is linearly equivalent to an effective divisor.
	
	\subsection{Linear equivalence of slope structures}
		
		We define a notion of linear equivalence for slope structures on a metric graph as follows.
		
		Let $\fS_1 = \mleft\{\ssS_1^x ; \ssS^\nu_1 \st x \in \Gamma, \, \nu \in \T_x(\Gamma)\mright\}$ and $\fS_2 = \mleft\{\ssS_2^x ; \ssS^\nu_2 \st x \in \Gamma, \, \nu \in \T_x(\Gamma)\mright\}$ be two slope structures on a metric graph $\Gamma$. We say $\fS_1$ and $\fS_2$ are \emph{linearly equivalent}, and write $\fS_1 \simeq \fS_2$, if there exists a rational function $f$ on $\Gamma$ such that for every point $x$ of $\Gamma$ and every $\nu \in \T_x(\Gamma)$, we have $\ssS^\nu_{1} = \ssS^\nu_{2} - \slope_\nu(f)$, and $\ssS^x_1 = \ssS^x_2 - \ssdelta_x(f)$. In this case, we write $\fS_1 = \fS_2 + \div(f)$. Note that if $\fS$ is a slope structure, then $\fS + \div(f)$ is a slope structure for every rational function $f \in \Rat(\Gamma)$.
	
	\subsection{Divisors endowed with a slope structure on $\Gamma$} \label{subsec:divisors_endowed_slope_structure}
		
		A \emph{divisor endowed with an $r$-slope structure of degree $d$} is a pair $(D, \fS)$ consisting of a divisor $D$ of degree $d$ and a slope structure $\fS$ of width $r$. We extend the definition of linear equivalence between slope structures to all pairs $(D, \fS)$ with $D$ a divisor of degree $d$ and $\fS$ an $r$-slope structure on $\Gamma$ by declaring that $(\ssD_1,\fS_1) \simeq (\ssD_2, \fS_2)$ if there exists a rational function $f$ on $\Gamma$ such that $\ssD_1 = \ssD_2 + \div(f)$ and $\fS_1 = \fS_2 + \div(f)$.
		
		\begin{defi} \label{def:divisor_endowed_slope_structure}
			A \emph{divisor class endowed with an $r$-slope structure of degree $d$} on $\Gamma$ is the linear equivalence class of a pair $(D, \fS)$ where $D$ is a divisor of degree $d$ and $\fS$ is an $r$-slope structure on $\Gamma$.
		\end{defi}
		
		We now define the space of rational functions relative to a divisor and a slope structure.
		
		\begin{defi}[Space of rational functions and linear system associated to a divisor endowed with a slope structure] \label{def:space_rational_functions}
			Let $(D, \fS)$ be a divisor endowed with a slope structure on $\Gamma$. We denote by $\Rat(D, \fS)$ the space of all $f \in \Rat(\fS)$ with the property that $D + \div(f) \geq 0$, and define the linear system $|(D, \fS)|$ associated to $(D, \fS)$ as the space of all effective divisors $E$ on $\Gamma$ of the form $D + \div(f)$ for some $f \in \Rat(D, \fS)$.
		\end{defi}
		
		\begin{remark}
			Recall that we define $\Rat(D)$ as the set of all functions $f \in \Rat(\Gamma)$ such that $D + \div(f) \geq 0$. In a similar way, we can define the space $\Rat\mleft(D, \fS^\ee\mright)$, see Section~\ref{subsec:rational_functions_compatible}.
		\end{remark}
		
		\begin{remark} \label{rem:linear_series_independence}
			Note that $|(D, \fS)|$ is independent of the choice of the pair $(D, \fS)$ in its linear equivalence class. Also note that if $D(x) > 0$ for some $x \in \Gamma$ in the interior of an edge of a model $G$ on which $\fS$ is defined, then, we have
			\[ |(D, \fS)| = |(D - (x), \fS)| + (x). \qedhere\]
		\end{remark}
		
		\begin{defi} \label{def:effective_divisor}
			A divisor endowed with a slope structure $(D, \fS)$ is called \emph{effective} if $\Rat(D, \fS)$ contains the null function.
		\end{defi}
		
		This is equivalent to asking that $D$ is effective and that we have $0 \in S^\nu$ and $\underline 0 \in S^x$ for every point $x$ and every $\nu \in \T_x(\Gamma)$.
		
		For future use, we make the following remark.
		
		\begin{remark} \label{rem:obvious_bound_slopes}
			Since $\Rat(D, \fS) \subset \mathcal{C}^0(\Gamma, \R)$, it is naturally endowed with the norm $\| \cdot \|_\infty$. The corresponding topology shall be used later on to study linear series (see Definition~\ref{def:grd}). We note that the slopes of all functions in $\Rat(D, \fS)$ are trivially bounded in magnitude by $\max_{1 \leq i \leq r} |\sss_i^e|$.
		\end{remark}

\section{Crude linear series} \label{sec:crude_linear_serues}
	
	We define \emph{crude linear series} which are the simplest notion of combinatorial linear series. The requirement in the definition is reminiscent of the rank condition on divisors on metric graphs. Moreover, it takes into account the data of the hypercube rank function on points.
	
	\subsection{Definition} \label{subsec:def_crude_linear_series}
		
		A \emph{crude linear series of degree $d$ and rank $r$}, or \emph{crude $\g^r_d$}, is the equivalence class of a divisor $D$ of degree $d$ endowed with a slope structure $\fS$ of width $r$ on $\Gamma$ subject to the following property:
		
		\noindent $(*)$ For every effective divisor $E$ on $\Gamma$ of degree $r$, there exists a rational function $f \in \Rat(D, \fS)$ such that
		
		\begin{enumerate}[label=(CL\arabic*)]
			\item \label{def:crude_linear_series_1} For every point $x \in \Gamma$, $\ssrho_x(\ssdelta_x(f)) \geq E(x)$, and in addition,
			
			\item \label{def:crude_linear_series_2} $D + \div(f) - E \geq 0$. \qedhere
		\end{enumerate}
		
		We call a crude linear series \emph{effective} if the underlying divisor endowed with the slope structure is so. \hfill $\diamond$
		
		\medskip
		
		We make a set of comments and provide examples in order to clarify the definition.
		
		First, note that Property~\ref{def:crude_linear_series_2} in Definition~\ref{subsec:def_crude_linear_series} implies that the Baker--Norine rank $r(D)$ is greater than or equal to $r$.
		
		Second, for given $E$ and $f$, \ref{def:crude_linear_series_1} does not necessarily imply \ref{def:crude_linear_series_2}, for example at points $x$ such that $D(x) < 0$ or such that $\ssrho_x$ is not standard. On the other hand, \ref{def:crude_linear_series_1} implies \ref{def:crude_linear_series_2} generically for all $E$ and $f$. More precisely, let $G = (V, E)$ be a graph model of $\Gamma$ such that $\fS$ comes from an $r$-slope structure on $G$, and $D$ has support on $V$. Then, for every point $x \in \Gamma \setminus V$ lying on an edge $\{u, v\}$, the rank function $\ssrho_x$ is standard. The first condition above is thus equivalent to $i + j \leq r - E(x)$ for $\sss_i^{uv}$ and $\sss_j^{vu}$ the two slopes of $f$ at $x$. In particular, since the slopes are all integral, it is easy to see that condition \ref{def:crude_linear_series_2}, written $\div_x(f) = - \sss_i^{uv} - \sss_j^{vu} \geq E(x)$, is automatically implied by condition \ref{def:crude_linear_series_1} for $x$. This means that \ref{def:crude_linear_series_1} implies \ref{def:crude_linear_series_2} outside the (finite) set of vertices of $G$. Note that \ref{def:crude_linear_series_1} can be \emph{strictly} stronger than \ref{def:crude_linear_series_2} in the interior of edges $e$ (as long as the possible slopes on $e$ do not form an integral interval, i.e., in the case there are gaps in $\ssS^e$.
		
		The relevance of \ref{def:crude_linear_series_1} will be justified in Section~\ref{sec:geometric_tropicalization}, which treats the geometric situation in which the slope structure comes from tropicalization.
		
		We finally note that the definition generalizes to $\Lambda$-rational divisors on $\Lambda$-metric graphs. In this case, we require $f$ to be in $\Rat_\Lambda(\Gamma)$.
		
		\begin{example} \label{ex:first_grd}
			Consider the metric graph $\Gamma$ below, with two edges of equal length this time.
			\begin{figure}[h!]
				\centering
				\begin{tikzpicture}[scale=0.7]
					\coordinate (A) at (0,0);
					\coordinate (B) at (3,0);
					
					\draw[postaction=decorate,decoration={markings,mark=at position 1/2 with {\arrow[red]{>}; }}]
					(A) .. controls (1,0.8) and (2,0.8) .. (B);
					
					\draw[postaction=decorate,decoration={markings,mark=at position 1/2 with {\arrow[red]{>};}}]
					(A) .. controls (1,-0.8) and (2,-0.8) .. (B);
					
					\foreach \c in {A, B} {
						\fill (\c) circle (2.5pt);
						}
					
					\node[left] at (A) {$u$};
					\node[right] at (B) {$v$};
				\end{tikzpicture}
			\end{figure}
			
			We allow the slopes $0 < 1 < 2$ on both edges in the direction of the arrows. We define the rank functions at $u$ and $v$ by the array
			\[ \begin{pmatrix}
				0 & 0 & \color{blue} 0 \\
				1 & \color{blue} 1 & 0 \\
				\color{blue} 2 & 1 & 0
			\end{pmatrix} \]
			(jumps depicted in blue), and the rank functions at all the other points of $\Gamma$ are chosen to be standard. This fully describes $\fS$. Consider the effective divisor $D = 4 \, (u)$. The pair $(D, \fS)$ is an effective crude linear series of degree $4$ and rank $2$. To see this, we need to check Property $(*)$ in Definition~\ref{subsec:def_crude_linear_series} for $E = (x) + (y)$, for points $x, y$ in $\Gamma$. This can be done by a case analysis depending on whether $x$ or $y$ coincide with a vertex, or they are in the interior of the same edge, or in the interior of two different edges of $\Gamma$. For example,
			\begin{itemize}
				\item if $x, y \notin \{u, v\}$ and $x, y$ are on the same edge, then we can take $f$ to behave as the following function (with slopes $2$, $1$ and $0$) on both edges:
				
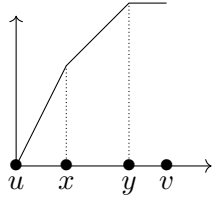
\begin{center}
					\begin{tikzpicture}[scale=0.7]
						\draw[->] (0,0) -- (2.6,0);
						\draw[->] (0,0) -- (0,2);
						\draw (0,0) node[below]{$u$} node{$\bullet$};
						\draw (2/3,0) node[below]{$x$} node{$\bullet$};
						\draw (6/4,0) node[below]{$y$} node{$\bullet$};
						\draw (2,0) node[below]{$v$} node{$\bullet$};
						\draw (0,0) -- (2/3,4/3) -- (6/4,26/12) -- (2,26/12);
						\draw[densely dotted] (2/3,0) -- (2/3,4/3);
						\draw[densely dotted] (6/4,0) -- (6/4,26/12);
					\end{tikzpicture}
				\end{center}
				
				\item or, if $x, y \notin \{u, v\}$ and $x, y$ are on different edges, we define $x'$ (resp. $y'$) to be the point of the edge containing $x$ (resp. $y$) symmetrical about the middle of this edge. Then, we can take $f$ to be the function			
				
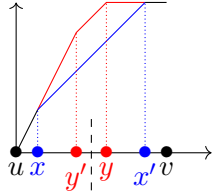
\begin{center}
					\begin{tikzpicture}[scale=0.9]
						\draw[->] (0,0) -- (2.6,0);
						\draw[->] (0,0) -- (0,2);
						\draw (0,0) node[below]{$u$} node{$\bullet$};
						\draw (2/7,0) node[below,blue]{$x$} node[color=blue]{$\bullet$};
						\draw (4/5,0) node[below,red]{$y'$} node[red]{$\bullet$};
						\draw (6/5,0) node[below,red]{$y$} node[red]{$\bullet$};
						\draw (12/7,0) node[below,blue]{$x'$} node[blue]{$\bullet$};
						\draw (2,0) node[below]{$v$} node{$\bullet$};
						\draw (0,0) -- (2/7,4/7);
						\draw[red] (2/7,4/7) -- (4/5,8/5) -- (6/5,2) -- (12/7,2);
						\draw[blue] (2/7,4/7) -- (12/7,2);
						\draw (12/7,2) -- (2,2);
						\draw[densely dotted,blue] (2/7,0) -- (2/7,4/7);
						\draw[densely dotted,red] (4/5,0) -- (4/5,8/5);
						\draw[densely dotted,red] (6/5,0) -- (6/5,2);
						\draw[densely dotted,blue] (12/7,0) -- (12/7,2);
						\draw[dashed] (1,-1/2) -- (1,1/2);
					\end{tikzpicture}
				\end{center}
				with values on one edge represented in blue, values on the other edge represented in red, and values common to both edges represented in black.
				
				\item We omit the remaining cases. 
			\end{itemize}

			On the same graph, $D = 2 \, (u) + 2\, (v)$ with slopes $-1 < 0 < 1$ on both edges, and the same slope structure as above, provides another crude linear series of degree $4$ and rank $2$.
		\end{example}
		
		\begin{example} \label{ex:second_grd}
			Here is another simple example that will be used later (see Example~\ref{ex:finiteness_condition}). We consider the metric path graph $\Gamma$ with three vertices depicted in the figure below, with arbitrary positive edge lengths.
			
			\begin{figure}[h!]
				\centering
				\begin{tikzpicture}
					\coordinate (A) at (0,0);
					\coordinate (B) at (1.5,0);
					\coordinate (C) at (4,0);
					
					\draw[postaction=decorate,decoration={markings,mark=at position 1/2 with {\arrow[red]{>}; }}]
					(A) -- (B);
					
					\draw[postaction=decorate,decoration={markings,mark=at position 1/2 with {\arrow[red]{>}; }}]
					(B) -- (C);
					
					\foreach \c in {A, B, C} {
						\fill (\c) circle (2.5pt);
						}
					
					\node[left] at (A) {$u$};
					\node[below] at (B) {$v$};
					\node[right] at (C) {$w$};
				\end{tikzpicture}
			\end{figure}
			
			We allow slopes $-2 < 0 < 2$ in the direction of the arrows and take $D = 2\, (u) + 4 \, (v) + 2 \, (w)$. We choose the rank function $\ssrho_v$ on the vertex $v$ to be the same as in Example~\ref{ex:first_grd}. The rank functions on $u$ and $w$ are automatically standard, since these vertices are of valence 1. Then, $(D, \fS)$ is a crude linear series of degree $8$ and rank $2$.
			
			Without changing $D$, we can also consider a slope substructure $\fS'$ of $\fS$ of rank one by allowing slopes $0 < 2$ on the edge $uv$ and slopes $-2 < 0$ on the edge $vw$. We adapt $\ssrho_v'$ at $v$ as follows (jumps depicted in blue):
			$\begin{pmatrix}
				0 & \color{blue} 0 \\
				\color{blue} 1 & 0
			\end{pmatrix}.$
			The pair $(D, \fS')$ is an effective crude linear series of degree $8$ and rank $1$.
			
			We divide the coefficients of $D$ to obtain the divisor $D' = 2 \, (v) + (u) + (w)$. We then consider the slope substructure $\fS''$ with allowed slopes $0 < 1$ on $uv$ and $-1 < 0$ on $vw$. This makes $(D', \fS'')$ an effective crude linear series of degree $4$ and rank $1$.
		\end{example}
	
	\subsection{Non-increasing property of slope vectors and a finiteness theorem} \label{subsec:slope_vectors_non_increasing}
		
		In the rest of this section, we prove two important results about crude linear series.
		
		Let $(D, \fS)$ be a crude linear series on $\Gamma$ of rank $r$. Let $G = (V, E)$ be a model of $\Gamma$ such that $\fS$ is defined on $G$ and $D$ is supported on $V$. Let $e = uv$ be an oriented edge of $G$. For each point $x$ in the interior of $e$, let $\nu \in \T_x(\Gamma)$ be the unit tangent vector consistent with the orientation of $e$. Let $\sss_0^\nu(x) < \sss_1^\nu(x) < \dots < \sss_r^\nu(x) $ be the corresponding slopes in $\ssS^e$. By an abuse of notation, for each point $y$ on $e$, we still denote by $\nu$ the tangent vector in $\T_y(\Gamma)$ parallel to $e$, and thus to $\nu \in \T_x(\Gamma)$, and denote by $\sss_0^\nu(y) < \sss_1^\nu(y) < \dots < \sss_r^\nu(y)$ the corresponding slopes.
		
		\begin{prop}[Non-increasing property of slope vectors] \label{prop:non_increasing_slope_vectors}
			Notation as above, the collection of vectors $(\sss_0^\nu, \dots, \sss_r^\nu)$, as a vector-valued function on the segment corresponding to $e$, forms a coordinate-wise non-increasing collection of vectors.
		\end{prop}
		
		In other words, for every $\varepsilon > 0$ small enough, denoting by $y = x + \varepsilon \nu$ the point at distance $\varepsilon$ from $x$ in the direction of $\nu$, we have
		\[ \sss_j^\nu(y) \leq \sss_j^\nu(x) \qquad \text{for every } j = 0, 1, \dots, r. \]
		
		\begin{proof}
			Since the vector of slopes is piecewise constant, we can suppose that $x$ is in the interior of $e$. Changing the model by adding $x$ as a new vertex of $\Gamma$, if necessary, we can suppose that $x$ is a vertex of $G$. Denoting by $e'$ the edge emanating from $x$ in the direction $\nu$ and by $e''$ the other edge incident to $x$, oriented toward $x$ (that is, with an orientation compatible with that of $e'$). To prove the proposition, it is sufficient to show that for every $j \in \{1, \dots, r\}$, $\sss_j^{e'}(x) \leq \sss_j^{e''}(x)$. We will use the divisorial rank property in the definition of crude linear series. In the following, $\overline{e''}$ denotes the edge $e''$ with the reverse orientation.
			
			Let $\ssp_1, \dots, \ssp_j$ be $j$ distinct points on $e'$ close enough to $x$, in this order away from $x$, such that the slope structure is constant between $x$ and $\ssp_j$. Likewise, let $\ssq_1, \dots, \ssq_{r - j}$ be $r - j$ distinct points on $e''$ close enough to $x$, in this order away from $x$, such that the slope structure is constant between $x$ and $\ssq_{r - j}$. Consider the effective divisor of degree $r$
			\[ E \coloneqq \sum_{i = 1}^j (\ssp_i) + \sum_{i = 1}^{r - j} (\ssq_i). \]
			By property $(*)$ in Definition~\ref{subsec:def_crude_linear_series}, there exists $f \in \Rat(D, \fS)$ such that $D + \div(f) - E$ is an effective divisor. The vector of outgoing slopes of $f$ around $x$, $\ssdelta_x(f)$, corresponds to some jump $\unda \in \ssJ_{\ssrho_x}$. Since by construction $D$ has no support between $\ssq_{r - j}$ and $\ssp_j$, the inequality $D + \div(f) - E \geq 0$ implies that $f$ has a positive order of vanishing on all the points $\ssp_i$ and $\ssq_i$. This in turn implies that $\ssa_{e'} \geq j$ and $\ssa_{\overline{e''}} \geq r - j$. The fact that $\div(f)(x) \geq 0$ implies that $\sss_{a_{e'}}^{e'}(x) + \sss_{a_{\overline{e''}}}^{\overline{e''}}(x) \leq 0$. Finally,
			\[ \sss_j^{e'}(x) \leq \sss_{a_{e'}}^{e'}(x) \leq - \sss_{a_{\overline{e''}}}^{\overline{e''}}(x) \leq - \sss_{r - j}^{\overline{e''}}(x) = \sss_j^{e''}(x). \]
		\end{proof}
		
		\begin{thm}[Theorem~\ref{thm:finiteness_slope_structures_intro} on the finiteness of slopes structures in crude linear series] \label{thm:finiteness_slope_structures}
			Let $\Gamma$ be a metric graph and let $D$ be a divisor on $\Gamma$. Let $G = (V, E)$ be a combinatorial model of $\Gamma$ that supports $D$ on its vertices. For each integer $r$, there exist finitely many subdivisions $\ssH_1, \dots, \ssH_k$ of $G$, and an $r$-slope structure $\fS_j$ defined on $\ssH_j$ for $j = 1, \dots, k$, such that every crude linear series $(D, \fS)$ of rank $r$ has a model $\ssH_j$ of $\Gamma$ among $\ssH_1, \dots, \ssH_k$ on which it is defined, and moreover, the equality $\fS = \fS_j$ holds.
		\end{thm}
		
		\begin{proof}
			Using \cite[Lemma~1.8]{gathmann2008riemann}, we infer that the slopes appearing in the slope structure $\fS$ are all bounded. Applying Proposition~\ref{prop:non_increasing_slope_vectors}, this implies that the number of graph models over which the slope structure is defined is finite, and there are only finitely many possibilities for rank functions on the vertices of each of these graph models. The result follows.
		\end{proof}


\section{Admissible semimodules}
	
	In this section, we introduce admissible semimodules of rational functions. The idea is to replace the full space of rational functions $\Rat(D, \fS)$ in a crude linear series (which is in general not closed, in the topological sense) by a closed semimodule, still enjoying the properties of crude linear series.
	
	\subsection{Semimodule structure on spaces of rational functions} \label{subsec:structure_rat}
		
		Let $\TF = (\R, \oplus, \odot)$ be the semifield of tropical numbers. 
		The space $\R^\Gamma$ of real-valued functions on $\Gamma$ is naturally a $\TF$-semimodule: for $f, g \in \R^\Gamma$ and $c \in \R$, we have operations of tropical addition and tropical multiplication by scalars
		\[ f \oplus g \coloneqq \min(f, g) \qquad \textrm{and} \qquad c \odot f \coloneqq f + c. \]
		
		The relevance of this discussion is in the following basic result, which shows that the space of rational functions associated to a divisor endowed with a slope structure is a $\TF$-semimodule.
		
		\begin{prop} \label{prop:structure_linear_series}
			The space $\Rat(\Gamma)$ is a $\TF$-semimodule. Both subsets $\Rat(D)$ and $\Rat(D, \fS)$, for a divisor $D$ and a divisor endowed with a slope structure, respectively, are $\TF$-semimodules.
		\end{prop}
		
		\begin{proof}
			For the first two statements see~\cite[Lemma~4]{haase2012linear}. In order to prove the last statement, we need to show that $\min(f, g)$ belongs to $\Rat(D, \fS)$ if $f$ and $g$ do. Let $v$ be a vertex of $\Gamma$ and let $\ssd_v$ be the valence of $v$. Denote by $\unda$, $\undb$ and $\undc$ the elements of $\hcube{r}{\ssd_v}$ such that $\sss_{\unda} = \ssdelta_v(f)$, $\sss_{\undb} = \ssdelta_v(g)$ and $\sss_{\undc} = \ssdelta_v(\min(f, g))$. If $f(v) = g(v)$, then we have $\ssdelta_v(\min(f, g)) = \min(\ssdelta_v(f), \ssdelta_v(g))$, which implies $\undc = \unda \wedge \undb$. Otherwise, we have $\undc = \unda$ or $ \undb$ depending on whether $f(v) < g(v)$ or $f(v) > g(v)$, respectively. Using that $\ssJ_{\ssrho_v}$ is stable under $\wedge$ (see Proposition~\ref{prop:jump_stability}), we conclude that, in either case, $\undc$ belongs to $\ssJ_{\ssrho_v}$, and the proposition follows.
		\end{proof}
		
		We introduce some terminology that we use later.
		
		\begin{defi} \label{nota:modification_semi_module}
			Let $S$ be a subset of $\R^\Gamma$ and $v \in \Gamma$. By $\ssS_v$, we mean the space of all functions $f$ of $S$ such that $f(v) = 0$. This will be used, in particular, when $S$ is a sub-semimodule $M$ of some linear series $\Rat(D, \fS)$.
			
			For a subset $S$ of $\R^\Gamma$ and $f \in \R^\Gamma$, we also define $S(-f) \coloneqq S - f = \mleft\{h - f \st h \in S\mright\}$.
		\end{defi}
		
		The latter definition mimics the linear equivalence relation between divisors or slope structures (Section~\ref{subsec:divisors_endowed_slope_structure}). Notice that if $M$ is a sub-semimodule of $\R^\Gamma$, then this is also true of $M(-f)$.
		
		\begin{defi}
			A semimodule $M \subseteq \R^\Gamma$ is called \emph{effective} if it contains the null function.
		\end{defi}
		
		\begin{remark} \label{rem:semi_module_make_effective}
			This definition extends Definition~\ref{def:effective_divisor} where $M = \Rat(D, \fS)$. Note that if $f \in M$, then $M(-f)$ is effective.
		\end{remark}
	
	\subsection{Restriction and extension of scalars}
		
		The above notions generalize easily to the case of $\Lambda$-rational divisors on $\Lambda$-metric graphs.
		
		If $\Gamma$ and $D$ are $\Lambda$-rational and if $M$ is a sub-$\TF$-semimodule of $\Rat(D, \fS)$, we define $\ssM_\Lambda \coloneqq M \cap \Rat!_\Lambda(\Gamma)$, the sub-$\Lambda$-semimodule of $\Rat!_\Lambda(D, \fS)$ made up of elements $f \in M$ which are $\Lambda$-rational. We say that $\ssM_\Lambda$ is obtained from $M$ by \emph{restriction of scalars}.
		
		Assume that $\Gamma$ and $D$ are $\Lambda$-rational. Let $\Lambda'$ be another sub-semi-field of $\R$ such that $\Lambda \subseteq \Lambda'$, and $M$ be a sub-$\Lambda$-semimodule of $\Rat!_\Lambda(D, \fS)$. We define $\ssub{M}!^{\Lambda'}$ to be the sub-$\Lambda'$-semimodule of $\Rat!_{\Lambda'}(D, \fS)$ generated by $M$. We say that $\ssub{M}!^{\Lambda'}$ is obtained from $M$ by \emph{extension of scalars}.
	
	\subsection{Finite generation and closedness} \label{subsec:finite_generation_closedness}
		
		A semimodule $M$ over a commutative semi-ring $R$ is \emph{finitely generated} if there exist $\ssf_1, \dots, \ssf_n \in M$ such that for all $g \in M$, there exist $\ssc_1, \dots, \ssc_n \in R$ such that $g = \bigoplus_{1 \leq i \leq n} \ssc_i \odot \ssf_i$.
		
		Using cut sets and extremal points, Haase, Musiker and Yu showed in \cite{haase2012linear} that for every divisor $D$, $\Rat(D)$ is a finitely generated $\TF$-semimodule. This is not necessarily the case for spaces of the form $\Rat(D, \fS)$ for crude linear series $(D, \fS)$ (see Example~\ref{ex:finiteness_condition}).
		
		A semimodule $M \subseteq \Rat(\Gamma)$ is called \emph{closed} if it is closed for the norm $\|\cdot\|_\infty$.
		
		We have the following basic result on the connection between finite generation and closedness properties.
		
		\begin{prop} \label{prop:finiteness_closedness}
			Let $M \subseteq \Rat(\Gamma)$ be a finitely generated semimodule. Then, $M$ is closed.
		\end{prop}
		
		\begin{proof}
			Let $(\ssf_n)_n$ be a sequence of functions of $M$ converging to a function $f$ in $\Rat(\Gamma)$ for the $\| \cdot \|_\infty$ topology. Assume that $M$ is generated by elements $\ssh_1, \dots, \ssh_r$, and write, for all $n$, $\ssf_n = \min_{1 \leq i \leq r} \mleft(\ssh_i + \ssc^n_i\mright)$ with $\ssc^n_i \in \R$. We can suppose that all $\ssh_i$ are zero at some point $v$. It follows that the sequences $\mleft(\ssc^n_i\mright)_n$ are bounded. By extraction, we can assume $\mleft(\ssc^n_i\mright)_n$ converges to some $\ssc_i \in \R$, implying that that for every $i$, $\ssf_n \xrightarrow{\| \cdot \|_\infty} \min_{1 \leq i \leq r} \mleft(\ssh_i + \ssc_i\mright)$, and thus $f \in M$.
		\end{proof}
		
		We next give an alternative characterization of closedness for semimodules $M \subseteq \Rat(D)$. First, notice that each effective divisor $D$ of degree $d$ is written in the form $D = (\ssx_1) + \cdots + (\ssx_d)$ for points $\ssx_j \in \Gamma$, and can be viewed as a point $[\ssx_1, \dots, \ssx_d]$ in the $d$-th symmetric product $\mathrm{Sym}^d(\Gamma)$ of the metric graph, which is a compact metric space. Furthermore, we can view $\mathrm{Sym}^d(\Gamma)$ as the subspace of $\Div^d(\Gamma)$ consisting of effective divisors of degree $d$. This defines a map
		\[ \varphi \colon \Rat(D) \longrightarrow \mathrm{Sym}^d(\Gamma) \hookrightarrow \Div^d(\Gamma) \]
		by $\varphi(f) \coloneqq D + \div(f)$. We set $|M| = \varphi(M)$, and view it in $\mathrm{Sym}^d(\Gamma)$. We endow $\Rat(D)$ with the $\|\cdot\|_\infty$ topology.
		
		\begin{prop} \label{prop:closedness_equivalence}
			The map $\varphi$ is continuous. Furthermore, a semimodule $M \subseteq \Rat(D)$ is closed if, and only if, $|M| \subseteq \mathrm{Sym}^d(\Gamma)$ is closed.
		\end{prop}
		
		\begin{proof}
			We omit the proof of the first assertion, and prove only the second.
			
			$\Longrightarrow$ If $M$ is closed then so is $M_v$, and $|M| = \varphi(M_v)$. Let $k$ be a universal bound for the slopes in $\fS$ (see~\cite[Lemma 1.8]{gathmann2008riemann}), and let $B$ be the space of continuous functions $\Gamma \to \R$ whose slopes are bounded by $k$. By the Arzel\'a--Ascoli theorem, $B$ is compact, and since $M_v \subseteq B$ is closed, we infer that $M_v$ is compact. Since $\varphi$ is continuous, $|M|$ is compact and thus closed in $\mathrm{Sym}^d(\Gamma)$.
			
			$\Longleftarrow$ Let $(\ssf_n)_n$ be a sequence of functions of $M$ converging to some $f \in \Rat(D)$. Then by continuity $\varphi(\ssf_n) \longrightarrow \varphi(f)$. By closedness, $\varphi(f) \in |M|$, so there exists some $g \in M$ such that $\varphi(f) = \varphi(g)$. Since $\div(f) = \div(g)$, $f$ and $g$ differ by some constant, and thus $f \in M$.
		\end{proof}
	
	\subsection{Admissible semimodules} \label{subsec:admissible_semimodules}
		
		Let $d$ be an integer and $D$ a divisor of degree $d$ on $\Gamma$.
		
 		A semimodule $M \subseteq \Rat(D)$ is called \emph{admissible of rank $r$} if it is closed for the topology induced by $\| \cdot \|_\infty$, and if there exists an $r$-slope structure $\fS$ such that $M \subseteq \Rat(D, \fS)$ and such that the following holds:
 		
		\smallskip
		
		$(**)$ For every effective divisor $E$ on $\Gamma$ of degree $r$, there exists $f \in M$ such that
		
		\begin{enumerate}[label=(AS\arabic*)]
			\item \label{def:admissible_semimodule_1} For every point $x \in \Gamma$, $\ssrho_x(\ssdelta_x(f)) \geq E(x)$; and in addition,
			
			\item \label{def:admissible_semimodule_2} $D + \div(f) - E \geq 0$.
		\end{enumerate}
		
		\smallskip
		
		When $\Gamma$ and $D$ are $\Lambda$-rational, we say a semimodule $M \subseteq \ssub{\Rat}!_\Lambda(D, \fS)$ is admissible if the extension $\ssM^\R \subseteq \Rat(D, \fS)$ is admissible. \hfill $\diamond$
		
		Note that if $M$ is effective as a semimodule, then so is $\Rat(D, \fS)$ and therefore $(D, \fS)$ is by definition effective as a crude linear series.
		
		Although it is not clear from the definition, we will show in the next section that admissibility is a property of the pair $(D, M)$, that is, both the rank $r$ and the slope structure $\fS$ (and thus the crude linear series $(D, \fS)$) can be extracted from $M$.
		
		Example~\ref{ex:finiteness_condition} shows that, in general, the semimodule $\Rat(D, \fS) \subseteq \Rat(D)$ for a crude linear series $(D, \fS)$ may not be closed in $\Rat(D)$, and therefore may not be admissible. Furthermore, a crude linear series $(D, \fS)$ of rank $r$ might not necessarily admit an admissible semimodule $M \subseteq \Rat(D, \fS)$ of the same rank, although we do not have an example thereof.
		
		\begin{remark} \label{rem:vector_slopes_increasing_limit}
			For further use, we note that if $(\ssf_n)_n$ is a sequence of functions in $\Rat(D)$ converging uniformly to a function $f$, then $f \in \Rat(D)$ and, for every $x \in \Gamma$, the inequality $\ssdelta_x(f) \leq \ssdelta_x(\ssf_n)$ holds coordinate-wise for all large $n$. This inequality can be strict in general, and the corresponding point of the hypercube may not be a jump of the rank function at $x$. This explains why a semimodule of the form $\Rat(D, \fS)$ is not necessarily closed.
		\end{remark}
	
	\subsection{Realization of jumps in admissible semimodules} \label{subsec:realization_jumps}
		
		Let $D$ be a divisor of degree $d$. In this section, we prove the following theorem.
		
		\begin{thm}[Realization of jumps in admissible semimodules] \label{thm:jumps_realized}
			Let $(D, \fS)$ be a crude linear series of degree $d$ and rank $r$, and let $M \subseteq \Rat(D, \fS)$ be an admissible semimodule of rank $r$. Let $v$ be a point of $\Gamma$ and let $\unda \in \ssJ_{\ssrho_v}$ be a jump. Then, there exists $f \in M$ such that $\ssdelta_v(f) = \sss_\unda$.
		\end{thm}
		
		For example, taking $E = r(v)$ for $v \in \Gamma$, we deduce from property $(**)$ the existence of a function $f \in M$ such that $\ssrho_v(\ssdelta_v(f)) \geq r$. This implies that $f$ takes all minimum slopes around $v$, that is, $\ssdelta_v(f) = \sss_{\underline 0}$. The proof for other values of $\unda$ is more involved.
		
		\smallskip
		
		This theorem immediately implies the following important results.
		
		\begin{cor}[Corollary~\ref{cor:slope_structure_retrieved_linear_series_intro}] \label{cor:slope_structure_retrieved_linear_series}
			Keeping the notation of Theorem~\ref{thm:jumps_realized}, the slope structure $\fS$ can be entirely retrieved from $M$.
		\end{cor}
		
		In light of this corollary, we will often drop $\fS$ and only say that $M$ is an admissible semimodule of $\mathrm{Rat}(D)$.
		
		\begin{cor} \label{cor:equivalence_slope_substructure}
			Let $(D, \fS)$ and $(D', \fS')$ be two crude linear series on $\Gamma$ of rank $r$ and $r'$, respectively, and $M'\subseteq \Rat(D', \fS')$ and $M \subseteq \Rat(D, \fS)$ be admissible. Then $\fS'$ is a slope substructure of $\fS$ provided that $M' \subseteq M$ holds.
		\end{cor}
		
		\begin{remark} \label{rem:minimal_divisor_admissible_semimodule}
			Note that the divisor $D$ cannot be uniquely determined by $M$, because if $M \subseteq \Rat(D, \fS)$, then $M \subseteq \Rat(D', \fS)$ for every divisor $D' \geq D$.
			However, there is a unique base-point free choice for $D$, that is, a unique minimal choice $\ssD_{\min}$ for $D$. Simply define the coefficient $\ssD_{\min}(x)$ of $\ssD_{\min}$ as the minimum possible integer such that all the functions $f \in M$ verify $\ord_x(f) + \ssD_{\min}(x) \geq 0$. This is well-defined, and every other choice $D$ with $M \subseteq \Rat(D, \fS)$ verifies $D \geq \ssD_{\min}$.
		\end{remark}
		
		The rest of the section is devoted to the proof of the realization theorem.
	
	\subsection{Proof of Theorem~\ref{thm:jumps_realized}}
		
		Consider a model $G = (V, E)$ such that $\fS$ is defined on $G$ and let $v$ be a vertex of $V$. Let $\grdeg = \ssd_v$. We show that for each jump $\unda \in \ssJ_{\ssrho_v}$, there exists $f \in M$ such that $\ssdelta_v(f) = \sss_\unda$. The idea is to use the divisorial rank condition $(**)$ satisfied by the admissible semimodule $M$ in order to find a sequence of functions with the desired slope at a given tangent direction at $v$. We then find $f$ as the limit of a subsequence, using the closedness of $M$.
		
		We start with a lemma.
		
		\begin{lemma} \label{lem:unique_jump_dominating_point}
			Let $\undq \in \hcube{r}{\grdeg}$ be any point in the hypercube at $v$ such that $\ssrho_v \mleft(\undq\mright) \geq 0$. Then, there exists a unique jump $\unda \in \ssJ_{\ssrho_v}$ of rank $\ssrho_v \mleft(\undq\mright)$ with $\unda \supfaceeq \undq$.
		\end{lemma}
		
		\begin{proof}
			We reason by contradiction to show existence. If no such jump exists, we can construct an increasing path (i.e., a finite sequence of points $\undx^i$ of $\hcube{r}{\grdeg}$ such that $\undx^i \subface \undx^{i + 1}$ for every $i$), starting at $\undq$ and made up of points of constant rank. This process necessarily ends up at the point $\undr = (r, \dots, r)$. We know that $\ssrho_v(\undr) = \ssrho_v \mleft(\undq\mright) = 0$. Therefore, $\undr$ is a jump, $\undr \supfaceeq \undq$ and they have the same rank, a contradiction.
			
			We now show uniqueness. Suppose there are two different jumps $\unda$ and $\undb$ such that $\unda \geq \undq$, $\undb \geq \undq$, $\rho_v(\unda) = \ssrho_v(\undb) = \ssrho_v \mleft(\undq\mright)$. Since $\unda \neq \undb$, then $\unda \wedge \undb$ is different from at least one among $\unda$ and $\undb$, say $\unda$. Proposition~\ref{prop:jumps_graded_poset} yields that $\ssrho_v(\unda \wedge \undb) > \ssrho_v(\unda)$, which is impossible because $\undq \leq \unda \wedge \undb \leq \unda$ implies $\ssrho_v(\unda) = \ssrho_v \mleft(\undq\mright) \geq \ssrho_v(\unda \wedge \undb) \geq \ssrho_v(\unda)$.
		\end{proof}
		
		We now come to the proof of the theorem. We start by defining an increasing path in $\hcube{r}{\grdeg}$ that starts at $\underline 0$, stays below $\unda$ and moves only along the first direction in the hypercube $\hcube{r}{\grdeg}$ at the beginning, then only along the second direction, and so on, until direction $\grdeg$.
		
		Let, for all $i \in \{1, \dots, \grdeg\}$, $\sse_i$ be the edge incident to $v$ corresponding to the direction $i$. For convenience, we also define, for all $i \in \{1, \dots, \grdeg\}$, $\sse_i' \coloneqq \sse_{\grdeg - i + 1}$, just reversing the order of the edges around $v$.
		
		Given a point $\undy \in \hcube{\grdeg}{r}$ with $\undy \subfaceeq \unda$ and a direction $\unde_i$ in the hypercube, we say that $\undy$ is \emph{a fall in the direction $\unde_i$ with respect to $\unda$} if $\ssy_i = 0$ or $\ssrho_v\mleft(\undy\mright) < \ssrho_v\mleft(\undy - \unde_i\mright)$. We say that $\undy$ is \emph{a largest fall in the direction $\unde_i$ with respect to $\unda$} if $\undy$ is a fall in the direction $\unde_i$ with respect to $\unda$ and if for all non-negative integers $n$ such that $\undy + n \, \unde_i \subfaceeq \unda$, we have $\ssrho_v\mleft(\undy + n \, \unde_i\mright) = \ssrho_v\mleft(\undy\mright)$. Saying in words, moving in the direction of $\unde_i$ from a largest fall, remaining bounded by $\unda$, the rank does not change.
		
		We construct our increasing path starting from $\underline 0$, going each time to a largest fall in the given direction relative to $\unda$ and remembering only the falls in that direction. More precisely, suppose that we have already built the path along directions $1, \dots, i - 1$ with $1 \leq i \leq \grdeg$, consisting of all the falls in the direction of $\unde_1$, then the falls in the direction of $\unde_2$, $\dots$, and the falls in the direction of $\unde_{i - 1}$. Therefore, the path currently ends at the point $\sum_{k = 1}^{i - 1} \sst_k \, \unde_k$, $\sst_k \in [r]$. We will now let the path continue only in the direction $i$, by adding multiples of $\unde_i$, adding the falls in the direction of $\unde_i$ to the path, until we reach the point $\undy = \sum_{k = 1}^{i - 1} \sst_k \, \unde_k + \sst_i \, \unde_i$ with $\undy \subfaceeq \unda$, which is a largest fall in the direction $\unde_i$ with respect to $\unda$. This way, we have built an increasing path starting at $\underline 0$, composed only of falls in some direction, staying below $\unda$ and moving successively in directions $\unde_1, \dots, \unde_{\grdeg}$. The turning point, where we move to a different direction, is a largest fall with respect to $\unda$ in that direction. The ending point of the path
		\[ \undz = \sum_{k = 1}^\grdeg \sst_k \, \unde_k \]
		is a largest fall in the direction $\unde_{\,\grdeg}$ with respect to $\unda$. In particular, note that $\undz \subfaceeq \unda$.
		
		For every $i \in \{1, \dots, \grdeg\}$, let $\ssell_i$ be the rank drop of the path in the direction $i$:
		\[ \ssell_i \coloneqq \ssrho_v\mleft(\sum_{k = 1}^{i - 1} \sst_k \, \unde_k\mright) - \ssrho_v\mleft(\sum_{k = 1}^i \sst_k \, \unde_k \mright). \]
		Then, $\sum_{i = 1}^\grdeg \ssell_i = r - \ssrho_v(\undz)$.
		
		\begin{lemma}
			Let $\undz$ be the endpoint of the increasing path constructed above. We have $\ssrho_v(\undz) = \ssrho_v(\unda)$. Therefore, $\ssrho_v(\unda) + \sum_{i = 1}^\grdeg \ssell_i = r$.
		\end{lemma}
		
		\begin{proof}
			Let $\undw \in \hcube{r}{\grdeg}$ be such that $\undz \subfaceeq \undw \subfaceeq \unda$ and let $i \in \{1, \dots, \grdeg\}$ be such that $\undw + \unde_i \subfaceeq \unda$. By construction, we know that
			\[ \ssrho_v\mleft(\sum_{k = 1}^i \sst_k \, \unde_k + \unde_i\mright) = \ssrho_v\mleft(\sum_{k = 1}^i \sst_k \, \unde_k\mright). \]
			This equality, together with supermodularity, implies that $\ssrho_v(\undw + \unde_i) = \ssrho_v(\undw)$. Applying this fact recursively yields $\ssrho_v(\undz) = \ssrho_v(\unda)$. The last equality follows from $\ssrho_v(\undz) + \sum_{i = 1}^\grdeg \ssell_i = r$.
		\end{proof}
		
		For every $i \in \{1, \dots, \grdeg\}$, let $\ssp_1^i, \dots, \ssp_{\ssell_i}^i$ be distinct points on the edge $\sse_i' = \sse_{n - i + 1}$, ordered increasingly with respect to their distance from $v$. Let
		\[ 
			E \coloneqq \ssrho_v(\unda) \, (v) + \sum_{i = 1}^\grdeg \sum_{j = 1}^{\ssell_i} \mleft(\ssp_j^i\mright).
		\]
		By the preceding lemma, this is an effective divisor of degree $r$. Since $M$ is admissible, there exists $f \in M$ such that $f(v) = 0$, and
		\begin{equation} \label{eq:rank_slope}
			\ssrho_v(\ssdelta_v(f)) \geq \ssrho_v(\unda) \quad \textrm{and} \quad \div(f) + D - E \geq 0.
		\end{equation}
		Let $\undb \coloneqq \ssdelta_v(f) \in \hcube{r}{\grdeg}$. Note that $\undb$ is a jump of $\ssrho_v$. The following lemma will imply that $\undb = \unda$.
		
		\begin{lemma}
			We have $\undb \supfaceeq \undz$.
		\end{lemma}
		
		\begin{proof}			 
			The first of the two properties in~\eqref{eq:rank_slope} tells us that $\ssrho_v(\undb) \geq \ssrho_v(\unda) = \ssrho_v(\undz)$.
			
			We will now show that the second property in~\eqref{eq:rank_slope}, when applied to a sequence of effective divisors of degree $r$ that starts at $E$, implies that $\undb \supfaceeq \undz$. The idea is to make all the points in the support of $E$ converge to $v$ one after the other, apply the admissibility of $M$, and define a decreasing sequence of jumps.
			
			Let $\mleft(\ssp_{1,n}^1\mright)_{n \geq 0}$ be a sequence of points of $\sse_1'$ starting at $\ssp_1^1$ and which converges to $v$. We replace $\ssp_1^1$ in $\ssE$ by $\ssp_{1,n}^1$ to obtain an effective divisor $\ssE_{1,n}^1$ of degree $r$. Formally, 
			\[ \ssE_{1,n}^1 \coloneqq \ssrho_v(\unda) \, (v) + \mleft(\ssp_{1,n}^1\mright) + \sum_{j = 2}^{\ssell_1} \mleft(\ssp_j^1\mright) + \sum_{i = 2}^{\grdeg} \sum_{j = 1}^{\ssell_i} \mleft(\ssp_j^i\mright). \]
			
			Applying the admissibility of $M$, we find a function $\ssf_{1,n}^1 \in M$ with $\ssf_{1,n}^1(v) = 0$ such that $\div(\ssf_{1,n}^1) + D - \ssE_{1,n}^1 \geq 0$. Replacing $\ssf_{1,n}^1$ by $\min\mleft(f, \ssf_{1,n}^1\mright) \in M$ ensures that for all $n$, $\ssdelta_v\mleft(\ssf_{1,n}^1\mright) \subfaceeq \ssdelta_v(f)$. Moreover, since $M$ is closed, we can extract a subsequence so that the sequence $\mleft(\ssf_{1,n}^1\mright)_n$ converges to an element $\ssf_1^1 \in M$. Let $\ssE_{1}^1 = \ssE - (\ssp^1_1) + (v)$.
			
			By Remark~\ref{rem:vector_slopes_increasing_limit}, for all $n$,
			\[ \ssdelta_v\mleft(\ssf_1^1\mright) \subfaceeq \ssdelta_v\mleft(\ssf_{1,n}^1\mright) \subfaceeq \ssdelta_v(f). \]
			Besides, the fact that $\ssub{\mleft(\ssdelta_v\mleft(\ssf_1^1\mright)\mright)}!_1 < \ssub{\mleft(\ssdelta_v\mleft(\ssf_{1,n}\mright)\mright)}!_1$ implies that $\ssdelta_v\mleft(\ssf_1^1\mright) \neq \ssdelta_v(f)$.
			
			What precedes yields a jump $\undb^1_1 \coloneqq \ssdelta_v\mleft(\ssf_1^1\mright) \subfaceeq \ssdelta_v(f) = \undb$, different from $\undb$ because the first coordinates verify $\ssub{\mleft(\undb^1_1\mright)}!_1<\ssb_1$. Therefore, $\undb^1_1\subface \undb$.
			
			We repeat the same process as above, starting from $\ssE^1_1$. We take a sequence $\mleft(\ssp_{2,n}^1\mright)_n$ of points of $\sse_1'$ starting at $\ssp_2^1$, and define
			\[ \ssE_{2,n}^1 \coloneqq (\ssrho_v(\unda) + 1) \, (v) + \mleft(\ssp_{2,n}^1\mright) + \sum_{j = 3}^{\ssell_1} \mleft(\ssp_j^1\mright) + \sum_{i = 2}^\grdeg \sum_{j = 1}^{\ssell_i} \mleft(\ssp_j^i\mright), \]
			and find an element $\ssf_{2,n}^1 \in M$ with $\ssf_{2,n}^1(v) = 0$ such that $\div(\ssf_{1,n}^1) + D - \ssE_{2,n}^1 \geq 0$.
			After taking the minimum $\min\mleft(\ssf^1_1, \ssf_{2,n}^1\mright) \in M$, and passing to a subsequence, we obtain a limit $\ssf^1_2\in M$. This yields a jump $\undb^1_2 \subfaceeq \undb^1_1$, different from $\undb^1_1$ because $\ssub{\mleft(\undb^1_2\mright)}!_1 < \ssub{\mleft(\undb^1_1\mright)}!_1$.
			
			We repeat the same process, exhausting first all the points $\ssp_j^1$ on $\sse_1'$ in the support of $E$, then all the points $\ssp_j^2$ on $\sse_2'$, and so on, until finally all the points $\ssp_j^\grdeg$ on the last edge. The above reasoning yields a decreasing path of jumps
			\[ \undb \supface \undb^1_1 \supface \undb^1_2 \supface \cdots \supface \undb^1_{\ssell_1} \supface \cdots \supface \undb^\grdeg_1 \supface \cdots \supface \undb^\grdeg_{\ssell_\grdeg}. \]
			Using Proposition~\ref{prop:jumps_graded_poset}, the sequence of ranks $\mleft(\ssrho_v(\undb^i_j)\mright)$ is increasing. As a consequence,
			\[ \ssrho_v(\undb^\grdeg_{\ssell_\grdeg}) \geq \ssrho_v\mleft(\undb\mright) + \sum_{i = 1}^\grdeg \ssell_i \geq r. \]
			We thus infer that $\undb^\grdeg_{\ssell_\grdeg} = \underline 0$ and all the rank differences between consecutive jumps in the sequence are exactly one.
			
			Reversing the order in the sequence of jumps constructed above, we get an increasing path of jumps starting at $\underline 0$ and ending at $\undb$, whereas we defined beforehand an increasing path starting at $\underline 0$ and ending at $\undz$. The two increasing paths have the same length, equal to $r - \ssrho_v(\unda)$. To show that $\undb \supfaceeq \undz$, we will prove that the path leading to $\undb$ remains greater than or equal to the path leading to $\undz$ at each step, that is, the $k$-th element of the former dominates the $k$-th element of the latter.
			
			We proceed by induction. The claim is true at the beginning because the starting point of both paths is $\underline 0$. Suppose that the inequality is true at some step $\ssj_0$ with $0 \leq \ssj_0 < r - \ssrho_v(\unda)$. We denote by $\undz_{\ssj_0}$ the current fall and by $\undb_{\ssj_0}$ the current jump. The inequality reads $\undb_{\ssj_0} \supfaceeq \undz_{\ssj_0}$. We suppose that the next fall $\undz_{\ssj_0 + 1}$ differs from $\undz_{\ssj_0}$ (only) in the direction $\ssi_0$. Let
			\[ \undc_{\ssj_0} \coloneqq \undz_{\ssj_0} + \mleft(\mleft(\undb_{\ssj_0}\mright)_{\ssi_0} - \mleft(\undz_{\ssj_0}\mright)_{\ssi_0}\mright) \, \unde_{\ssi_0} \]
			be the point obtained by starting at $\undz_{\ssj_0}$ and moving in the direction $\ssi_0$ as much as possible without overtaking $\undb_{\ssj_0}$ along this axis.
			
			Since both paths are parametrized by the same integers $\ssell_i$, we know that the next jump $\undb_{\ssj_0 + 1}$ will also differ from $\undb_{\ssj_0}$ (at least) in the direction $\ssi_0$. The latter statement implies that $\undb_{\ssj_0 + 1} \supfaceeq \undb_{\ssj_0} + \unde_{\ssi_0}$. Since $\undb_{\ssj_0}$ is a jump, by supermodularity and using $\undb_{\ssj_0} \supfaceeq \undc_{\ssj_0}$, we get that
			\[ \ssrho_v\mleft(\undc_{\ssj_0} + \unde_{\ssi_0}\mright) = \ssrho_v\mleft(\undc_{\ssj_0}\mright) - 1 \]
			so $\undc_{\ssj_0} + \unde_{\ssi_0} \supfaceeq \undz_{\ssj_0}$ is a fall in the direction $\ssi_0$ which coincides with $\undz_{\ssj_0}$ in all directions but $\ssi_0$. Therefore,
			\[ \undz_{\ssj_0 + 1} \subfaceeq \undc_{\ssj_0} + \unde_{\ssi_0} \subfaceeq \undb_{\ssj_0} + \unde_{\ssi_0} \subfaceeq \undb_{\ssj_0 + 1}. \]
			We have proved the claim for $\ssj_0 + 1$. We infer that $\undb \supfaceeq \undz$, as desired.
		\end{proof}	
		
		\begin{proof}[Proof of Theorem~\ref{thm:jumps_realized}]
			The preceding lemma shows that $\undb \supfaceeq \undz$, and thus $\ssrho_v(\undb) \leq \ssrho_v(\undz) = \ssrho_v(\unda)$. Combined with the inequality $\ssrho_v(\undb) \geq \ssrho_v(\unda)$, we deduce that $\ssrho_v(\undb) = \ssrho_v(\unda)$. On the other hand, $\undb \supfaceeq \undz$ is a jump of $\ssrho_v$. Consequently, $\unda$ and $\undb$ are two jumps of $\ssrho_v$ which dominate $\undz$ and have the same rank as $\undz$. The uniqueness in Lemma~\ref{lem:unique_jump_dominating_point} implies that $\undb = \unda$, which finishes the proof.
		\end{proof}


\section{Combinatorial limit linear series} \label{sec:linear_series_definition}
	
	In this section, we define combinatorial limit linear series on metric graphs.
	
	\subsection{Tropical rank} \label{subsec:tropical_rank}
		
		In linear algebra, the dimension $d$ of a vector space can be characterized as either the size of a minimal generating set or as the least integer such that every collection of $d + 1$ elements is linearly dependent. For spaces of tropical functions, however, there is \emph{a priori} no direct link between the two above notions. The former one corresponds to the finite generation property in semimodules, discussed in Section~\ref{subsec:admissible_semimodules}. We will need a second notion of finiteness based on tropical independence from~\cite{JP14, JP16}, see also~\cite{AGG09} where a finite version of this was developed.
		
		\begin{defi}[Tropical rank] \label{def:tropical_rank}
			Let $M \subseteq \Rat(\Gamma)$ be a semimodule. We call \emph{tropical rank} of $M$ the least integer $r \in \Z_{\geq 0}$ such that for all elements $\ssf_0, \dots, \ssf_{r + 1} \in M$, there exist $\ssc_0, \dots, \ssc_{r + 1} \in \R$ such that for each $x \in \Gamma$, the minimum in
			\[ \min_{0 \leq i \leq r + 1} (\ssf_i(x) + \ssc_i) \]
			is achieved at least twice, that is, for at least two indices $i \in [r + 1]$. If such $r$ does not exist, we say that $M$ is of infinite tropical rank.
			
			For a $\Lambda$-semimodule, the tropical rank is defined the same way by imposing $\ssc_i \in \Lambda$.
		\end{defi}
	
	\subsection{Limit linear series: Definition}	
		
		\begin{defi} \label{def:grd}
			A \emph{combinatorial limit linear series} of rank $r$ and degree $d$ on a metric graph $\Gamma$, more simply called a \emph{linear series} or a $\g^r_d$, is a pair $(D, M)$ consisting of a divisor $D$ of degree $d$ and an admissible semimodule $M \subseteq \Rat(D)$ of rank $r$ which is moreover finitely generated and has tropical rank $r$. The linear series is called \emph{effective} if $M$ is so.
			
			\smallskip
			
			When $\Gamma$ and $D$ are $\Lambda$-rational, we define a $\Lambda$-linear series, or simply $\Lambda$-$\g^r_d$, as a finitely generated semimodule $M \subseteq \ssub{\Rat}!_\Lambda(D)$ of tropical rank $r$ such that the semimodule $\ssM^\R \subseteq \Rat(D)$, obtained by extension of scalars, is admissible of rank $r$.
		\end{defi}
		
		Note that $\ssM^\R$ is finitely generated, and by Proposition~\ref{prop:finiteness_closedness}, it is automatically closed.
		
		\begin{remark}
			For a linear series $(D, M)$, by Corollary~\ref{cor:slope_structure_retrieved_linear_series}, there is a unique slope structure $\fS$ with $M \subseteq \Rat(D,\fS)$, that is, $\fS$ is entirely determined by $M$. This explains why we do not include the data of $\fS$ in the linear series.
		\end{remark}
		
		The relevance of this will be explained by Theorem~\ref{thm:specialization_linear_series_grd}, which states that a semimodule $M \subseteq \Rat(D)$ that comes from the tropicalization of a linear series is finitely generated, is of tropical rank $r$, and also verifies property $(**)$ in Section~\ref{subsec:admissible_semimodules}. In fact, such an $M$ verifies the stronger property $\threestars$ below. For the notion of slope substructure, we refer to Section~\ref{subsec:slope_substructures}.
		
		\begin{defi}[Refined linear series] \label{def:refined_linear_series}
			A \emph{refined linear series}, or \emph{refined $\g^r_d$}, is a pair $(D, M)$ which is a $\g^r_d$ and which in addition verifies the following stronger version of $(**)$:
			
			\smallskip
			
			$\threestars$ For every effective divisor $E$ of degree $s \leq r$ on $\Gamma$, there exists a semimodule $\ssM_E \subseteq M$ such that $(D, \ssM_E)$ is a $\g^{r - s}_d$ and for every element $f \in \ssM_E$, we have
			
			\begin{itemize}
				\item[(1)] $\ssrho_x(\ssdelta_x(f)) \geq E(x)$ for each point $x \in \Gamma$, and in addition,
				
				\item[(2)] $D + \div(f) - E \geq 0$. \qedhere
			\end{itemize}
		\end{defi}
		
		Each linear series $(D, M)$ with underlying slope structure $\fS$ gives rise to a linear system $|(D, M)| \subseteq |(D, \fS)|$, defined as follows.
		
		\begin{defi}[Linear system associated to a $\g^r_d$] \label{def:space_rational_functions_grd}
			For a linear series $(D, M)$, we define the linear system $|(D, M)|$ as the set of all effective divisors $E$ on $\Gamma$ of the form $D + \div(f)$ for $f \in M$.
		\end{defi}
		
		\begin{remark} \label{rem:comparison_tropical_rank}
			For a finitely generated admissible semimodule $M \subseteq \Rat(D)$ of rank $r$, being a linear series, that is, having tropical rank $r$, is equivalent to having tropical rank \emph{at most} $r$. In fact, the tropical rank of an admissible semimodule of given rank $r$ is always at least $r$. To see this, we observe using Theorem~\ref{thm:jumps_realized} that at every point of $\Gamma$, in every direction, each of the $r + 1$ slopes in the slope structure underlying $M$ is taken by an element of $M$. Since the $r + 1$ elements of $M$ that take different slopes along that direction are tropically independent, this implies that the tropical rank of $M$ is at least $r$.
		\end{remark}
	
	\subsection{Examples}
		
		We give two simple examples of linear series, and refer to Sections~\ref{subsubsec:elementary_g_2_4_dipole_graph} and~\ref{subsubsec:involved_g_2_4_dipole_graph} for more examples.
		
		\begin{example}[A $\g^1_2$ on the barbell graph]	\label{ex:linear_series_1}		
			Consider the barbell graph $\Gamma$ with edges of arbitrary length, see Figure~\ref{fig:barbell_clls}. This metric graph has genus two and the canonical divisor $K$ has rank one. We define a linear series $(K, M)$, $M \subseteq \Rat(K)$, of rank one on $\Gamma$. Allow slopes $-1 < 1$ on the middle edge and, for $i = 1, 2$, allow slopes $0 < 1$ on both oriented edges $\ssu_i \ssv_i$, in the direction of the arrows.
			
			\begin{figure}[h!]
				\centering
				\begin{tikzpicture}[scale=0.7]
					\coordinate (A) at (0,0);
					\coordinate (B) at (.8,0);
					\coordinate (C) at (2.8,0);
					\coordinate (D) at (4,0);
					\coordinate (E) at (-0.8,0);
					\coordinate (F) at (5.2,0);
					
					\begin{scope}[decoration={markings,mark=at position 0.5 with {\arrow[red]{>}}},
					] 
						\draw[postaction={decorate}] (B) -- (C);
						\draw[postaction={decorate}] (B) arc (0:180:0.8);
						\draw[postaction={decorate}] (B) arc (0:-180:0.8);
						\draw[postaction={decorate}] (C) arc (180:0:1.2);
						\draw[postaction={decorate}] (C) arc (-180:0:1.2);
					\end{scope}
					
					\foreach \c in {B,C,E,F} {
						\fill (\c) circle (2.5pt);
					}
					
					\node[right] at (E) {$\ssv_1$};
					\node[left] at (B) {$\ssu_1$};
					\node[right] at (C) {$\ssu_2$};
					\node[left] at (F) {$\ssv_2$};
					
					\node[above right] at (B) {$1$};
					\node[above left] at (C) {$1$};
					
					\node[below] at (1.8,0) {$(-1, 1)$};
					\node[above] at (0,0.8) {$(0, 1)$};
					\node[below] at (0,-0.8) {$(0, 1)$};
					\node[above] at (4,1.2) {$(0, 1)$};
					\node[below] at (4,-1.2) {$(0, 1)$};
				\end{tikzpicture}
				\caption{The barbell graph, the canonical divisor and the slope structure $\fS$.}
				\label{fig:barbell_clls}
			\end{figure}
			
			We define suitable rank functions on vertices, as follows. For $i = 1, 2$, we endow $\ssv_i$ with the rank function on $\hcube{1}{2} = [1]^2$ defined by the array
			$
			\begin{pmatrix}
				0 & \color{blue} 0 \\
				\color{blue} 1 & 0
			\end{pmatrix}
			$, and endow $\ssu_i$ with the rank function on $\hcube{1}{3} = [1]^3$ whose restrictions to $[1]^2 \times \{0\}$ and $[1]^2 \times \{1\}$ are defined by the following two matrices, respectively:
			$
			\begin{pmatrix}
				0 & \color{blue} 0 \\
				\color{blue} 1 & 0
			\end{pmatrix}, \;
			\begin{pmatrix}
				-1 & -1 \\
				\color{blue} 0 & -1
			\end{pmatrix}
			$.
			Here, jumps are depicted in blue. For the two vertices $\ssu_1$ and $\ssu_2$, the third coordinate in each of the two rank functions corresponds to the middle edge of $\Gamma$. We further endow each other point of $\Gamma$ with the standard rank function. Set $M \coloneqq \Rat(K, \fS)$. The pair $(K, M)$ is a linear series of degree two and rank one on $\Gamma$.
			
			Note that $M$ is not effective, and thus the constant functions do not belong to $M$. Moreover, $(K, M)$ is the unique $\g^1_2$ on $\Gamma$ with the underlying divisor $K$. This shows that the canonical divisor $K$ on $\Gamma$ is not realizable, that is, $K$ is not the tropicalization of a divisor $\mathcal K \in |\omega!_C|$ in the canonical linear system $|\omega!_C|$ on a smooth proper curve $C$ over an algebraically closed field with a non-trivial non-Archimedean valuation (otherwise, the constant functions would belong to the semimodule $M$).
		\end{example}
		
		\begin{example} \label{ex:linear_series_2}
			Consider the following metric graph $\Gamma$ with two edges of equal length.
			\begin{figure}[h!]
				\centering
				\begin{tikzpicture}
					\coordinate (A) at (0,0);
					\coordinate (B) at (1.5,0);
					\coordinate (C) at (3,0);
					
					\draw[postaction=decorate,decoration={markings,mark=at position 1/2 with {\arrow[red]{>}; }}]
					(A) -- (B);
					
					\draw[postaction=decorate,decoration={markings,mark=at position 1/2 with {\arrow[red]{>}; }}]
					(B) -- (C);
					
					\foreach \c in {A, B, C} {
						\fill (\c) circle (2.5pt);
						}
					
					\node[left] at (A) {$u$};
					\node[below] at (B) {$v$};
					\node[right] at (C) {$w$};
				\end{tikzpicture}
			\end{figure}
			
			Let $D = (u) + (w)$ and define a slope structure $\fS$ by allowing slopes $0 < 1$ on the edge $uv$ and slopes $-1 < 0$ on the edge $vw$, in the direction of the arrows. Let $M$ be the sub-semimodule of $\Rat(D, \fS)$ made up of all functions which are symmetric with respect to $v$. We define $\ssrho_v$ by the array $\begin{pmatrix}
				0 & \color{blue} 0 \\
				\color{blue} 1 & 0
			\end{pmatrix}$.
			Here, the jumps of $\ssrho_v$ are depicted in blue. We endow every other point of $\Gamma$ with the standard rank function. The pair $(D, M)$ is a $\g^1_2$.
		\end{example}
		
		In the next two sections, we provide a classification of $\g^1_d$'s.


\section{Reduced divisors} \label{sec:reduced_divisors}
	We establish an extension of the machinery of reduced divisors to linear series. The results of this section are valid without extra effort for the linear system $|(D, M)|$ associated to a pair $(D, M)$ consisting of a divisor $D$ of degree $d$ and an admissible semimodule $M$ of rank $r$, so we present the results in this generality. Moreover, replacing $(D, M)$ with a linearly equivalent admissible pair, we can assume for the full section that $M$ is effective.
	
	\subsection{Reduced divisors in the chip-firing context}
		
		We briefly recall the definition of reduced divisors in the ``chip-firing'' context. In terms of the chip-firing game, the $x$-reduced divisor is obtained from $D$ by firing chips the closest possible to $x$.
		
		More formally, a \emph{cut} $X$ in a metric graph $\Gamma$ is a compact subset of $\Gamma$ with finitely many connected components. The (finite) set of boundary points of $X$, denoted by $\partial X$, is the set of all points of $X$ which are in the closure of the complement of $X$ in $\Gamma$. For a point $x \in \partial X$, we denote by $\out_X(x)$ the set of all outgoing branches (given by unit tangent vectors) at $x$ from $X$, and by $\outval_X(x)$ the number of such branches.
		
		\begin{defi}[Reduced divisors in the classical context] \label{def:reduced_divisor_classical}
			Let $x$ be a point of $\Gamma$. A divisor $D$ on $\Gamma$ is \emph{$x$-reduced} if the following two conditions are met:
			
			\begin{itemize}
				\item[(1)] For each $y \in \Gamma \setminus \{x\}$, we have $D(y) \geq 0$, and
				
				\item[(2)] For each cut $X \subset \Gamma \setminus \{x\}$, there exists a point $y \in \partial X$ with $D(y) < \outval_X(y)$. \qedhere
			\end{itemize}
		\end{defi}
		
		For each divisor $D$ and each point $x$ in $\Gamma$, there exists a unique $x$-reduced divisor linearly equivalent to $D$, which we denote by $D_x$. For graphs, this was proved by Baker and Norine in~\cite{baker2007riemann}. The extension to metric graphs was given in~\cite{MZ08}, see also~\cite{amini2013reduced}.
	
	\subsection{Reduced divisors in linear systems} \label{subsec:definition_reduced_divisors}
		
		We need the following definition.
		
		\begin{defi} \label{def:reduced_f}
			Let $M \subseteq \Rat(D)$ be an admissible semimodule of rank $r$, and let $v$ be a point of $\Gamma$. We define the rational function $\ssf_v^M$, denoted $\ssf_v$ when the context is clear, by
			\[ \ssf_v^M(x) \coloneqq \inf_{g \in M} \mleft[g(x) - g(v)\mright] = \inf_{f \in \ssM_v} f(x) \]
			for every point $x$ of $\Gamma$.
		\end{defi}
		
		\noindent (Recall that $M_v$ is the set of all functions of $M$ which vanish at $v$. In particular, $\ssf_v(v) = 0$.)
		
		\begin{prop} \label{prop:f_v}
			The function $\ssf_v^M$ is well-defined and belongs to $M_v$.
		\end{prop}
		
		\begin{proof}
			We can show that $\ssf^M_v$ can be written as the uniform limit of a decreasing sequence of functions in $M_v$, essentially following the proof of Dini's theorem. Since the slopes of functions in $M_v$ are bounded, this set is uniformly Lipschitz, so $\ssf^M_v$ is itself Lipschitz and therefore continuous. To conclude, since $M$ is closed, $\ssf^M_v$ belongs to $M$.		
		\end{proof}
		
		\begin{remark} \label{rem:reduced_function_minimal_slopes}
			A very special case of Theorem~\ref{thm:jumps_realized} implies that $\ssf^M_v$ takes minimum slopes around $v$: for each $\nu \in \T_v(\Gamma)$, $\mathrm{slope}_\nu(\ssf_v^M) = \sss_0^\nu$, see the discussion right after that theorem. \qedhere
		\end{remark}
		
		\begin{defi}[Reduced divisor] \label{def:reduced_divisor}
			Let $(D, M)$ be a pair consisting of a divisor $D$ of degree $d$ and an admissible semimodule $M \subseteq \Rat(D)$. The effective divisor defined by
			\[ \ssD_v^M \coloneqq D + \div(\ssf_v^M), \]
			simply denoted by $\ssD_v$ when $M$ is contextually clear, is called the \emph{$v$-reduced divisor linearly equivalent to $D$ in the linear system $|(D, M)|$}. Denoting by $\fS$ the slope structure underlying $M$, we denote the slope structure $\fS + \div(\ssf_v)$ by $\fS_v$. The effective pair $(\ssD_v, \fS_v)$ is then equivalent to the pair $(D, \fS)$. Finally, $\ssf_v$ gives rise to a modification of $M$, which is denoted by $M(-\ssf_v)$ according to Definition~\ref{nota:modification_semi_module}, and which is effective.
		\end{defi}
		
		We will see in Proposition~\ref{prop:reduced_divisor_unsaturated_cut} below that $\ssD_v$ is the \emph{unique} effective divisor in the linear system $|(D, M)|$ without unsaturated cuts (see Definition~\ref{def:unsaturated_cut}).

	\begin{example} \label{ex:finiteness_condition}
		The following example shows that the closedness condition on $M$ is necessary for Proposition~\ref{prop:f_v} to hold, and so for Definition~\ref{def:reduced_divisor} to make sense. Consider the metric graph $\Gamma$ and crude linear series of degree four and rank one defined at the end of Example~\ref{ex:second_grd}.
		\begin{figure}[h!]
			\centering
			\begin{tikzpicture}
				\coordinate (A) at (0,0);
				\coordinate (B) at (1.5,0);
				\coordinate (C) at (4,0);
				
				\draw[postaction=decorate,decoration={markings,mark=at position 1/2 with {\arrow[red]{>}; }}]
				(A) -- (B);
				
				\draw[postaction=decorate,decoration={markings,mark=at position 1/2 with {\arrow[red]{>}; }}]
				(B) -- (C);
				
				\foreach \c in {A, B, C} {
					\fill (\c) circle (2.5pt);
					}
				
				\node[left] at (A) {$u$};
				\node[below] at (B) {$v$};
				\node[right] at (C) {$w$};
			\end{tikzpicture}
		\end{figure}
		
		We allow slopes $0 < 1$ on the edge $uv$ and slopes $-1 < 0$ on the edge $vw$, in the direction of the arrows. We set $D = (u) + (w)$ and we define $\ssrho_v$ by the array $\begin{pmatrix}
			0 & \color{blue} 0 \\
			\color{blue} 1 & 0
		\end{pmatrix}$.
		Here, the jumps of $\ssrho_v$ are depicted in blue. This makes $(D, \fS)$ a crude linear series of degree two and rank one. We take $M=\Rat(D, \fS)$. The function $\ssf_v$ 
		\begin{center}
			\begin{tikzpicture}
				\draw[->] (0,0) -- ({1.6},0);
				\draw[->] (0,-1) -- (0,1/2);
				\draw (0,0) node[above left]{$u$} node{$\bullet$};
				\draw (1/2,0) node[above]{$v$} node{$\bullet$};
				\draw (1.33,0) node[above]{$w$} node{$\bullet$};
				\draw[domain=0:1.33] plot(\x, {min(\x-1/2, -\x+1/2)});
			\end{tikzpicture}
		\end{center}
		takes slopes $1$ and $-1$ on the segments $[u, v]$ and $[v, w]$, respectively, and belongs to $\Rat(D, \fS)$. But if we choose any $x \neq v$, then $\ssf_x \notin \Rat(D, \fS)$.
		\begin{center}
			\begin{tikzpicture}[scale=0.7]
				\draw[->] (0,0) -- (4.9,0);
				\draw[->] (0,-2.8) -- (0,3/2);
				\draw (0,0) node[above left]{$u$} node{$\bullet$};
				\draw (4,0) node[above]{$w$} node{$\bullet$};
				\draw[very thick,color=blue] (0,-1) -- (1,0) -- (3/2,0) -- (4,-5/2) node[midway,below]{$\boldsymbol{\ssf_x}$};
				\draw[color=red] (3/2,0) -- (5/2,0);
				\draw[color=red] (5/2,0) -- (4,-3/2);
				\draw[color=red] (4,-3/2) node[right,scale=0.7]{$\ssf_1$};
				\draw[color=red] (2,0) -- (4,-2);
				\draw[color=red] (4,-2) node[right,scale=0.7]{$\ssf_2$};
				\draw[color=red] (1.8,0) -- (4,-2.2);
				\draw[color=red] (4,-2.2) node[right,scale=0.7]{$\ssf_3$};
				\draw[color=red] (1.6,0) -- (4,-2.4);
				\draw[color=blue] (1,0) node[above]{$x$} node{$\bullet$};
				\draw (3/2,0) node[above]{$v$} node{$\bullet$};
			\end{tikzpicture}
		\end{center}
		
		Indeed, consider a sequence of functions $(\ssf_i)_i$ which coincide with $\ssf_x$ on the segment $[u, v]$ and whose graphs on the segment $[v, w]$ are represented in thin red in the above figure, converging uniformly to $\ssf_x$ as $i \to +\infty$. All the functions $\ssf_i$ belong to $\ssub{\Rat(D, \fS)}!_x$; but $\ssf_x$ does not itself belong to $\Rat(D, \fS)$ since it does not take symmetrical slopes around $v$, a mandatory condition given the expression of $\ssrho_v$. In particular, this shows that $\Rat(D, \fS)$ is not closed in $\Rat(D)$.
	\end{example}
	
	\subsection{Coefficient at the base-point}
		
		We have the following useful result.
		
		 \begin{prop} \label{prop:coeff}
			Let $(D, M)$ be pair consisting of a divisor $D$ and an admissible semimodule $M \subseteq \Rat(D)$ of rank $r$. Let $\fS$ be the underlying slope structure. For every point $x \in \Gamma$, we have $\ssD_x(x) = D(x) - \sum_{\nu \in \T_x(\Gamma)} \sss_0^\nu.$
			In addition, this quantity is greater than or equal to $r$.
		\end{prop}
		
		\begin{proof}
			This is a direct consequence of Remark~\ref{rem:reduced_function_minimal_slopes}.
		\end{proof}
	
	\subsection{Unsaturated cuts}
		
		We provide an alternative characterization of reduced divisors.
		
		\begin{defi}[Unsaturated cut with respect to an admissible pair] \label{def:unsaturated_cut}
			Let $(D, M)$ be a pair consisting of a divisor $D$ of degree $d$ and an admissible semimodule $M\subseteq \Rat(D)$. Let $\fS$ be the underlying slope structure. Let $v$ be a point of $\Gamma$. Consider a cut $X$ in $\Gamma$ and assume that the point $v$ does not belong to $X$.
			
			We say that $X$ is \emph{unsaturated with respect to $v$ and $M$} if, for a sufficiently small $\varepsilon > 0$, there exists a function $f \in M$ which satisfies the following properties:
			
			\begin{itemize}
				\item for every point $x \in \partial X$, $f$ is linear of positive slope $\sss^\nu > 0$ on a small segment $\ssI_\nu$ on each adjacent outgoing branch $\nu \in \out_X(x)$;
				
				\item $f$ is identically equal to $-\varepsilon$ on $X$; and
				
				\item $f$ is zero everywhere else.
			\end{itemize}
			
			In this case, we say that $f$ \emph{fires} the unsaturated cut $X$ \emph{at level $\varepsilon$}. Otherwise, if no such $f$ exists, $X$ is called \emph{saturated}.
		\end{defi}
		
		Note that if $f \in M$ fires the cut $X$ at level $\varepsilon$ for some $\varepsilon > 0$, then for every $\varepsilon' \in (0, \varepsilon)$, the element $f' \coloneqq \min(f + \varepsilon - \varepsilon', 0)$ of $M$ fires $X$ at level $\varepsilon'$. Therefore, the saturation property for cuts is well-defined.
		
		\begin{remark} \label{rem:unsaturated_cut_vectors_slopes}
			If $X$ is an unsaturated cut with respect to $v$ and $M$, then for every point $x \in \partial X$, there exists an element $\unds(x) = \prod_{\nu \in \T_x(\Gamma)} \sss^\nu \in \ssS^x$ such that $\sss^\nu \geq 0$ for all $\nu \in \T_x(\Gamma)$, with equality $\sss^\nu = 0$ if and only if $\nu \in \T_x(\Gamma) \setminus \out_X(x)$, and, in addition,
			\begin{equation} \label{eq:inequality_unsaturated_cut_slopes}
				D(x) - \sum_{\nu \in \T_x(\Gamma)} \sss^\nu \geq 0.
			\end{equation}
			Such vectors of slopes $\unds(x)$ are provided by the vector of slopes $\partial_x(f)$ of any function $f \in M$ which fires the unsaturated cut $X$. The data of these vectors of slopes, together with a small number $\varepsilon > 0$, determine the function $f$ entirely.
		\end{remark}
		
		\begin{remark} \label{rem:unsaturated_cut_boundary_finite}
			Inequality~\eqref{eq:inequality_unsaturated_cut_slopes} in the previous remark implies that each point $x$ of the boundary $\partial X$ of the cut $X$ is in the support of the divisor $D$. Since the support of $D$ is finite, for given $v$ and $M$, there are finitely many unsaturated cuts with respect to $v$ and $M$.
		\end{remark}
		
		\begin{prop}[Characterization of reduced divisors by unsaturated cuts] \label{prop:reduced_divisor_unsaturated_cut}
			Let $(D, M)$ be a pair consisting of a divisor $D$ of degree $d$ and an admissible semimodule $M \subseteq \Rat(D)$. Then, $D$ is $v$-reduced if, and only if, there is no unsaturated cut with respect to $v$ and $M$.
		\end{prop}
		
		\begin{proof}
			We prove the equivalence of the negations.
			
			First, assume that there is an unsaturated cut $X$ with respect to $v$ and $M$. By definition, this means that $v \notin X$ and that there exists $\varepsilon > 0$ and a function $f \in M$ which fires $X$ at level $\varepsilon$. Obviously, $f \leq 0$ and $f$ takes negative values on a non-empty set of points. Moreover, if $\varepsilon$ is small enough, then $f(v) = 0$. This shows that $D$ is not $v$-reduced.
			
			To prove the other direction, assume that $D$ is not $v$-reduced. Then, there exists a non-constant function $f \in M \setminus \{0\}$ such that $f(v) = 0$ and $f(x) \leq 0$ for every $x \in \Gamma$. Let $X$ be the set of points of $\Gamma$ where $f$ takes its minimum value. This is a compact set that does not contain $v$ with finitely many connected components, given that $f$ is piecewise linear. Therefore, $X$ is a cut, and the properties of $f$ imply that it is unsaturated with respect to $v$ and $M$, as required.
		\end{proof}
	
	\subsection{Behavior of reduced divisors with respect to the base point} \label{subsec:behavior_reduced_divisors_with_respect_base_point}
		
		We describe how reduced divisors behave under an infinitesimal change of the base point. This will be used to prove Theorem~\ref{thm:red_continuous}.
		
		Let $(D, M)$ be a pair consisting of a divisor $D$ of degree $d$ and an admissible semimodule $M \subseteq \Rat(D)$ of rank $r$ with underlying slope structure $\fS$. Let $v$ be a point of $\Gamma$, $\nu \in \T_v(\Gamma)$, and $e$ the oriented edge of $\Gamma$ parallel to $\nu$, all fixed for the remainder of this section. We give an explicit description of $\ssD_u$ for $u$ in a small segment $\ssI^e \subset e$ with an endpoint equal to $v$.
		
		Replacing $(D, \fS)$ with the linearly equivalent $(\ssD_v, \fS_v)$ and $M$ with $M(-\ssf_v)$, we can assume that $D$ is $v$-reduced. Note that $M$ remains effective. In particular, the zero function belongs to $M$. Also, adapting the combinatorial model accordingly, we suppose that the slope structure $\fS$ is constant on every edge. Now, there exists a sufficiently small segment $\ssI^e$ on $e$ adjacent to $v$ which does not contain any point of the support of $D$ apart from $v$. For a point $u \neq v$ on this segment, we have $D(u) = 0$ and $\ssD_u(u) > 0$, and thus $D \neq \ssD_u$. We infer that $D$ is not $u$-reduced. It follows from Proposition~\ref{prop:reduced_divisor_unsaturated_cut} that there exist unsaturated cuts with respect to $u$ and $M$. Since $D$ is $v$-reduced, we infer that every such cut $Y$ contains $v$. In addition, $\nu$ must belong to $\out_Y(v)$, and therefore $v \in \partial Y$, since otherwise, the boundary of that cut would contain a point of $e$ between $v$ and $u$ (see Remark~\ref{rem:unsaturated_cut_boundary_finite}). This would be impossible by the assumption made on the support of $D$. We have proved
		
		\begin{claim} \label{claim:point_boundary_unsaturated_cut}
			For every unsaturated cut $Y$ with respect to $u$ and $M$, we have $v \in \partial Y$.
		\end{claim}
		
		Since $D$ is $v$-reduced, we have $\sss_0^\nu = 0$ for every $\nu \in \T_v(\Gamma)$ (see Remark~\ref{rem:reduced_function_minimal_slopes}). It follows that $\sss_r^\nu > 0$, and, by Proposition~\ref{prop:coeff} and the definition of slope structures, the coefficient of $\ssD_u$ at $u$ is precisely equal to $\sss_r^\nu + \sss_0^\nu = \sss_r^\nu$, given that $-\sss_r^\nu$ and $\sss_0^\nu = 0$ are the smallest possible slopes at $u$ in the direction of $v$ and away from $v$, respectively. We now claim
		
		\begin{claim} \label{claim:existence_unsaturated_cut_maximal_slope}
			There exists an unsaturated cut $Y$ with respect to $u$ and $M$ for which we can choose a function $f \in M$ firing $Y$ at some level $\varepsilon > 0$ with the additional property that its slope along the tangent vector $\nu$ is the maximal possible slope $\sss_r^\nu$.
		\end{claim}
		
		\begin{proof}
			Let $f = \ssf_u$ be the element of $M$ with $D + \div(f) = \ssD_u$ and $f(u) = 0$, as in Definition~\ref{def:reduced_f}. Since $-\sss_r^\nu$ is the minimum slope at $u$ along $e$ in the direction of $v$, Remark~\ref{rem:reduced_function_minimal_slopes} implies that $f$ takes slope $-\sss_r^\nu$ away from $u$ on a sufficiently small segment included in $[v, u]$ containing $u$. The fact that $\ssD_u$ is effective then implies that $f$ has constant slope along $[v, u]$. Consequently, $f$ has slope $\sss_r^\nu$ in the direction of $u$ on the whole segment $[u, v]$, and in particular along the tangent vector $\nu$.
			
			Now, let $\eta$ be the distance between $v$ and $u$ on $\ssI^e$, and set $\varepsilon \coloneqq \sss^\nu_r \eta$. Then, we have $f(v) = -\varepsilon$. Also, note that, by the hypothesis made on $D$, $\fS$ and $M$, we have $\ssf_v \equiv 0$ on $\Gamma$. Therefore, following Definition~\ref{def:reduced_f}, we get $f \geq f(v) + \ssf_v = -\varepsilon$ on $\Gamma$. Consequently, $-\varepsilon = f(v)$ is the minimum of $f$ on $\Gamma$.
			
			Consider the set $Y$ of all points of $\Gamma$ where $f$ takes this minimum value. $Y$ is a compact set with finitely many connected components. Moreover, $Y$ does not contain $u$ and, for a sufficiently small $\varepsilon' \in (0, \varepsilon]$, the function $f' \coloneqq \min(f + \varepsilon - \varepsilon', 0)$ fires $Y$ at level $\varepsilon'$. 
			Since in addition the slope of $f$ (and thus the slope of $f'$) along $\nu$ is equal to $\sss_r^\nu$, the unsaturated cut $Y$ has the desired properties.
		\end{proof}
		
		Now consider the family $\mathcal X$ of all the unsaturated cuts $Y$ with respect to $u$ and $M$ such that $Y$ verifies the properties of Claim~\ref{claim:existence_unsaturated_cut_maximal_slope}. Let $X \coloneqq \bigcup_{Y \in \mathcal X} Y$.
		
		Since there are finitely many unsaturated cuts with respect to $u$ and $M$ (see Remark~\ref{rem:unsaturated_cut_boundary_finite}), the family $\mathcal X$ is finite. Therefore, $X$ is still compact and has finitely many connected components. i.e., it is itself a cut.
		
		\begin{fact} \label{fact:set_uns_cuts}
			We notice that if we choose another point $u' \neq v$ on the segment $[v, u]$, a cut $Y$ is unsaturated with respect to $u'$ if, and only if, it is unsaturated with respect to $u$. Verifying the stronger property of Claim~\ref{claim:existence_unsaturated_cut_maximal_slope} also remains unchanged. This shows that $\mathcal X$ (and thus $X$) does not change when we choose $u$ to be even closer to $v$, a fact that will be used in the discussions below and in the proof of Theorem~\ref{thm:exprd}.
		\end{fact}
		
		\begin{claim} \label{claim:union_cuts_is_cut}
			Notation as above, the cut $X$ itself belongs to $\mathcal X$.
		\end{claim}
		
		\begin{proof}
			To prove that $X$ is unsaturated with respect to $u$ and $M$, we first note that $u \notin X$.
			
			Then, denote by $\ssY_1, \dots, \ssY_N$ the elements of $\mathcal X$. By assumption, for every $i \in \{1, \dots, N\}$, there exists a function $\ssf_i \in M$ such that for $\ssvarepsilon_i > 0$ small enough, $\ssf_i$ fires the unsaturated cut $\ssY_i$ at level $\ssvarepsilon_i$, and such that the slope of $\ssf_i$ along $\nu$ is equal to $\sss_r^\nu$. For a given $i$ and for a fixed $\ssvarepsilon_i > 0$, there are finitely many such functions $\ssf_i \in M$, one for each choice of joint outgoing slopes (see Remark~\ref{rem:unsaturated_cut_vectors_slopes}). By choosing the minimum of these possible choices, we ensure that $\ssf_i$ takes the minimal possible set of joint slopes on edges leaving $\ssY_i$, under the constraint, if $v \in \partial \ssY_i$, of having slope $\sss_r^\nu$ on $\nu$. Furthermore, we replace all the numbers $\ssvarepsilon_i$ with $\varepsilon \coloneqq \min_i \ssvarepsilon_i$.
			
			For each $i \in \{1, \dots, N\}$, each $x \in \partial \ssY_i$, and each $\nu' \in \out_{\ssY_i}(x)$, we call $\ssI_{\nu'}$ the small segment parallel to $\nu'$ on which $\ssf_i$ takes a positive slope away from $\ssY_i$. We then reduce $\varepsilon > 0$ again so that the outgoing segments $\ssI_{\nu'}$ are essentially disjoint and therefore the functions $\ssf_i$ fire in a ``decoupled'' way. More precisely, up to reducing $\varepsilon$ further, we can assume that for every $i, j \in \{1, \dots, N\}$, every $x \in \partial \ssY_i$, $y \in \partial \ssY_j$ and every $\nu' \in \out_{\ssY_i}(x)$, $\nu'' \in \out_{\ssY_j}(y)$, if $\mathring \ssI_{\nu'} \cap \mathring \ssI_{\nu''} \neq \varnothing$, then $x = y$ (and thus one of $\ssI_{\nu'}, \ssI_{\nu''}$ contains the other).
			
			For a sufficiently small $\varepsilon$ chosen this way, the function $\ssf^\varepsilon \coloneqq \min_i \ssf_i$ belongs to $M$ and fires $X$ at level $\varepsilon$. Indeed,
			
			\begin{itemize}
				\item on $X = \bigcup_i \ssY_i$, $\ssf^\varepsilon$ takes value $-\varepsilon$, and away from $X$ and the union of $\ssI_{\nu'}$, for $x \in \partial \ssY_i$, $i = 1, \dots, N$, and $\nu' \in \out_{\ssY_i}(x)$, it is identically zero.
				
				\item on a segment of the form $\ssI_{\nu'}$ based at some $x \in \partial X$, the slope of $\ssf^\varepsilon$ away from $X$ is the minimum of the slopes of the functions $\ssf_i$ for $i$ such that $x \in \partial \ssY_i$.
			\end{itemize}
			
			This shows that $X$ is unsaturated. Moreover, the slope of every $\ssf_i$ along $\nu$ is $\sss_r^\nu$, which concludes. Note that $f$ inherits the slope minimality properties of the functions $\ssf_i$, which will be useful in the proof of Theorem~\ref{thm:exprd}.
		\end{proof}
		
		Now, as in the proof of Claim~\ref{claim:existence_unsaturated_cut_maximal_slope}, let $\eta$ be the distance between $v$ and $u$ on $\ssI^e$, and set $\ssvarepsilon_0 \coloneqq \sss^\nu_r \eta$. We look at which number between $\ssvarepsilon_0$ and the number $\varepsilon$ appearing in the proof of Claim~\ref{claim:union_cuts_is_cut} is smaller. If $\ssvarepsilon_0 < \varepsilon$, we freely reduce $\varepsilon$ to $\ssvarepsilon_0$ in the proof of Claim~\ref{claim:union_cuts_is_cut}. If, on the contrary, $\varepsilon < \ssvarepsilon_0$, we move $u$ closer to $v$ so that both numbers are equal; we can do this properly and without self-reference because, as $u$ moves closer to $v$, the family $\mathcal X$ remains constant, by Fact~\ref{fact:set_uns_cuts}. At this point, $X$ is an unsaturated cut with respect to $u$ and $M$, and there exists a function $\ssf^u \coloneqq \ssf^\varepsilon \in M$, constructed in the proof of Claim~\ref{claim:union_cuts_is_cut}, with fires $X$ at level exactly $\varepsilon = \ssvarepsilon_0$, with constant slope $\sss_r^\nu$ on the whole segment $[v, u]$. Moreover, by construction, among the functions with the previous properties, $\ssf^u$ has the smallest possible set of joint slopes on outgoing tangent vectors.
		
		\begin{remark} \label{rem:continuity}
			Given the way it was constructed in the proof of Claim~\ref{claim:union_cuts_is_cut}, and thanks to Fact~\ref{fact:set_uns_cuts}, the function $\ssf^u$ is entirely determined by $\ssvarepsilon_0$, and its non-zero slopes are independent of $u$. More precisely, if $u' \in (v, u]$, denoting by $\eta'$ the distance between $v$ and $u'$, and letting $\ssvarepsilon_0' \coloneqq \sss_r^\nu \eta'$, then the function $\ssf^{u'}$ is given by $\ssf^{u'} = \min(\ssf^u + \ssvarepsilon_0 - \ssvarepsilon_0', 0)$. This shows that $\ssf^u$ depends continuously on $u$.
		\end{remark}
		
		We can now formulate the main theorem of this section.
		
		\begin{thm} \label{thm:exprd}
			The $u$-reduced divisor with respect to $M$ is $\ssD_u = D + \div(\ssf^u)$.
		\end{thm}
		
		\begin{proof}
			Let $\ssf_u$ be the rational function in $M$ which defines the $u$-reduced divisor $\ssD_u$, and which takes value zero at $u$ (see Definition~\ref{def:reduced_f}). We will prove that in fact $\ssf_u = \ssf^u$. Note that $\ssf_u$ verifies the following properties.
			
			\begin{itemize}
				\item $\ssf_u$ is linear on the segment $[v, u]$ with slope $\sss_r^\nu$.
			\end{itemize}
			
			This follows from Remark~\ref{rem:reduced_function_minimal_slopes} and the fact that $\ssD_u$ is effective.
			
			\begin{itemize}
				\item $\ssf_u$ takes it minimum value at $v$, and its maximum value at $u$.
			\end{itemize}
			
			That $\ssf_u$ takes its maximum value at $u$ comes from the definition of the $u$-reduced divisors, given that $0 \in M$. That $\ssf_u$ takes its minimum value at $v$ was shown in the proof of Claim~\ref{claim:existence_unsaturated_cut_maximal_slope}.
			
			\begin{itemize}
				\item Let $\ssX_0$ be the set of all points where $\ssf_u$ takes its minimum value. Then, $\ssX_0 = X$.
			\end{itemize}
			
			Obviously, $\ssX_0$ is an unsaturated cut with respect to $u$ and $M$ which in addition verifies the property of Claim~\ref{claim:existence_unsaturated_cut_maximal_slope} (it is the set $Y$ in the proof of that statement). Therefore $\ssX_0 \in \mathcal X$. This shows by Claim~\ref{claim:union_cuts_is_cut} that $\ssX_0 \subseteq X$. The reverse inclusion comes from the definition of reduced divisors, which implies $\ssf_u \leq \ssf^u$, globally on $\Gamma$. This shows that $\ssX_0 = X$.
			
			\begin{itemize}
				\item Up to reducing $\eta$ further (i.e., up to moving $u$ even closer to $v$), $\ssf_u$ and $\ssf^u$ coincide on all segments of length $\eta$ around $X$.
			\end{itemize}
			
			Note that $\ssf_u = \ssf^u$ on $X$ by the previous claim. Then, the slope minimality property of $\ssf^u$ (see the proof of Claim~\ref{claim:union_cuts_is_cut}) ensures, on the one hand, that the slopes taken by $\ssf^u$ on tangent vectors leaving $X$ do not exceed those taken by $\ssf_u$. On the other hand, the inequality $\ssf_u \leq \ssf^u$ holds on $\Gamma$, as previously noted. All in all, we can move $u$ closer to $v$ (and thus reduce $\eta$), without changing $X$ (see Fact~\ref{fact:set_uns_cuts}), and keeping the properties of $\ssf_u$ and $\ssf^u$ established so far in the current proof, in such a way that $\ssf_u$ and $\ssf^u$ coincide on all segments of length $\eta$ leaving $X$.
			
			\begin{itemize}
				\item $\ssf_u = \ssf^u$ everywhere.
			\end{itemize}
			
			To finish the proof of the theorem, let $X'$ be the closure of the complement of $X \cup \bigcup_{x \in \partial X, \nu' \in \out_X(x)} \ssI_{\nu'}$. We need to show that $\ssf_u = \ssf^u$ on $X'$. In other words, we need to show that $\ssf_u \equiv 0$ on $X'$. Suppose this is not the case, and consider the minimum locus $Y$ of $\ssf_u$ on $X'$. Note that $Y$ lies in the interior of $X'$, i.e., $Y \cap \partial X' = \varnothing$, and $v \notin Y$. 
			This shows that $Y$ is an unsaturated cut with respect to $v$ and $M$, contradicting Proposition~\ref{prop:reduced_divisor_unsaturated_cut}, given that $D$ is $v$-reduced.
		\end{proof}
		
		\begin{remark}
			Note that in the proof of Theorem~\ref{thm:exprd}, we did not refer to the admissibility of $M$, but this does not contradict Example~\ref{ex:finiteness_condition}. Indeed, Proposition~\ref{prop:f_v} states that, assuming admissibility, $\ssD_v$ exists for all $v$, whereas Theorem~\ref{thm:exprd} only implies that if $\ssD_v$ \emph{does exist} for some $v$, then $\ssD_u$ exists for all $u$ in a neighborhood of $v$, and behaves as stated in the theorem. In other words, the set of $u$ for which $\ssD_u$ exists is an open subset of $\Gamma$.
		\end{remark}
	
	\subsection{Continuity of the map to $|(D, M)|$ defined by reduced divisors}
		
		Theorem~\ref{thm:exprd} has the following direct consequence, which will be crucial in the next section.
		
		\begin{thm} \label{thm:red_continuous}
			Let $D$ be a divisor of degree $d$ and $M \subset \Rat(D)$ an admissible semimodule on $\Gamma$. The map 
			\[ \Red \colon \Gamma \rightarrow |(D, M)| \subseteq \mathrm{Sym}^d(\Gamma), \qquad \qquad v \mapsto \ssD_v, \]
			is continuous and non-contracting.
		\end{thm}
		
		\begin{proof}
			By Proposition~\ref{prop:coeff}, the coefficient of $\ssD_v$ at $v$ is precisely $D(v) - \sum_{\nu \in \T_v(\Gamma)} \sss_0^\nu > 0$ for every point $v$ of $\Gamma$. This obviously shows that $\Red$ cannot be constant on any segment of positive length, proving that it is non-contracting. Furthermore, it follows from Theorem~\ref{thm:exprd} and Remark~\ref{rem:continuity} that the function $\ssf_u = \ssf^u$ depends continuously on the point $u \in \Gamma$. Since $\ssD_u = D + \div(\ssf^u)$, the continuity of $\Red$ follows.
	 \end{proof}
	 
	 \begin{remark}
	 	The above result is a generalization of \cite[Theorem~3]{amini2013reduced}. Note that the map $\Red$ is also integer affine linear with respect to the natural integer affine structures on $\Gamma$ and $\mathrm{Sym}^d(\Gamma)$ as described in \cite[Section~2]{amini2013reduced}. Thus, the image of $\Red$ is a metric graph. In the case $(D, M)$ is a $\g^1_d$, we will see in Proposition~\ref{prop:T_no_cycles} that the image is a metric tree.
	\end{remark}


\section{Classification of $\g^1_d$'s} \label{sec:g1d}
	
	In this section, we consider the case $r = 1$, and prove, roughly speaking, that the data of a $\g^1_d$ on $\Gamma$ is equivalent to the data of a finite harmonic map to a metric tree (see Theorem~\ref{thm:g1d} below for a precise statement). Then, we formulate a smoothing theorem for combinatorial $\g^1_d$'s. In this regard, $\g^1_d$'s on metric graphs are well-behaved, and our theorem can be regarded as a generalization of the Eisenbud--Harris smoothing result for their limit $\g^1_d$'s~\cite[Proposition~3.1]{EH86}.
	
	\smallskip
	
	Let $(D, M)$ be a $\g^1_d$ on $\Gamma$ with underlying slope structure $\fS$ of width one. By Remarks ~\ref{rem:linear_series_independence}
	and~\ref{rem:semi_module_make_effective}, we can assume that $(D, M)$ is effective. This implies that for every point $x$ and every outgoing tangent direction $\nu \in \T_x(\Gamma)$ at $x$, one of the two integers $\sss_0^\nu$ or $\sss_1^\nu$ is equal to zero.
	
	Assume that $(D, \fS)$ is defined on a model $G = (V, E)$ of $\Gamma$.
	
	\begin{defi}[Orientation associated to a $\g^1_d$]
		Notation as above, we define an orientation of $G$ in such a way that the edge $\{u, v\}$ gets orientation $uv$ if $\sss_0^{uv} = 0 < \sss_1^{uv}$.
	\end{defi}
	
	Let $\rho$ be a rank function on $\hcube{1}{\grdeg}$. The point $\underline 0$ is the only point of $\hcube{1}{\grdeg}$ of rank one (Remark~\ref{rem:link_coordinate_rank}). Besides, the set $\ssJ_\rho$ of jumps of $\rho$ contains the point $\underline 0$ (because $\rho(\unde_i) = r - 1 = 0$ for all $i$) and every other point $\unda \neq \underline 0$ has at least one coordinate equal to one. For each $\unda \in \ssJ_\rho$, denote by $\ssP{\unda}$ the subset of $\{1, \dots, \grdeg\}$ consisting of all the indices $i$ with $\ssa_i = 1$ (the support of $\unda$). Denote by $\ssmathcalP_\rho$ the collection of all sets $\ssP{\unda}$ for $\unda \in \ssJ_\rho \setminus \{\underline 0\}$. By Lemma~\ref{lem:partition} and Remark~\ref{rem:partition_1}, $\ssmathcalP_\rho$ provides a partition of $\{1, \dots, \grdeg\}$.
	
	\begin{prop} \label{prop:rank_function_1_representable}
		Every rank function on $\hcube{1}\grdeg$ is representable, and the field can be chosen to be of characteristic zero.
	\end{prop}
	
	\begin{proof}
		Chose $\k$ to be any infinite field of characteristic zero. Let $\ssP{\rho} = \{\ssP 1, \dots, \ssP s\}$ be the partition of $\{1, \dots, \grdeg\}$ as previously described, each $\ssP i$ being a subset of $\{1, \dots, \grdeg\}$. In the plane $\k^2$, let $\ssL_1, \dots, \ssL_s$ be distinct lines. Now let, for each $i \in \{1, \dots, \grdeg\}$, $\ssfilt_i^1 \coloneqq \ssL_{\tau(i)}$, where $\tau(i) \in \{1, \dots, \grdeg\}$ is the integer such that $i \in \ssP{\tau(i)}$. This, in turn, defines $\grdeg$ complete flags in $\k^2$, completing the lines with $\k^2$ and $(0)$. The collection of flags $\ssfilt_1^\bullet, \dots, \ssfilt_\grdeg^\bullet$ represents $\rho$.
	\end{proof}
	
	For every point $p$ of $\Gamma$ (resp. $G$) of valence $d_p$, consider the rank function $\ssrho_p$ on $\hcube{1}{d_p}$, and define $\ssmathcalP_p$ as the partition of $\T_p(\Gamma)$ (resp. $\mathbb E_p$) given by the previous proposition. Note that $S^p $ consists of the point $(\sss_0^\nu)_{\nu \in \T_p(\Gamma)}$ (resp. $(\sss_0^e)_{e \in \mathbb E_p}$), and vectors $(\sss_{\ssa_\nu}^\nu)$ (resp. $(\sss_{\ssa_e}^e)$) for $\unda \in \hcube{1}{d_p}$ with $\ssP{\unda} \in \ssmathcalP_p$.
	
	\begin{prop} \label{prop:negative_slopes_partition}
		Let $v$ be a vertex of $G$, and denote by $\ssI_v \subset \mathbb E_p$ the set of all the edges $vw$ with $\sss_0^{vw} < 0$. Then we have $\ssI_v \in \ssmathcalP_v$ or $\ssI_v = \emptyset$.
	\end{prop}
	
	\begin{proof}
		We know by hypothesis that $0 \in M$. Let $\unda$ be the element of $\hcube{1}{\ssd_v}$ such that $\ssdelta_v(0) = \sss_{\unda}$. Then, if $vw$ is an edge incident to $v$, since exactly one of the possible slopes along $vw$ is zero, we have $\sss_0^{vw} < 0 \Leftrightarrow \ssa_{vw} = 1$, which concludes.
	\end{proof}
	
	A consequence of this proposition is the following corollary.
	
	\begin{cor} \label{cor:unique_set_partition_negative_slopes}
		Let $p$ be a point of $\Gamma$ and let $\T_p(\Gamma) = \{\nu_1, \dots, \nu_{\ssd_p}\}$. For each $P \in \ssmathcalP_p$, all the coordinates $\sss_j^{\nu_i}$ with $i \in P$ and $j \in \{0, 1\}$ have the same sign. In addition, there exists at most one $P \in \ssmathcalP_p$ such that all $\sss_0^{\nu_i}$ for $i \in P$ are negative (and, if it exists, this is $\ssI_p$).
	\end{cor}
	
	\begin{propdefi} \label{prop_def:unique_g11_tree}
		Let $T$ be a metric tree. Then, there exists a unique effective $\g^1_1$ on $T$ up to linear equivalence.
	\end{propdefi}
	
	\begin{proof}
		Let $(D, N)$ be an effective $\g^1_1$ on $T$ with underlying slope structure $\fS$. Then, we can assume that $D = (v)$ for some vertex $v$ of the tree. For every point $y$ in the tree, let $\ssf_{v \rightarrow y}$ to be the unique function taking value zero at $v$ with $\div(\ssf_{v \rightarrow y}) + (v) = (y)$. Since $N$ has rank one and $N \subseteq \Rat(D)$, $\ssf_{v \rightarrow y}$ must belong to $N$. This shows that $\Rat(D) = N$. By~\cite{haase2012linear}, $\Rat(D)$ is finitely generated (for every metric graph). It is easy to see that $N$ has tropical rank one and from this, the proposition follows.
	\end{proof}
	
	The following is the main theorem of this section. We recall that a map between metric graphs is called a \emph{morphism} if it is integer affine with respect to the natural integer affine structures of the metric graphs, and it is called \emph{finite} if the preimage of every point is finite.
	
	\begin{thm} \label{thm:g1d}
		Let $(D, M)$ be an effective $\g^1_d$ on a metric graph $\Gamma$ with underlying slope structure $\fS$. We suppose that $D$ is $\ssx_0$-reduced for some point $\ssx_0 \in \Gamma$. Then, we have the following results.
		
		\begin{itemize}
			\item The image of the map $\Red \colon \Gamma \rightarrow |M| \subseteq \mathrm{Sym}^d(\Gamma)$ is a metric tree $T$, and $\Red$ is a finite morphism. Moreover, we have $M = \Red^*(N)$, where $N$ is the semimodule on $T$ defined in Proposition--Definition~\ref{prop_def:unique_g11_tree} using the point $\Red(\ssx_0) \in T$.
			
			\item There exist a tropical modification $\alpha \colon \widetilde \Gamma \longrightarrow \Gamma$ of $\Gamma$ and a finite harmonic morphism $\varphi \colon \widetilde \Gamma \rightarrow T$ of degree $d$ such that $\varphi_{|\Gamma} = \Red$.
		\end{itemize}
		
		Conversely, let $\psi \colon \wtGamma \rightarrow T$ be a finite harmonic morphism from a tropical modification $\wtGamma$ of $\Gamma$ to a metric tree $T$, and denote by $\tau \colon \wtGamma \to \Gamma$ the contraction map. We put the semimodule $N$ on $T$ defined in Proposition--Definition~\ref{prop_def:unique_g11_tree} using the point $\ssy_0 = \psi(\ssx_0)$. Then, $\mleft(\tau\mleft(\psi^{-1}(\ssy_0)\mright), \psi^*(N)\mright)$ is an effective $\g^1_d$ on $\Gamma$ for $d$ equal to the degree of $\psi$.
	\end{thm}
	
	\begin{remark}
		The theorem shows that the combinatorial $\g^1_d$'s are the precise analogues of algebraic geometric linear series of rank one. In the context of linear series on metrized complexes introduced in~\cite{AB15}, a smoothing theorem in rank one was previously obtained by Luo and Manjunath in~\cite{LM18}. However, in that context, there are obstructions to smoothing, and the results are different.
	\end{remark}
	
	Using~\cite[Theorem~3.11]{amini2015lifting2} and~\cite[Theorem~7.7]{amini2015lifting1}, we deduce the following smoothing theorem. See Section~\ref{sec:geometric_tropicalization} for more details about Berkovich curves and tropicalization of linear series.
	
	\begin{thm}[Smoothing theorem~\ref{thm:g1d_smoothing_intro} for $\g^1_d$'s] \label{thm:g1d_smoothing}
		Every effective $\g^1_d$ $(D, M)$, with $M \subseteq \Rat(D)$, on $\Gamma$ is smoothable. That is, there exists a smooth proper curve $X$ over an algebraically closed field with a non-trivial non-Archimedean valuation such that $\Gamma$ is a skeleton of $X^{\an}$, a divisor $E$ on $X$, and a vector subspace $\VS \subseteq H^0(X, \mathcal O(E))$ of rank one on $X$ such that $(D, M)$ is the tropicalization of the $\g^1_d$ $(E, \VS)$ from $X$ to $\Gamma$.
	\end{thm}
	
	\begin{remark}
		This crucially uses the fact that rank functions of rank one are geometric in the terminology of Section~\ref{sec:rank} (see~\cite{amini2015lifting1} and~\cite{amini2015lifting2}). This is applied to every $\ssrho_v$ with $v \in V$ for $G = (V, E)$ a model a $\Gamma$ on which $(D, \fS)$ is defined.
	\end{remark}
	
	The rest of this section is devoted to the proof of Theorem~\ref{thm:g1d}.
	
	\subsection{Proof of Theorem~\ref{thm:g1d}}
		
		To this end, we will define an equivalence relation and a partial preorder on the points of $\Gamma$, show that the equivalence classes correspond to the map $\Red$, and prove that the quotient $\Gamma/\sim$ is a metric tree.
		
		Let $(D, M)$ be an effective $\g^1_d$ with underlying slope structure $\fS$ on $\Gamma$. We define an equivalence relation $\sssim_M$ on the set of points of $\Gamma$ as follows. For two points $x, y \in \Gamma$, we write $x \sssim_M y$ (or simply $x \sim y$) if, for all $f \in M$, we have $f(x) = f(y)$. We also define a partial order by writing $x \ssleq_M y$ (or simply $x \leq y$) if, for all $f \in M$, we have $f(x) \leq f(y)$. Note that, if $x$ and $x'$ belong to the same edge $e$, then $x$ and $x'$ are comparable for $\ssleq_M$: in this case, the comparability is given by the orientation of $e$. We also write $x \ssle_M y$ (or simply $x < y$) the corresponding strict partial preorder, when $x \leq y$ and there exists some $f \in M$ such that $f(x) < f(y)$. The statement $x \ssle_M y$ is equivalent to $x \ssleq_M y$ and $x \ssnsim_M y$.
		
		It is immediate that $\ssleq_M$ is reflexive and transitive, and $\ssle_M$ irreflexive, transitive and asymmetric. Also, $\ssleq_M$, resp. $\ssle_M$, induces a well-defined partial order, resp. a well-defined strict partial order, on the quotient $\Gamma/\sim$.
		
		We first show that $\sim$ is a well-behaved equivalence relation.
		
		\begin{prop} \label{prop:equivalence_classes_finite}
			For all $x \in \Gamma$, the class of $x$ under $\sim$ is finite.
		\end{prop}
		
		\begin{proof}
			Let $x$ and $x'$ belong to a common edge $e$. Suppose without loss of generality that $x \leq x'$. Then, by Theorem~\ref{thm:jumps_realized}, there exists $f \in M$ that takes all minimal slopes around $x'$ (for example $f = \ssf_{x'}$), in particular in the direction of $x$ in which the minimal slope is negative. Thus $f(x) < f(x')$ and $x \nsim x'$. This shows that we can only have $x \sim x'$ if $x$ and $x'$ do not belong to the same edge, from which the result follows.
		\end{proof}
		
		In particular, this shows that $\Red$ is a finite morphism. We next show that the equivalence relation corresponds to the map $\Red$.
		
		\begin{prop} \label{prop:equivalence_reduced_map}
			Let $x, y \in \Gamma$. Then $x \sim y \Longleftrightarrow \Red(x) = \Red(y)$.
		\end{prop}
		
		\begin{proof}
			$\Longrightarrow$ If $x \sim y$, then the sets $\ssM_x$ and $\ssM_y$ coincide, which shows, according to Definition~\ref{def:reduced_f}, that $\ssf_x = \ssf_y$ and, therefore, $\Red(x) = \Red(y)$.
			
			\smallskip
			
			$\Longleftarrow$ If $\Red(x) = \Red(y)$, by Proposition~\ref{prop:functions_equal_divisor}, there exists a constant $c$ such that $\ssf_y = \ssf_x + c$. Suppose that $c > 0$. Then $\ssf_y(x) = \ssf_x(x) + c = c > 0$, which is impossible because for all $z \in \Gamma$, $\ssf_z \leq 0$. So we have, in fact, $\ssf_x = \ssf_y$ by symmetry. Now, let $g$ be any function of $M$ such that $g(x) = 0$. We need to show that $g(y) = 0$. First, we note that
			\[ g(y) \geq \ssf_x(y) = \ssf_y(y) = 0. \]
			Second, let $h \coloneqq g - g(y)$. $h$ belongs to $M$ and verifies $h(y) = 0$, so we know that $h \geq \ssf_y = \ssf_x$, and thus, $-g(y) = h(x) \geq \ssf_x(x) = 0$, which concludes.
		\end{proof}
		
		The map $\Red$ is thus equal to the projection map $\Gamma \rightarrow \Gamma/\sim$. We denote $\Gamma/\sim$ by $T$ and transfer naturally the orientation of the edges of $\Gamma$ to the edges of $T$.
		
		\smallskip

		The goal now is to prove that $T$ is a metric tree. This amounts to showing that $T$ has no cycle ($T$ is already a metric graph). The proof uses the condition that the tropical rank is one.

		\begin{lemma} \label{lem:indegree_two}
			Let $\overline x \in T$. Then $\overline x$ has in-degree at most one in $T$.
		\end{lemma}
		
		\begin{proof}
			Suppose, on the contrary, that $\overline x$ has in-degree at least two, that is, using the orientation of $G$ defined by $\fS$, that two distinct edges $\overline e$ and $\overline e'$ incident to $\overline x$ are oriented toward $\overline x$. There are two cases.
			
			\begin{enumerate}[label=(\roman*)]
				\item $\overline e$ and $\overline e'$ originate from two edges of $G$, $e$ and $e'$ respectively, which are incident to the same preimage $x$ of $\overline x$. 
			\end{enumerate}
			
			We show that this cannot happen by proving that, on the contrary, $e$ and $e'$ will be glued together by $\Red$ in the following way: there exist $\lambda, \varepsilon > 0$ such that
			
			\begin{itemize}
				\item $\varepsilon$ and $\lambda \cdot \varepsilon$ are smaller than $\ell(e), \ell(e')$, and
				
				\item for every choice of $\eta \in (0, \varepsilon)$, for every $y \in e$ at distance $\eta$ from $x$, for $y' \in e'$ at distance $\lambda \cdot \eta$ from $x$, we have $y \sim y'$.
			\end{itemize}
			
			We know that $\ssf_x$ takes all minimum slopes around $x$, and, in particular, takes negative slopes on $e$ and $e'$ away from $x$. These two slopes can be different, but since they are both negative, the image by $\ssf_x$ identifies points of $e$ and $e'$ close to $x$ as described above for some fixed $\varepsilon>0$, with a dilation factor $\sslambda = \frac{\sss_0^e}{\sss_0^{e'}}$. The null function $0$ also identifies $e$ and $e'$ with the same dilation factor (in fact, with any dilation factor).

			Let now $g$ be any function of $M$. Without loss of generality, we assume that $g(x) = 0$. Since $M$ is assumed to have tropical rank one, there exist $c, d \in \R$ such that for all $z \in \Gamma$, the minimum in
			\[ \min(0, \ssf_x(z) + c, g(z) + d) \]
			is attained at least twice. We will show that the image by $g$ identifies $e$ and $e'$ in the same way, getting a contradiction.
			
			We will consider three cases, depending on whether $c$ is negative, positive, or zero.

			Let us first assume that $c$ is negative. This implies that, for $y$ on $e$ or $e'$, $\ssf_x(y) + c < 0$ and thus $g(y) + d = \ssf_x(y) + c$. So $g = \ssf_x$ on $e$ and $e'$ (evaluate at $x$) and thus identifies $e$ and $e'$ close to $x$ with dilation factor $\sslambda$. 
			
			Now, let us assume that $c$ is positive. Then, for $y$ on $e$ or $e'$ close to $x$, $\ssf_x(y) + c$ is still positive and thus $g(y) + d = 0$. In fact, $d$ must be zero, and so $g$ is null close to $x$ on $e$ and $e'$. Therefore, $g$ identifies $e$ and $e'$ with dilation factor $\sslambda$.
			
			The last case is when $c = 0$. For $y$ on $e$ different from $x$, we have $\ssf_x(y) = g(y) + d$ and thus, by continuity, $g$ identifies $e$ and $e'$ with dilation factor $\sslambda$.
			
			We have thus shown that in all cases, $g$ identifies the two edges oriented toward $x$. Since this is true for all $g \in M$, we infer tha $\overline x$ cannot have in-degree at least two.
			
			\begin{enumerate}[label=(\roman*),resume]
				\item $\overline e$ and $\overline e'$ originate from two edges of $G$, $e$ and $e'$ respectively, which are incident respectively to two different preimages $x$ and $x'$ of $\overline x$. Since, by definition, $\ssf_x = \ssf_{x'}$ and $g(x) = g(x')$ for all $g \in M$, we can use exactly the same argument, \emph{mutadis mutandis}. \qedhere
			\end{enumerate}
		\end{proof}
		
		\begin{remark} \label{rem:edges_partition_glued}
			The above proof shows that the branches incident to some $x \in \Gamma$ that are glued by $\Red$ are exactly those belonging to the same set in the partition $\ssmathcalP_x$.
		\end{remark}
		
		\begin{prop} \label{prop:T_no_cycles}
			$T$ has no cycles, and therefore is a metric tree.
		\end{prop}
		
		\begin{proof}
			First, we claim that $T$ can have no oriented cycle. Indeed, Proposition~\ref{prop:equivalence_classes_finite} states that if $x$ and $y$ are the vertices of an edge oriented from $x$ to $y$, then $x < y$ (using the strict partial order induced on $T$). In an oriented cycle, we would get a strict inequality of the form $x < x$, which is absurd.
			
			To conclude, we show that $T$ cannot have any cycle. Indeed, a cycle endowed with an orientation of its edges, if it is not an oriented cycle, gives an in-degree two to at least one of its vertices, which is impossible by Lemma~\ref{lem:indegree_two}.
		\end{proof}
		
		We thus conclude that $T$ is an acyclic metric graph, that is, a metric tree.

		\begin{proof}[Proof of Theorem~\ref{thm:g1d}: $M = \Red^*(N)$]
			We show the assertion in Theorem~\ref{thm:g1d} that $M = \Red^*(N)$, where $N$ is associated to $\Red(\ssx_0) \in T$.
			
			(The inclusion $\subseteq$) Let $f \in M$. By definition of $\Red$, $f$ is constant on the equivalence classes for $\sssim_M$, so it can be written in the form $f = g \circ \Red$, with $g$ a function $T \rightarrow \R$. It is straightforward that $g$ is continuous and affine linear with integral slopes (see Proposition~\ref{prop:equivalence_classes_finite} which implies that $\Red$ is a local homeomorphism on its image). By definition of the metric on $T$, $g$ has slopes zero and one compatible with $\fS^\ee$, with $\fS$ the slope structure underlying $N$. To show that $g$ is compatible with $D$ and $\fS^\vv$, we simply look at the possible sets of slopes of $f$ around a point $y \in \Gamma$ in the cases (i) and (ii) explored in the proof of Proposition--Definition~\ref{prop_def:unique_g11_tree}.
			
			(The inclusion $\supseteq$) The other way around, let $g \in N$ and $f \coloneqq g \circ \Red$. We have to show that $f \in M$. Since $g$ belongs to $N$, it is of the form $g = \ssf_y = \ssf_{\Red(\ssx_0) \rightarrow y} \in \Rat((\Red(\ssx_0)), \fS)$ with $y \in T$. Let $w \in \Red^{-1}(y)$. All the functions $\ssf_{\ssx_0 \rightarrow w}$ for such a $w$ are equal thanks to Proposition~\ref{prop:equivalence_reduced_map}. The fact that $f$ belongs to $M$ then comes from the equality $f = \ssf_w$ which is implied by Theorem~\ref{thm:jumps_realized}.
		\end{proof}
		
		To finish the first part of the theorem, we need to show that the morphism $\Red$ from $\Gamma$ to $T$ can be resolved to a finite harmonic morphism of degree $d$ by a tropical modification of $\Gamma$. Given a point $\bar x \in T$ and $\nu \in \T_{\bar x}(T)$, we denote by $T_{\bar x, \nu}$ the metric subtree of $T$ consisting of the point $\bar x$ and all the points $z$ of $T$ which have the property that the unique path from $\bar x$ to $z$ in $T$ has tangent vector at $\bar x$ equal to $\nu$. For the language of harmonic morphisms and degrees of maps between metric graphs, we refer the reader to~\cite[Section~2.1]{amini2015lifting1}.
		
		\begin{proof}[Proof of Theorem~\ref{thm:g1d}: resolution to a finite harmonic morphism]
			Let $\bar x$ be a point of $T$, and consider a point $x$ of $\Gamma$ with $\Red(x) = \bar x$. We denote by $\Red_* \colon \T_x(\Gamma) \to \T_{\bar x}(T)$ the induced map on tangent vectors. For each unit tangent vector $\nu \in \T_{\bar x}(T)$, consider the degree of $\Red$ at $x$ above $\nu$, denoted by $\deg_\nu\Red(x)$ and defined as the sum
			\[
				\deg_\nu \Red(x) \coloneqq \sum_{\substack{\mu \in \T_x(\Gamma) \\ \Red_*(\mu) = \nu}} \slope_\mu(\Red)(x).
			\]
			Define the degree of $\Red$ at $x$, denoted by $\deg \Red(x)$, to be the maximum quantity $\deg_\nu \Red(x)$ for all $\nu \in \T_{\bar x}(T)$.
			
			We define the tropical modification $\wtGamma_0$ as follows. For each point $\bar x \in T$ and each point $x \in \Gamma$ as above, and each $\nu \in \T_{\bar x}(T)$ with $\deg_\nu \Red(x) < \deg \Red(x)$, we take $\deg \Red(x) - \deg_\nu \Red(x)$ copies of the subtree $T_{\bar x, \nu}$ of $T$, glue them to the point $x$ by identifying the point $\bar x$ of 	$T_{\bar x, \nu}$ in each copy with the point $x$. We then naturally extend $\Red$ to each copy of $T_{\bar x, \nu}$ by the identity map on $T_{\bar x, \nu}$. Since $\deg \Red(x) - \deg_\nu \Red(x)$ is non-zero only for finitely many pairs $(\bar x, \nu)$, we obtain a tropical modification $\wtGamma_0$ of $\Gamma$ and a map $\phi!_0 \colon \wtGamma_0 \to T$ which coincides with $\Red$ on the metric subgraph $\Gamma \hookrightarrow \wtGamma_0$. By construction, the map $\phi!_0$ is harmonic of degree $\deg \Red(x)$ at all points $x$ of $\Gamma$. It is also harmonic of degree one at all other points of $\wtGamma_0 \setminus \Gamma$. Note that $\phi!_0$ is of degree $\ssd_0$ at most $d$. In fact, in the gluing process above, in the construction of $\wtGamma_0$, when we add $\deg \Red(x) - \deg_\nu \Red(x)$ at $x$, the degree remains bounded by $d$. Indeed, for each point $\bar y$ in the subtree $T_{\bar x, \nu}$, and each $y \in \Gamma$ with $\Red(y) = \bar y$, the coefficient of $x$ in $D_y^M$ is at least $\deg \Red(x) - \deg_\nu \Red(x)$, using the explicit description of the variation of reduced divisors with respect to the base point provided in Section~\ref{subsec:behavior_reduced_divisors_with_respect_base_point}. This implies that the fiber $\Red^{-1}(\bar y)$, counted by multiplicity, is of size bounded $d - \deg \Red(x) + \deg_\nu \Red(x)$, so that after modification, the degree remains bounded by $d$. Treating the points one by one in the gluing process, and proceeding by induction, we obtain $\ssd_0 \leq d$.
			
			To conclude, let $\wtGamma$ be the tropical modification of $\wtGamma_0$ obtained by plugging $d - \ssd_0$ copies of $T$ at some point $y$ of $\wtGamma_0$, by identification of the point $\phi_0(y)$ in each copy with $y$. We extend $\phi!_0$ by the identity map on these copies to get the harmonic map $\phi \colon \wtGamma \to T$ of degree $d$, as required.
		\end{proof}
		
		Now we show the second part of Theorem~\ref{thm:g1d}: if $\psi \colon \wtGamma \rightarrow T$ is any finite harmonic morphism to a metric tree $T$ on which we put the semimodule $N$ using the point $\ssy_0 = \psi(\ssx_0)$, then $\mleft(\tau\mleft(\psi^{-1}(\ssy_0)\mright), \psi^*(N)\mright)$ is an effective $\g^1_d$ on $\Gamma$.
		
		\begin{proof}[Proof of Theorem~\ref{thm:g1d}: the converse]
			It will be enough to show that $\mleft(\psi^{-1}(\ssy_0), \psi^*(N)\mright)$ is an effective $\g^1_d$ on $\wtGamma$.

			We first have to define the pullback of the slope structure $\fS'$ of $N$ by $\psi$ to $\wtGamma$. Let $D$ be the divisor $\psi^{-1}(\ssy_0)$. Its degree is $d$, the degree of $\psi$. At some point $x \in \Gamma$, in the direction $\nu \in \T_x(\Gamma)$, the non-zero possible slope is defined to be the relative slope of $\psi$ in the direction $\nu$: this defines $\fS^\ee$. We now define $\fS^\vv$ around $x$ by saying that the jumps of $\ssrho_x$ are exactly the vector $\underline 0$, and the vectors having ones for all edges belonging to a certain complete set of edges identified by $\psi$, and zero on all other edges, which entirely defines $\ssrho_x$. We have defined the pair $(D, \fS)$.
			
			We now show that $\psi^*(N)$ is a semimodule included in $\Rat(D, \fS)$. Firstly, it is stable by the two tropical operations since $N$ is. Secondly, we show that $\psi^*(N) \subseteq \Rat(D, \fS)$. If $f$ is a function of $\psi^*(N)$, we can write it $f = g \circ \psi$ with $g \in N$. It is automatic by the construction of $\fS$ on $\wtGamma$ that $f$ is compatible with $\fS^\ee$ and $\fS^\vv$. The fact that $D + \div(f) \geq 0$ comes from the harmonicity of $\psi$.
			
			We now check property $(**)$ of Section~\ref{subsec:admissible_semimodules}. Let $x \in \Gamma$ and $E = (x)$. Let $y = \psi(x)$. Then, the function $\ssf_{\ssy_0 \rightarrow y} \circ \psi$ has the required properties (1) and (2). This function belongs to $\psi^*(N)$. Therefore, $(D, \psi^*(N))$ verifies property $(**)$ of Section~\ref{subsec:admissible_semimodules}. Finally, the finite generation of $\psi^*(N)$ and the fact that $\psi^*(N)$ has tropical rank one follow from the same properties for $N$. The same is true for effectivity. We have proved that $\mleft(D, \psi^*(N)\mright)$ is an effective $\g^1_d$.
		\end{proof}
		
		This finishes the classification of $\g^1_d$'s on metric graphs.


\section{Limit linear series on the skeleton of a Berkovich curve} \label{sec:geometric_tropicalization}
	
	In this section, we show that the tropicalization of linear series on $\curve$ gives combinatorial limit linear series.
	
	Let $\K$ be an algebraically closed field with a non-trivial non-Archimedean valuation $\val$ and let $\curve$ be a smooth proper curve over $\K$. We assume that $\K$ is complete with respect to $\val$. Let $\vR$, $\fm$, and $\k = \vR/\fm$ be the valuation ring, the maximal ideal of $\vR$, and the residue field, respectively. Denote by $\Lambda$ the value group of $\val$. We also denote by $|\cdot|$ the corresponding norm on $\K$, so that $\val(\cdot) = -\log|\cdot|$. Let $\K(\curve)$ be the function field of $\curve$.
	
	A semistable vertex set for $\curve^{\an}$ is a finite set of type 2 points $V$ in $\curve^{\an}$ such that the complement $\curve^{\an} \setminus V$ is a disjoint union of finitely many open annuli and infinitely many open disks. A semistable vertex set $V$ gives rise to a skeleton $\Gamma$ for $\curve^{\an}$, defined as the union in $\curve^{\an}$ of $V$ and the skeleta of the open annuli in $\curve^{\an} \setminus V$. The canonical metric on the skeleta of the open annuli gives the skeleton a metric graph structure, naturally embedded in $\curve^{\an}$. The underlying graph $G = (V, E)$ has vertex set $V$ and edge set $E$ in bijection with the set of open annuli in $\curve^{\an} \setminus V$. There is an edge between a pair of vertices $v$ and $u$ in $V$ for each open annulus whose closure contains the points $v$ and $u$. Moreover, the edge length function $\ell \colon E \to (0, +\infty)$ associates to each edge of $G$ the modulus of the corresponding annulus.
	
	Let $\Gamma$ be a metric graph skeleton of $\curve^{\an}$ with underlying graph $G = (V, E)$ and denote by $\tau \colon \curve^{\an} \to \Gamma$ the canonical retraction map. We call $\tau$ the tropicalization map. We get a tropicalization map $\tau_* \colon \Div(\curve) \to \Div_\Lambda(\Gamma)$ that sends a divisor $\varD = \sum_{x \in \curve(\K)} a_x (x)$ on $\curve$ to the divisor $\tau_*(\varD) = \sum_{x \in \curve(\K)} a_x (\tau(x))$.
	
	We denote by $\valuation_x \colon \K(\curve) \to \R \cup \{+\infty\}$ the valuation of a point $x \in \curve^{\an} \setminus X(\K)$ with $\valuation_x(f) = +\infty$ only if $f = 0$. The residue field of this valuation is denoted by $\k(x)$. We also denote by $|\cdot|_x = \exp(-\valuation_x)$ the corresponding norm.
	
	For each non-zero $\varf \in \K(\curve)$, we define the tropicalization of $\varf$, denoted $\Trop(\varf) \colon \Gamma \to \R$, as the map that sends each $x \in \Gamma \subseteq \curve^{\an} \setminus \curve(\K)$ to $\valuation_x(\varf)$. This induces a tropicalization map $\Trop \colon \K(\curve) \setminus \{0\} \to \Rat_\Lambda(\Gamma)$.
	
	Let $\varD$ be divisor of degree $d$ and rank at least $r$ on $\curve$, and let $(\mathcal O(\varD), \varH)$ be a $\g^r_d$ on $\curve$. We identify $\varH$ with a subspace of $\K(C)$ of dimension $r + 1$.
	
	Let $\Gamma$ be a skeleton of $\curve^{\mathrm{an}}$. We define
	\[ M \coloneqq \Trop(\varH) = \mleft\{\Trop(\varf) \st \varf \in \varH \setminus \{0\}\mright\}. \]
	
	\begin{thm} \label{thm:specialization_linear_series_grd}
		The pair $(D, M)$ is a refined $\Lambda$-rational $\g^r_d$ on $\Gamma$.
	\end{thm}
	
	We will call the linear equivalence class of the pair $(D, M)$ \emph{the combinatorial limit linear series on $\Gamma$ induced by $(\mathcal O(\varD), \varH)$}. It is easy to see that $(D, M)$ is effective provided that $\varH$ contains constants.
	
	\subsection{Reduction}	
		
		For each point $x$ of type 2 in $\curve^{\an}$, the extension $\k(x)/\k$ is of transcendence degree one. Let $\rmC_x$ be the smooth projective curve over $\k$ with function field $\k(x)$. For the point $x$, the valuation $\valuation_x$ has the same value group as $\val$. For each nonzero $ \varf \in \K(\curve)$, choosing $a \in \K$ with $\val(a) = \valuation_x(\varf)$, we get that $\ssa^{-1} \varf$ has valuation $\valuation_x(\ssa^{-1} \varf) = 0$, and therefore gives an element in the residue field $\k(x)$ that we denote by $\sstildef_x$. We call this the \emph{reduction of $\varf$ at $x$}, which is nonzero and defined only up to multiplication by a non-zero scalar in $\k$. For a vector subspace $\varH \subset \K(X)$ of dimension $r + 1$, denote by $\sstildeH_x \subset \k(x)$ the $\k$-vector subspace spanned by the reductions $\sstildef_x$ of elements $\varf \in \varH$ \cite[Section 4.4]{AB15}. By \cite[Lemma 4.3]{AB15}, $\sstildeH_x$ has dimension $r + 1$ over $\k$.
	
	\subsection{Slope structure and crude linear series coming from tropicalization} \label{subsec:linear_series_specialization}
		
		For a point $x$ in $\Gamma$ of type 2 in $\curve^{\an}$, each unit tangent direction $\nu \in \T_x(\Gamma)$ gives a point $\ssp_x^\nu \in \rmC_x(\k)$. By the slope formula~\cite{BPR}, for any non-zero $\varf \in \K(\curve)$, we have $\slope_\nu(\Trop(\varf)) = \ord_{\ssp_x^\nu}\mleft(\sstildef_x\mright)$. Moreover, as a consequence of the slope formula~\cite{AB15}, we get 
		\[ \tau_*(\div(\varf)) = \div(\Trop(\varf)). \]
		(Note that there is a sign difference between our definition of the divisor of a rational function and that of~\cite{AB15}.)
		
		If $\varH \subset \K(\curve)$ is a $\K$-vector subspace of dimension $r + 1$, for any unit tangent vector $\nu \in \T_x(\Gamma)$, we get a collection of integers $\slope_\nu(\Trop(\varf)) = \ord_{\ssp_x^\nu}\mleft(\sstildef_x\mright)$, $\varf \in \varH$. Since $\sstildeH_x$ has dimension $r + 1$, this collection has size $r + 1$. This means that the collection of slopes $\slope_{\nu}(h)$, for $h \in M = \Trop(\varH)$, has size $r + 1$. For each unit tangent vector $\nu$, we order the slopes $\slope_{\nu}(h)$, for $h \in M$, in the form $s^\nu_0 < s^\nu_1 < \cdots < s_r^{\nu}$. Define $S^\nu \coloneqq \ssub{\{\sss_i^ \nu\}}!_{i \in [r]}$. In addition, the collection of points $\ssp_x^\nu \in \ssC_x(\k)$ for $\nu \in \T_x(\Gamma)$ defines a geometric rank function $\ssrho_x$ associated to the corresponding filtrations on $\sstildeH_x$ as in Section~\ref{subsubsec:geometric_rank_functions}.
		We define $S^x$ as the set of jumps of $\ssrho_x$. We have the following theorem, which can be regarded as a refinement of~\cite[Theorem~5.9]{AB15}.
		
		\begin{prop}[Slope structures induced by tropicalization of rational functions] \label{thm:specialization_linear_series_divisor}
			Notation as above, let $\varH \subset \K(\curve)$ be a $\K$-vector space of dimension $r + 1$. Let $\Gamma$ be a skeleton of $\curve^{\an}$. There exists a semistable vertex set $V$ for $\curve$ such that $\Sigma(\curve, V) = \Gamma$, and such that the slopes of the tropicalizations of rational functions $\varf$ in $\varH$ along edges in $\Gamma$ define a slope structure of width $r$ on $\Gamma$.
		\end{prop}
		
		\begin{proof}
			We already defined $S^x$ and $S^\nu$ for type 2 points of $\Gamma \subseteq \curve^{\an}$ and $\nu \in \T_x(\Gamma)$. We show that the definitions can be extended to all points of $\Gamma$, and that the collection $\fS = \mleft\{S^x; S^\nu\mright\}_{x \in \Gamma, \, \nu \in \T_x(\Gamma)}$ is induced by a simple graph model of $\Gamma$ (or, equivalently, by a semistable vertex set of $\curve^{\an}$). To show this, let $x$ be a point of type 2, and let $\nu \in \T_x(\Gamma)$ be a tangent direction at $x$ in $\Gamma$. Let $\varf_0, \dots, \varf_r$ be a basis of $\varH$ such that the reductions $\ssf_{0} = \tilde \varf_0, \dots, f_{r} = \tilde \varf_r$ to $\k(x)$ yield the orders of vanishings $\sss_0^\nu, \dots, \sss_r^\nu$ at the point $\ssp_x^\nu$, respectively. By the slope formula, the slope of $\Trop(\varf_j)$ along $\nu$ at $x$ coincides with $\sss_j^\nu$. There thus exists a half segment $\ssI_x^\nu = [x, y^\nu]$ on the edge supporting the point $x$ and the tangent direction $\nu$ such that the slope of the function $\Trop(\varf_j)$ along $\nu$ at every point of $\ssI_x^\nu$ is $\sss_j^\nu$. Since these are $r + 1$ integers, they are all the possible slopes along $\nu$ in $\ssI_x^\nu$ of tropicalizations of functions in $\varH$. We extend the definition of $S^\bullet$ to every point in the interior of the segment $\ssI_x^\nu$ by taking these slopes and by declaring the rank function to be standard. Applying now the compactness of $\Gamma$, we deduce a finite covering of $\Gamma$ by segments $\ssI_x^\nu$, from which we deduce the statement in the proposition.
		\end{proof}
		
		The fact that $\tau(\div(\varf)) = \div(\Trop (\varf))$ then shows that $(D, \fS)$ has rank $r$, as defined in Definition~\ref{subsec:def_crude_linear_series}. This is analogous to the proof of the specialization theorem for metrized complexes in~\cite{AB15}, we thus omit the details. This shows that
		
		\begin{prop}
		 	The pair $(D, \fS)$ defined by tropicalization is a crude linear series on $\Gamma$ of degree $d$ and rank $r$.
		\end{prop}
		
		Note that for two linearly equivalent divisors $\varD \sim \varD'$ on $\curve$, and $\varH$ a subspace of the space of global sections of the corresponding line bundles $\mathcal O(\varD) = \mathcal O(\varD')$ of dimension $r + 1$, the two pairs $(D, \fS)$ and $(D', \fS')$ are linearly equivalent.
	
	\subsection{$M \subseteq \Rat_\Lambda(D)$ is a semimodule}
		
		We show that $M \subseteq \Rat_\Lambda(D)$ is a semimodule. Changing the vertex set if necessary, we can assume without loss of generality that the support of $D$ is included in the set of vertices $V$. Let $f$ and $g$ be two elements in $M$, and let $\lambda \in \Lambda$. By definition, we can write $f = \Trop(\varf)$ and $g = \Trop(\varg)$ with $\varf, \, \varg \in \varH \setminus \{0\}$. Let $h \coloneqq \min(f + \lambda, g)$ for $\lambda \in \Lambda$. We show the existence of $\varh \in \varH$ such that $\Trop(\varh) = h$. This will show that $M$ is a $\Lambda$-semimodule. Since $\lambda \in \Lambda = \val(\K \setminus \{0\})$, we can write $\lambda = \val(\alpha)$ with $\alpha \in \K$. We take $\varh \coloneqq \alpha \varf + \varg$. The goal is to show that replacing, if necessary, $\alpha$ with $\alpha \beta$ for $\beta \in \K^\times$ of valuation zero, we have $h = \Trop(\varh)$.
		
		Let $x \in \Gamma$ be such that $f(x) + \lambda \neq g(x)$. Then, by the enhanced non-Archimedean triangular inequality (that is, if $\val(a) \neq \val(b)$, then $\val(a + b) = \min(\val(a), \val(b))$), we get automatically that $\valuation_x(\varh) = \min(f(x) + \lambda, g(x)) = h(x)$.
		
		Let now $\Gamma_0 \coloneqq \{x \in \Gamma, \, f(x) + \lambda = g(x)\}$. Since $f$ and $g$ are piecewise affine linear, $\Gamma_0$ can be written as a union of finitely many segments of $\Gamma$ so that both $f$ and $g$ are affine on each of these segments. Let now $I$ be any of them, whose extremities we denote $x$ and $y$, on an edge $e$ of $\Gamma$. A segment of $\Gamma$ is a segment of an edge of $G$; it can be reduced to just a point. Let $\nu$ be the tangent direction in $\T_x(\Gamma)$ pointing toward $y$; if $x = y$, we choose an arbitrary $\nu \in \T_x(\Gamma)$. Consider the point $\ssp_x^\nu$ of $\ssC_x(\k)$ corresponding to $\nu$. By the slope formula and the fact that $f + \lambda = g$ on all $I$, we have
		\[ \ord_{\ssp_x^\nu}\mleft(\widetilde \varf\mright) = \slope_\nu(f) = \slope_\nu(g) = \ord_{\ssp_x^\nu}\mleft(\widetilde \varg\mright). \]
		Up to multiplying $\alpha$ by some element $\beta$ of $\K^\times$ of valuation zero, we can ensure that
		\[ \ord_{\ssp_x^\nu}\mleft(\widetilde{\alpha \varf + \varg}\mright) = \ord_{\ssp_x^\nu}\mleft(\widetilde \varf\mright) = \ord_{\ssp_x^\nu}\mleft(\widetilde \varg\mright) = \ord_{\ssp_x^\nu}\mleft(\widetilde \alpha \widetilde \varf + \widetilde \varg\mright). \]
		Note that only one value for the reduction $\widetilde \alpha$ in $\k$ is forbidden.
		
		Using again the slope formula for $\widetilde{\alpha \varf + \varg}$ yields $\ord_{\ssp_x^\nu}\mleft(\widetilde \alpha \widetilde \varf + \widetilde \varg\mright) = \mathrm{slope}_\nu(h)$, so that
		\[ \mathrm{slope}_\nu(h) = \mathrm{slope}_\nu(f) = \mathrm{slope}_\nu(g). \]
		
		We can do the same at the other extremity $y$ of $I$, and ensure that locally, starting at either extremity of $I$, $h$ has the same slope as $f$ and $g$. Since $f$ and $g$ are linear on $I$, $D$ has no support in the interior of $I$, and $h$ coincides with $f + \lambda$ and $g$ on the extremities of $I$, we must have $h = f + \lambda = g$ on the full interval $I$.
		
		Since this can be done for each of the finitely many segments composing $\Gamma_0$, forbidding at most one value for the reduction $\tilde \alpha \in \k$ each time, and since $\k$ is algebraically closed (and thus infinite), there is some $\alpha \in \K$ such that $\val(\alpha) = \lambda$ and
		\[ \Trop(\alpha \varf + \varg) = \min(f + \lambda, g) = h \]
		on all $\Gamma$. We have shown that $h \in M$, so $M$ is a semimodule of $\Rat_\Lambda(D, \fS)$.
	
	\subsection{Finite generation property}
			
		\begin{prop}\label{prop:finite_generation}		
			The semimodule $M = \Trop(\varH)$ is finitely generated.
		\end{prop}
		
		We suppose that the slope structure $\fS$ is defined on the model $G = (V, E)$ of $\Gamma$ associated to the semistable vertex $V$.
		
		Let $\ssA^1 = \Spec(\K[\varT])$ and let $\ssA^{1, \an}$ be its Berkovich analytification. Let $A(\rho)$ be the closed annulus in $\ssA^{1, \an}$ of center $0$ with outer radius one and inner radius $\rho \in (0, 1)$,
		\[ A(\rho) = \mleft\{x \in \ssA^{1, \an} \st \rho \leq |\varT|_x \leq 1 \mright\}. \]
		Let $R(\rho)$ be the ring of analytic functions on $A(\rho)$. An analytic function $f$ on $A(\rho)$ admits a formal power series expansion
		\[ f = \sum_{n \in \Z} a_n \varT^n \]
		with $\lim_{n \to \pm \infty} |a_n| s^n = 0$ for all $s \in [\rho, 1]$. The skeleton of $A(\rho)$ is a closed interval, which can be identified with $I \coloneqq [0, -\log \rho]$: each point $q$ in this interval corresponds to the norm $|\cdot|_{\zeta_q}(f) = \sup_{n \in \Z} |a_n| \exp(-q n) = \max_{n \in \Z} |a_n| \exp(-q n)$ on any analytic function $f$ as above. The tropicalization of an analytic function $f$ is the function $\Trop(f)$ on the interval $I$ given by
		\[
			\Trop(f)(q) = \min \{\val(a_n) + n q \st n \in \Z \} \qquad \forall q \in I.
		\]
		
		Each edge $e$ in $E$ is the skeleton of one of the annuli in the complement $\curve^{\an} \setminus V$, that we denote by $A(e)$. We have $A(e) \simeq A(\rho_e)$ with $\rho_e = \exp(-\ell_e)$, where $\ell_e$ is the length of $e$ in $\Gamma$. Using this identification, $e$ is identified in $\Gamma$ with the interval $[0, \ell_e]$.
		
		\begin{proof}[Proof of Proposition~\ref{prop:finite_generation}]
			By restriction, each element $\varf \in \varH$ gives rise to an analytic function on the annulus $A(e)$ that we denote by $\varf_e$. The tropicalization of $\varf$ restricted to $e$ coincides with the tropicalization of $\varf_e$, and has slopes among $s_0^e < \dots < s_r^e$. It follows that, taking the analytic development $\varf_e = \sum_{n \in \Z} a_n^e(\varf) \varT^n$ in $A(e)$, with $a_n^e(\varf) \in \K$, the tropicalization of $\varf_e$ is entirely defined by the truncation
			\[ \underline{\varf}_e \coloneqq \sum_{n = s^e_0}^{s^e_e} a_n^e(\varf) \varT^n. \]
			Consider the $\K$-linear map 
			\begin{align*} 
				\varphi \colon \varH \longrightarrow \prod_e \K^{[s^e_0, s^r_e]}, \qquad \qquad
				\varf \longmapsto \mleft(\mleft(a^e_n(\varf)\mright)_{n = s^e_0, \dots, s^r_e}\mright)_{e \in E}.
			\end{align*}
			We infer that the tropicalization $\Trop(\varf)$, for $\varf$ non-zero in $\varH$, is entirely determined by tropicalization 
			\[
				\Trop(\varphi(\varf)) \coloneqq \mleft(\mleft(\val(a^e_n(\varf))\mright)_{n = s^e_0, \dots, s^r_e}\mright)_{e \in E} \in \prod_e \mathbb{T}^{[s^e_0, s^r_e]}
			\]
			via the expression
			\[
				\Trop(\varf) \rest{e}(q) = \min \mleft\{\val(a_n^e(\varf)) + n q \st n = s^e_0, \dots, s^e_r \mright\} \qquad \forall q \in [0, \ell_e].
			\]
			Let $\varH' = \varphi(\varH)$, and $M' = \Trop(\varH') \coloneqq \mleft\{\Trop(\varphi(\varf)) \st \varf \in \varH \setminus\{0\}\mright\} \subset \prod_e \mathbb{T}^{[s^e_0, s^r_e]}$. From the above discussion we deduce the existence of a surjective morphism of semimodules $\varphi^{\trop} \colon M' \to M$. We conclude by observing that $M'$, being the tropicalization of a linear subspace of $\K^n$ for $n = \sum_e (s^e_r - s^e_0 + 1)$, is finitely generated, see~\cite[Theorem~3.1]{Spe08} and \cite[Chap.~III, Theorem~1.2.2]{Gaubert-thesis}. A generating set for $M'$ gives a generating set for $M$ via the surjection $\varphi^\trop$. The proposition follows.
		\end{proof}
	
	\subsection{Proof of Theorem~\ref{thm:specialization_linear_series_grd}}
		
		We have proved nearly all the properties to show that $(D, M)$ is a linear series, and by~\cite{JP14}, $M$ has tropical rank $r$.
		We only need to show that property $\threestars$ in Definition~\ref{def:refined_linear_series} is verified. Let $E$ be an effective $\Lambda$-rational divisor on $\Gamma$ of degree $s \leq r$. Choose $\varE$ to be any lift of $E$ to the curve $\curve$. Let $\ssH_{\varE} \subset H$ be a subspace of $H$ of rank $r - s$ such that for each $\varf \in \ssH_\varE$, we have $\div(\varf) + \varD - \varE \geq 0$. This space exists since $\varE$ imposes at most $s$ linear constraints on $H$. Let $\fS_E$ be the slope structure of width $r - s$ induced by the tropicalization of $\ssH_\varE$, which is a slope substructure of $\fS$. Let $\ssM_E \coloneqq \Trop(\ssH_\varE)$. Then, $(D, \ssM_E)$ is a linear series of rank $r - s$, $\ssM_E \subset M$, and the properties $(1)$ and $(2)$ of $\threestars$ are verified.
		
		At this point, the proof of Theorem~\ref{thm:specialization_linear_series_grd} is complete. \qed

\section{Examples and discussions}
	
	We provide a set of examples and some complementary results.
	
	\subsection{The dipole graph} \label{subsec:slope_structures_on_dipole_graph}
		
		We consider the dipole graph with four edges (of unit length to simplify the notation). The genus is $g = 3$ and the rank of the canonical divisor $K$ is $r = 2$. Denote by $u$ and $v$ the two vertices and by $\sse_1, \sse_2, \sse_3$ and $\sse_4$ the four edges of $\Gamma$ (see Figure~\ref{fig:dipole_slope_structure_divisor}).
		
		\subsubsection{A crude linear series of degree 4 and rank 2} \label{subsubsec:crude_linear_series_dipole_graph}
			
			For $i = 1, 2, 3, 4$, let $\sst_i \in \mleft(0, \, \frac 1 6\mright]$. For each choice of the $\sst_i$'s, we will construct a slope structure $\fS$ of width two.
			
			Let the $\sst_i$'s be fixed. For each $i$, we endow the edge $\sse_i$ with the slope sets $0 < 1 < 2$ on the interval $\mleft[0, \, \frac 1 2 - \sst_i\mright]$ in the direction of the arrows, $-1 < 0 < 1$ on the interval $\mleft[\frac 1 2 - \sst_i, \, \frac 1 2 + \sst_i\mright]$, and $-2 < -1 < 0$ on the interval $\mleft[\frac 1 2 + \sst_i, \, 1\mright]$.
			
			\begin{figure}[h!]
				\centering
				\begin{tikzpicture}[scale=0.8]
					\coordinate (A) at (0,0);
					\coordinate (B) at (4,0);
					
					\draw[postaction=decorate,decoration={markings,mark=at position 1/6 with {\arrow[red]{>}; \node[above left]{$(0, 1, 2)$};},mark=at position 1/3 with {\fill circle[radius=2.5pt];},mark=at position 1/2 with {\arrow[red]{>}; \node[above]{$(-1, 0, 1)$};},mark=at position 2/3 with {\fill circle[radius=2.5pt];},mark=at position 5/6 with {\arrow[red]{>}; \node[above right]{$(-2, -1, 0)$};}}]
					(A) .. controls (0,1.5) and (4,1.5) .. (B);
					
					\draw[postaction=decorate,decoration={markings,mark=at position 3/16 with {\arrow[red]{>}; \node[below left]{$(0, 1, 2)$};},mark=at position 3/8 with {\fill circle[radius=2.5pt];},mark=at position 1/2 with {\arrow[red]{>}; \node[below]{$(-1, 0, 1)$};},mark=at position 5/8 with {\fill circle[radius=2.5pt];},mark=at position 13/16 with {\arrow[red]{>}; \node[below right]{$(-2, -1, 0)$};}}]
					(A) .. controls (0,-1.5) and (4,-1.5) .. (B);
					
					\draw[postaction=decorate,decoration={markings,mark=at position 9/40 with {\arrow[red]{>};},mark=at position 9/20 with {\fill circle[radius=2.5pt];},mark=at position 1/2 with {\arrow[red]{>};},mark=at position 11/20 with {\fill circle[radius=2.5pt];},mark=at position 31/40 with {\arrow[red]{>};}}]
					(A) .. controls (1,0.5) and (3,0.5) .. (B);
					
					\draw[postaction=decorate,decoration={markings,mark=at position 5/24 with {\arrow[red]{>};},mark=at position 5/12 with {\fill circle[radius=2.5pt];},mark=at position 1/2 with {\arrow[red]{>};},mark=at position 7/12 with{\fill circle[radius=2.5pt];},mark=at position 19/24 with {\arrow[red]{>};}}]
					(A) .. controls (1,-0.5) and (3,-0.5) .. (B);
					
					\foreach \c in {A,B} {
						\fill (\c) circle (2.5pt);
					}
					\node[left] at (A) {$2$};
					\node[right] at (B) {$2$};
				\end{tikzpicture}
				\caption{A slope structure of width two on the dipole graph.}
				\label{fig:dipole_slope_structure_divisor}
			\end{figure}
			
			We now define suitable rank functions. Endow the eight points $\frac 1 2 - \sst_i, \frac 1 2 + \sst_i$ on the edges $\sse_i$ for $i = 1, 2, 3, 4$, respectively, with the rank function on $[2]^2$ defined by the array
			$\begin{pmatrix}
				0 & \color{blue} 0 & -1 \\
				1 & \color{blue} 1 & \color{blue} 0 \\
				\color{blue} 2 & 1 & 0
			\end{pmatrix}$
			(jumps in blue) and all the other points of $\Gamma$, including $u$ and $v$, with the standard rank function. This defines a slope structure $\fS$ of width two on $\Gamma$, and $(K, \fS)$ is a crude linear series of rank two on $\Gamma$.
		
		\subsubsection{An elementary $\g_4^2$ on the dipole graph} \label{subsubsec:elementary_g_2_4_dipole_graph}
			
			We will present an example of a combinatorial limit linear series of degree $4$ and rank $2$ on the dipole graph with four edges of unit length, the same graph as in Section~\ref{subsubsec:crude_linear_series_dipole_graph}. We consider the canonical divisor $K$ on $\Gamma$, and keep the notation of Section~\ref{subsubsec:crude_linear_series_dipole_graph}. The linear series will be based on a degenerate version of the slope structure defined in that section, essentially corresponding to the limit $t_i = 0$ for every $i = 1, 2, 3, 4$.
			
			More precisely, as shown in Figure~\ref{fig:dipole_elementary_clls}, we endow, for each $i$, the edge $\sse_i$ with the slope sets $0 < 1 < 2$ on the interval $\mleft[0, \frac 1 2 \mright]$, and $-2 < -1 < 0$ on the interval $\mleft[\frac 1 2, 1 \mright]$, in the direction of the arrows. We moreover endow every point of $\Gamma$, including $u$, $v$ and the middle point $\ssm_i$ of each edge $\sse_i$, with the standard rank function. This defines a slope structure $\fS$ of width two.
			
			\begin{figure}[h!]
				\centering
				\begin{tikzpicture}[scale=0.8]
					\coordinate (A) at (0,0);
					\coordinate (B) at (4,0);
					\coordinate (C) at (2,1.125);
					\coordinate (D) at (2,-1.125);
					
					\draw[postaction=decorate,decoration={markings,mark=at position 1/4 with {\arrow[red]{>}; \node[above left]{$(0, 1, 2)$};},mark=at position 1/2 with {\fill circle[radius=2.5pt];},mark=at position 3/4 with{\arrow[red]{>}; \node[above right]{$(-2, -1, 0)$};}}]
					(A) .. controls (0,1.5) and (4,1.5) .. (B);
					
					\draw[postaction=decorate,decoration={markings,mark=at position 1/4 with {\arrow[red]{>}; \node[below left]{$(0, 1, 2)$};},mark=at position 1/2 with {\fill circle[radius=2.5pt];},mark=at position 3/4 with{\arrow[red]{>}; \node[below right]{$(-2, -1, 0)$};}}]
					(A) .. controls (0,-1.5) and (4,-1.5) .. (B);
					
					\draw[postaction=decorate,decoration={markings,mark=at position 1/4 with {\arrow[red]{>};},mark=at position 1/2 with {\fill circle[radius=2.5pt];},mark=at position 3/4 with {\arrow[red]{>};}}]
					(A) .. controls (1,0.5) and (3,0.5) .. (B);
					
					\draw[postaction=decorate,decoration={markings,mark=at position 1/4 with {\arrow[red]{>};},mark=at position 1/2 with{\fill circle[radius=2.5pt];},mark=at position 3/4 with {\arrow[red]{>};}}]
					(A) .. controls (1,-0.5) and (3,-0.5) .. (B);
					
					\foreach \c in {A,B} {
						\fill (\c) circle (2.5pt);
					}
					\node[left] at (A) {$2$};
					\node[right] at (B) {$2$};
				\end{tikzpicture}
				\caption{A slope structure of width two on the dipole graph.}
				\label{fig:dipole_elementary_clls}
			\end{figure}
			
			Denote by $\Rat^{\mathrm{Sym}}(K)$ (resp., $\Rat^{\mathrm{Sym}}(K, \fS)$) the set of functions $f \in \Rat(K)$ (resp., $f \in \Rat(K, \fS)$) whose restriction to each $\sse_i$ is symmetric with respect to the middle point $\ssm_i$. It is not difficult to see that
			\[ \Rat^{\mathrm{Sym}}(K, \fS) = \Rat^{\mathrm{Sym}}(K) \subsetneq \Rat(K, \fS) \subsetneq \Rat(K), \]
			where both inclusions are moreover closed. The linear series will be defined as
			\[ M \coloneqq \Rat^{\mathrm{Sym}}(K, \fS) = \Rat^{\mathrm{Sym}}(K). \]
			
			Then, $M$ is a sub-semimodule of $\Rat(K, \fS)$ of rank two, which here also implies the same for $\Rat(K, \fS)$ and $\Rat(K)$. The symmetries of the functions $f \in M$ then imply that $M$ can be viewed as a \emph{complete} linear series: it is isomorphic to the linear series $\Rat(2 \, (u))$ on the tree $T$ obtained as the left half of $\Gamma$, see Figure~\ref{fig:dipole_elementary_clls_symmetry_tree} below.
			
			\begin{figure}[h!]
				\centering
				\begin{tikzpicture}[scale=0.5]
					\coordinate (A) at (0,0);
					\coordinate (B) at (2,1.125);
					\coordinate (C) at (2,-1.125);
					\coordinate (D) at (2,0.375);
					\coordinate (E) at (2,-0.375);
					
					\draw (A) .. controls (0,0.82) and (1.15,1.125) .. (B);
					\draw (A) .. controls (0,-0.82) and (1.15,-1.125) .. (C);
					\draw (A) .. controls (0.5,0.25) and (1.22,0.375) .. (D);
					\draw (A) .. controls (0.5,-0.25) and (1.22,-0.375) .. (E);
					
					\foreach \c in {A,B,C,D,E} {
						\fill (\c) circle (2.5pt);
					}
					\node[left] at (A) {$2$};
				\end{tikzpicture}
				\caption{The tree $T$ on which the complete linear series $\Rat(2 \, (u))$ is defined.}
				\label{fig:dipole_elementary_clls_symmetry_tree}
			\end{figure}
			
			Now that we have essentially reduced to the case of a complete linear series, we use~\cite[Theorem 6]{haase2012linear} and \cite[Corollary~9]{haase2012linear}, which show that $M$ is generated by its finitely many extremal points, and provide a characterization of the extremal points which enables to enumerate them all.
			Using this, $M$ is seen to be generated by the functions defined as follows. For every $i = 1, 2, 3, 4$, let $\ssf_i$ be the function on $\Gamma$ whose graph on the edge $\sse_i$ is
			\begin{center}
				\begin{tikzpicture}[scale=0.5]
					\draw[->] (0,0) -- (2.6,0);
					\draw[->] (0,0) -- (0,2.6);
					\draw (0,0) node[below]{$u$} node{$\bullet$};
					\draw (2,0) node[below]{$v$} node{$\bullet$};
					\draw (0,0) -- (1,2) -- (2,0);
				\end{tikzpicture}
			\end{center}
			with slopes $2$ and $-2$, and which is identically zero on all other edges. For every choice of indices $1 \leq i < j \leq 4$, let $\ssg_{i,j}$ be the function on $\Gamma$ whose graph on the edges $\sse_i$ and $\sse_j$ is
			\begin{center}
				\begin{tikzpicture}[scale=0.5]
					\draw[->] (0,0) -- (2.6,0);
					\draw[->] (0,0) -- (0,1.6);
					\draw (0,0) node[below]{$u$} node{$\bullet$};
					\draw (2,0) node[below]{$v$} node{$\bullet$};
					\draw (0,0) -- (1,1) -- (2,0);
				\end{tikzpicture}
			\end{center}
			with slopes $1$ and $-1$, and which is identically zero on both other edges. Then $M$ is generated by the functions $\ssf_i$ and $\ssg_{i,j}$, consisting of $10$ generators.
			Finally, a case-by-case analysis on the generators shows that $M$ has tropical rank two, concluding the proof that $M$ is a $\g_4^2$.
		
		\subsubsection{A more involved $\g_4^2$ on the dipole graph} \label{subsubsec:involved_g_2_4_dipole_graph}
			
			Keeping the same dipole graph as in the preceding two sections, we now present another example of a combinatorial limit linear series of degree $4$ and rank $2$. Unlike in Section~\ref{subsubsec:elementary_g_2_4_dipole_graph}, we here specialize the slope structure defined in Section~\ref{subsubsec:crude_linear_series_dipole_graph} to the choice of parameters $\sst_i = \frac 1 6$, $i = 1, 2, 3, 4$.
			
			Unlike the linear series in Section~\ref{subsubsec:elementary_g_2_4_dipole_graph}, symmetric functions will not be sufficient to get a $\g_4^2$ compatible with $\fS$, but we will constrain the functions in another way. Let
			$M \subsetneq \Rat(K, \fS)$
			be the subset of functions $f \in \Rat(K, \fS)$ which have slope $2$ along \emph{at most one tangent} vector based at $u$ or $v$. Equivalently, a function $f \in \Rat(K, \fS)$ belongs to $M$ if, and only if, on every edge $\sse_i$, if $f$ has slope $2$ close to $u$ on $\sse_i$, then its slope close to $v$ on $\sse_i$ is less than $2$ (in fact, it consequently has to be $1$).
			
			It is easily verified that $M$ is a sub-semimodule of $\Rat(K, \fS)$. We will show that it is a linear series of degree $4$ and rank $2$. Firstly, $M$ can be shown to have rank $2$; the key observation here is that functions of $\Rat(K, \fS)$ which take slope $2$ on both endpoints of some edge are not needed.
			
			To prove that $M$ is finitely generated, we then divide the functions $f \in M$ into two categories, and find a finite number of generators for each category, following the kind of strategy deployed in Section~\ref{subsubsec:elementary_g_2_4_dipole_graph}.	Namely,	we consider, on the one hand, the functions which are non-constant on exactly two distinct edges. Those functions cannot use the slope $2$. Like in Section~\ref{subsubsec:elementary_g_2_4_dipole_graph}, it is easy to see that they are consequently generated by the functions $\ssg_{i,j}$, for indices $1 \leq i < j \leq 4$, whose graph on the edges $\sse_i$ and $\sse_j$ is
			\begin{center}
				\begin{tikzpicture}[scale=0.5]
					\draw[->] (0,0) -- (2.6,0);
					\draw[->] (0,0) -- (0,1.6);
					\draw (0,0) node[below]{$u$} node{$\bullet$};
					\draw (2,0) node[below]{$v$} node{$\bullet$};
					\draw (0,0) -- (1,1) -- (2,0);
				\end{tikzpicture}
			\end{center}
			with slopes $1$ and $-1$, and which is identically zero on the other edges. And, on the other hand, we look at the sub-semimodule $M' \subseteq M$ of functions which are non-constant on at most one edge. For example, denote by $\ssM'_1$ the sub-semimodule of $M'$ consisting of the functions which are constant on all edges $\sse_i$, $i = 2, 3, 4$. Since all functions in $\Rat(K)$ have the same value on $u$ and $v$, and thanks to the slope constraint defining $M$, it follows that $\ssM'_1$ is isomorphic to the \emph{complete} linear series $\Rat(3 \, (u)) = \Rat(3 \, (u), \fS)$ on the cycle obtained by identifying the endpoints of $\sse_1$ and deleting the other edges, see Figure~\ref{fig:dipole_involved_clls_cycle}.
			
			\begin{figure}[h!]
				\centering
				\begin{tikzpicture}
					\coordinate (A) at (0,-0.5);
					\coordinate (B) at ($(30:0.5)$);
					\coordinate (C) at ($(150:0.5)$);
					
					\draw[postaction=decorate,decoration={markings,mark=at position 1/2 with {\arrow[red]{>}; \node[right]{$(0, 1, 2)$};}}]
					(A) arc (-90:30:0.5);
					
					\draw[postaction=decorate,decoration={markings,mark=at position 1/2 with {\arrow[red]{>}; \node[left]{$(0, 1, 2)$};}}]
					(A) arc (270:150:0.5);
					
					\draw[postaction=decorate,decoration={markings,mark=at position 1/2 with {\arrow[red]{>}; \node[above]{$(-1, 0, 1)$};}}]
					(B) arc (30:150:0.5);
						
					\foreach \c in {A,B,C} {
						\fill (\c) circle (2.5pt);
					}
					\node[below=2pt] at (A) {$3$};
					\node[above=2pt] at (A) {$u$};
				\end{tikzpicture}
				\caption{A cycle of length one with slope structure inherited from $\fS$.}
				\label{fig:dipole_involved_clls_cycle}
			\end{figure}
			
			As in Section~\ref{subsubsec:elementary_g_2_4_dipole_graph}, reducing to a complete linear series provides a set of generators for $\ssM_1'$, obtained as the extremal points. It follows that $\ssM'_1$ is generated by the function $\ssh_1$ whose graph on the edge $\sse_1$ is
			\begin{center}
				\begin{tikzpicture}[scale=0.5]
					\draw[->] (0,0) -- (2.6,0);
					\draw[->] (0,0) -- (0,1.933);
					\draw (0,0) node[below]{$u$} node{$\bullet$};
					\draw (2,0) node[below]{$v$} node{$\bullet$};
					\draw (0,0) -- (2/3,4/3) -- (2,0);
				\end{tikzpicture}
			\end{center}
			with slopes $2$ on $\mleft[0, \frac 1 3 \mright]$ and $-1$ on $\mleft[ \frac 1 3, 1\mright]$, and which is identically zero on both other edges, together with the function $\ssh'_1$ obtained from $\ssh_1$ by symmetry with respect to the midpoint of each edge. Likewise, for $i = 2, 3, 4$, we define functions $\ssh_i$ and $\ssh'_i$. Then the set of all these functions generates $M'$.
			
			All in all, it follows that
			\[ M = \mleft\langle \ssg_{i,j}, \; \ssh_k \st 1 \leq i < j \leq 4, \; 1 \leq k \leq 4 \mright\rangle, \]
			which amounts to $14$ generators. An elementary check on these generators shows that $M$ has tropical rank two and therefore is a $\g_4^2$.
	
	\subsection{Realizability and genera} \label{subsec:realizability_question}
		
		In the context of tropicalization of linear series, the metric graph $\Gamma$ can be enriched with a supplementary numerical data consisting of a genus function $\g \colon V \to \Z_{\geq 0}$ associating to each vertex $v$ the genus of the curve $C_v$; the pair $(\Gamma, \g)$ is called an \emph{augmented metric graph}. The study of linear series on metric graphs is not dependent on the data of a genus function. However, considering genus functions becomes important in making connection to geometry when asking realizability questions. Indeed, for a slope structure to be realizable, that is, to be the tropicalization of a linear series on a curve defined over a non-Archimedean field, it is necessary that the divisor defined on $\Gamma$ by
		\begin{equation} \label{eq:clls_weierstrass_divisor_clls_discussions}
			\mu^\alg_\W(x) \coloneqq (r + 1) \, K(x) + \frac{r(r + 1)}{2} \, (\ssd_x + \g(x) - 2) - \sum_{\nu \in \T_x(\Gamma)} \, \sum_{j = 0}^r s^\nu_j
		\end{equation}
		be effective (here, $\ssd_x$ is the valence of $x$ and $\g(x)$ denotes the genus of the curve $C_x$ corresponding to $x$). See Section~\cite[Section 5]{AGR23}, more specifically Proposition~5.12 in \emph{loc.~cit.}. The realizability of a linear series on an augmented graph will depend on the genus function $\g$. On the one hand, the general form of Equation~\eqref{eq:clls_weierstrass_divisor_clls_discussions} puts constraints on $\g$ for the linear series to be realizable; and, on the other hand, the property of each rank function $\rho_v$ being geometric will depend on the genus $\g(v)$ of $C_v$.
		
		To give an example, consider the metric graph $\Gamma$ below with arbitrary edge lengths, and the divisor $K$ with coefficient one at the trivalent vertices. Consider a family of slope structures on $\Gamma$, as follows.
		
		First, for each bridge edge oriented outwards (towards the adjacent circle), allow slopes $-1 < 1 < 3$. Divide each circle into three equal parts, in a way compatible with the position of the attachment points. On each circle, on the two edges adjacent to the attachment points, allow slopes $0 < 1 < 2$ away from the attachment point, and on the remaining edge, allow slopes $-1 < 0 < 1$. We endow the vertices of $\Gamma$ with the following rank functions. The central vertex gets the standard rank function. The three attachment points are endowed with the rank function on $[2]^3$ whose restrictions to $[2]^2 \times \{0\}$, $[2]^2 \times \{1\}$ and $[2]^2 \times \{2\}$ are given respectively by the three matrices
		\[ \begin{pmatrix}
			0 & \color{blue} 0 & -1 \\
			1 & \color{blue} 1 & \color{blue} 0 \\
			\color{blue} 2 & 1 & 0
		\end{pmatrix}, \;
		\begin{pmatrix}
			-1 & -1 & -1 \\
			0 & \color{blue} 0 & -1 \\
			\color{blue} 1 & 0 & -1
		\end{pmatrix} \; \text{and} \;
		\begin{pmatrix}
			-1 & -1 & -1 \\
			-1 & -1 & -1 \\
			\color{blue} 0 & -1 & -1
		\end{pmatrix}. \]
		Jumps are depicted in blue and the third coordinate corresponds to the edge connecting the attachment point to the central vertex. Finally, the six other vertices of $\Gamma$ are endowed with the rank function on $[2]^2$ defined by the first of the three arrays above.
		The standard rank function is imposed on all other points of $\Gamma$. This defines a slope structure of width two $\fS$ on $\Gamma$, and $M \coloneqq \Rat(K, \fS)$ is an admissible semimodule.
		
		\begin{figure}[h!]
			\centering
			\begin{tikzpicture}[scale=0.8]
				\coordinate (A) at (0,0);
				\coordinate (B) at (0,1);
				\coordinate (C) at (-0.87,-0.5);
				\coordinate (D) at (0.87,-0.5);
				\coordinate (Bc) at (0,1.5);
				\coordinate (Cc) at (-1.3,-0.75);
				\coordinate (Dc) at (1.3,-0.75);
				\coordinate (Bf) at ($(Bc) + (30:0.5)$);
				\coordinate (Bg) at ($(Bc) + (150:0.5)$);
				\coordinate (Cf) at ($(Cc) + (150:0.5)$);
				\coordinate (Cg) at ($(Cc) + (-90:0.5)$);
				\coordinate (Df) at ($(Dc) + (270:0.5)$);
				\coordinate (Dg) at ($(Dc) + (30:0.5)$);
				
				\draw[postaction=decorate,decoration={markings,mark=at position 0.5 with{\arrow[red]{>}},mark=at position 0.5 with {\node[right] {$(0, 1, 2)$};}}]
				(B) arc (-90:30:0.5);
				
				\draw[postaction=decorate,decoration={markings,mark=at position 0.5 with{\arrow[red]{>}},mark=at position 0.5 with {\node[left] {$(0, 1, 2)$};}}]
				(B) arc (270:150:0.5);
				
				\draw[postaction=decorate,decoration={markings,mark=at position 0.5 with{\arrow[red]{>}},mark=at position 0.5 with {\node[above] {$(-1, 0, 1)$};}}]
				(Bf) arc (30:150:0.5);
				
				\begin{scope}[decoration={markings,mark=at position 0.5 with {\arrow[red]{>}}}] 
					\draw[postaction={decorate}] (A) -- node[label=right:{$(-1,1,3)$}] {} (B);
					\draw[postaction={decorate}] (A) -- (C);
					\draw[postaction={decorate}] (A) -- (D);
					\draw[postaction={decorate}] (C) arc (30:150:0.5);
					\draw[postaction={decorate}] (C) arc (30:-90:0.5);
					\draw[postaction={decorate}] (Cf) arc (150:270:0.5);
					\draw[postaction={decorate}] (D) arc (150:30:0.5);
					\draw[postaction={decorate}] (D) arc (150:270:0.5);
					\draw[postaction={decorate}] (Df) arc (-90:30:0.5);
				\end{scope}
					
				\foreach \c in {A,B,C,D,Bf,Bg,Cf,Cg,Df,Dg} {
					\fill (\c) circle (2.5pt);
				}
				\node[below=2pt] at (A) {$1$};
				\node[above=2pt] at (B) {$1$};
				\node[below right] at (C) {$1$};
				\node[below left] at (D) {$1$};
			\end{tikzpicture}
			\caption{A slope structure of width two on the three-cycle graph.}
			\label{fig:three_cycles}
		\end{figure}
		
		\medskip
		
		Now consider the example where, instead of allowing the slopes $-1 < 1 < 3$ on the three central edges, we allow the slopes $-1 < 1 < 3$ on (possibly trivial) intervals incident to the central vertex on these edges, and the slopes $-1 < 1 < 2$ on the rest of the edges, with the same choice of rank functions. Choosing the lengths of these intervals independently yields a three-parameter family of slope structures on $\Gamma$. 
		
		Among these possibilities, with genus function 0, only the one where all three intervals carrying the slopes $-1 < 1 < 3$ are trivial can be realizable.

\bibliographystyle{alpha}
\bibliography{Bibliography}

\end{document}